\documentclass[11pt,reqno]{amsart}
\usepackage{fixltx2e}                    
\usepackage{mathpazo}  
\usepackage[utf8]{inputenc}            
\usepackage{amsmath}                     
\usepackage{amssymb, latexsym, stmaryrd, amsthm, dsfont, amsfonts, amsbsy,amsthm, amsmath, mathrsfs}            
\usepackage{mathtools}                   
\usepackage{bm}                          
\usepackage{enumerate}                   
\usepackage{verbatim}                    
\usepackage{url}   
\usepackage{lscape}                      
\usepackage{soul}

\usepackage[margin=1.3in]{geometry}

\usepackage{microtype}  
\usepackage[all, knot]{xy}
        \xyoption{arc} 
        \xyoption{web}                 
\makeatletter                            
\def\MT@register@subst@font{\MT@exp@one@n\MT@in@clist\font@name\MT@font@list
 \ifMT@inlist@\else\xdef\MT@font@list{\MT@font@list\font@name,}\fi}
\makeatother

\usepackage[pdftex,bookmarks,bookmarksnumbered,linktocpage,   
         colorlinks,linkcolor=blue,citecolor=blue]{hyperref}

\newcommand{\bit}{\begin{itemize}}    
\newcommand{\eit}{\end{itemize}}
\newcommand{\ben}{\begin{enumerate}}
\newcommand{\een}{\end{enumerate}}
\newcommand{\benormal}{\ben[\normalfont 1.]}   

\let\enormal\een
\newcommand{\benroman}{\ben[\normalfont (i)]}  
\let\eroman\een

\newcommand{\bde}{\begin{description}}
\newcommand{\ede}{\end{description}}

\newcommand{\?}{\ensuremath{\mkern0.4\thinmuskip}}   

\let\leq=\leqslant
\let\geq=\geqslant
\let\Box=\square                            

\let\epsilon=\varepsilon
\let\Lambda\varLambda
\let\Gamma\varGamma
\let\Delta\varDelta
\let\Lambda\varLambda
\let\Omega\varOmega
\let\Theta\varTheta
\let\Xi\varXi
\let\Pi\varPi
\let\Sigma\varSigma

\let\class=\mathsf                              
\let\oper=\mathbb                               

\bmdefine{\A}{A}                                
\bmdefine{\2}{2}
\bmdefine{\B}{B}
\bmdefine{\D}{D}
\bmdefine{\M}{M}                                
\bmdefine{\LLL}{L}                              
\bmdefine{\Fm}{Fm}                              
\bmdefine{\zerou}{[0{,}1]}  
\bmdefine{\T}{T}                                

\newcommand{\LL}{\mathscr{L}}                   

\newcommand{\HHSinv}{{\oper{H}_{\textsc{s}}^{-1}}}

\newcommand{\Var}{\mathnormal{V\mkern-.8\thinmuskip ar}} 
      
\newcommand{\Alg}{\class{Alg}}

\bmdefine{\boldstar}{\mathchoice{\textstyle*}{\textstyle*}{\textstyle*}{\scriptstyle*}}

\bmdefine{\btau}{\tau}                                  
\bmdefine{\brho}{\rho}                                  

\newcommand{\sineq}{\mathrel{\dashv\mkern1.5mu\vdash}}  

\bmdefine{\leibniz}{\Omega}        

\bmdefine{\frege}{\Lambda}         

\makeatletter
\newcommand{\tarskidsp}{\mathord%
   {\m@th\raisebox{0pt}[0pt][0pt]{$\stackrel%
   {\raisebox{-2.7pt}[0ex][0pt]{$\displaystyle \,\?\thicksim$}}%
   {\displaystyle\leibniz}$}}}
\newcommand{\tarskitxt}{\mathord%
   {\m@th\raisebox{0pt}[0pt][0pt]{$\stackrel%
   {\raisebox{-2.7pt}[0ex][0pt]{$\,\?\thicksim$}}{\displaystyle\leibniz}$}}}
\newcommand{\tarskiscr}{\mathord%
   {{\m@th\raisebox{0pt}[0pt][0pt]{$\stackrel%
   {\raisebox{-2.4pt}[0ex][0pt]{$\scriptstyle \,\?\thicksim$}}%
   {\scriptstyle\leibniz}$}}}}
\newcommand{\tarskiscrscr}{\mathord%
   {{\m@th\raisebox{0pt}[0pt][0pt]{$\stackrel%
   {\raisebox{-2pt}[0ex][0pt]{$\scriptscriptstyle \,\?\thicksim$}}%
   {\scriptscriptstyle\leibniz}$}}}}
\newcommand{\tarski}{\@ifnextchar ^ %
   {\mathchoice{\tarskidsp\kern-.07em}{\tarskitxt\kern-.07em}%
   {\tarskiscr\kern-.07em}{\tarskiscrscr\kern-.07em}}%
   {\mathchoice{\tarskidsp}{\tarskitxt}{\tarskiscr}{\tarskiscrscr}}}
\makeatother

\theoremstyle{theorem}
\newtheorem{Theorem}{Theorem}[section]
\newtheorem{Proposition}[Theorem]{Proposition}
\newtheorem{Lemma}[Theorem]{Lemma}
\newtheorem{Corollary}[Theorem]{Corollary}
\newtheorem{Claim}[Theorem]{Claim}
\newtheorem{Notation}[Theorem]{Convention}

\theoremstyle{definition}
\newtheorem{law}[Theorem]{Definition}
\newtheorem{exa}[Theorem]{Example}

\theoremstyle{remark}

\newtheorem{problem}{\bf Problem}
\newtheorem{Remark}[Theorem]{Remark}


\begin{document}
\title[On equational completeness theorems]{On equational completeness theorems}


\author{T. Moraschini}
\address{Department of Philosophy, University of Barcelona, Carrer de Montalegre $6$, $08001$, Barcelona, Spain}
\email{tommaso.moraschini@ub.edu}
\date{\today}

\begin{abstract}
A logic is said to admit an equational completeness theorem when it can be interpreted into the equational consequence relative to some class of algebras. We characterize logics admitting an equational completeness theorem that are either locally tabular or have some tautology.\ In particular, it is shown that a protoalgebraic logic admits an equational completeness theorem precisely when its has two distinct logically equivalent formulas.\ While the problem of determining whether a logic admits an equational completeness theorem is shown to be decidable both for logics presented by a finite set of finite matrices and for locally tabular logics presented by a finite Hilbert calculus, it becomes undecidable for arbitrary logics presented by finite Hilbert calculi.
\end{abstract}

\maketitle

\section{Introduction}

A propositional logic $\vdash$ admits an \textit{equational (soundness and) completeness theorem} if there are a set of equations $\btau(x)$ and a class of algebras $\class{K}$ such that for every set of formulas $\Gamma \cup \{ \varphi \}$,
\begin{align*}
\Gamma \vdash \varphi \Longleftrightarrow& \text{ for every }\A \in \class{K} \text{ and } \vec{a} \in A,\\
&\text{ if }\A \vDash \btau(\gamma^{\A}(\vec{a})) \text{ for all }\gamma \in \Gamma, \text{ then }\A \vDash \btau(\varphi^{\A}(\vec{a})).
\end{align*}
In this case, $\class{K}$ is said to be an \textit{algebraic semantics} for $\vdash$ (or a $\btau$\textit{-algebraic semantics} for $\vdash$) \cite{BP89}. Accordingly, a logic admits an equational completeness theorem precisely when it has an algebraic semantics.\footnote{For related studies of completeness theorems involving equations, see for instance  \cite{BP89,CzJa00,He93a,TMo15,Ra06a}.}

For instance, in view of the well-known equational completeness theorem of classical propositional logic ${\bf CPC}$ with respect to the variety of Boolean algebras, stating that for every set of formulas $\Gamma \cup \{ \varphi\}$,
\begin{align*}
\Gamma \vdash_{{\bf CPC}} \varphi \Longleftrightarrow& \text{ for every Boolean algebra }\A \text{ and } \vec{a} \in A,\\
&\text{ if }\A \vDash \gamma^{\A}(\vec{a}) \thickapprox 1 \text{ for all }\gamma \in \Gamma, \text{ then }\A \vDash \varphi^{\A}(\vec{a}) \thickapprox 1,
\end{align*}
the variety of Boolean algebras is an algebraic semantic for ${\bf CPC}$.




Despite the apparent simplicity of the concept, intrinsic characterizations of logics with an algebraic semantics have proved elusive, partly because equational completeness theorem can take \textit{nonstandard} forms. For instance, by Glivenko's theorem \cite{Gli29}, 
for every set of formulas $\Gamma \cup \{ \varphi \}$,
\[
\Gamma \vdash_{{\bf CPC}} \varphi \Longleftrightarrow\text{ } \{ \lnot \lnot \gamma \colon \gamma \in \Gamma \} \vdash_{{\bf IPC}} \lnot \lnot \varphi,
\]
where ${\bf IPC}$ stands for intuitionistic propositional logic. Since Heyting algebras form an $\{ x \thickapprox 1 \}$-algebraic semantics for ${\bf IPC}$, one obtains
\begin{align*}
\Gamma \vdash_{{\bf CPC}} \varphi \Longleftrightarrow&\text{ for every Heyting algebra }\A \text{ and } \vec{a} \in A,\\
&\text{ if }\A \vDash \lnot \lnot \gamma^{\A}(\vec{a}) \thickapprox 1 \text{ for all }\gamma \in \Gamma, \text{ then }\A \vDash \lnot \lnot \varphi^{\A}(\vec{a})\thickapprox 1.
\end{align*}
Consequently, the variety of Heyting algebras is also an algebraic semantics for ${\bf CPC}$, although certainly not the intended one \cite[Prop.\ 2.6]{BlRe03}. 

As it happens, not only there is no easy escape from nonstandard equational completeness theorems, but, sometimes, these are the sole possible ones. It might be convenient to illustrate this point with a simple example, namely the $\langle \land, \lor \rangle$-fragment ${\bf CPC}_{\land\lor}$ of ${\bf CPC}$. As we proceed to explain, this fragment lacks any standard equational completeness theorem, i.e., one with respect to the variety of distributive lattices $\class{DL}$, but admits a nonstandard one.


\begin{proof}[Proof sketch.]
Suppose, with a view to contradiction, that ${\bf CPC}_{\land\lor}$ admits a standard equational completeness theorem. Then there exists a set of equations $\btau(x)$ such that $\class{DL}$ is a $\btau$-algebraic semantics for ${\bf CPC}_{\land\lor}$. As every equation in which only a single variable occurs is valid in $\class{DL}$, we obtain $\class{DL} \vDash \btau(x)$. Together with the assumption that $\class{DL}$ is a $\btau$-algebraic semantics for ${\bf CPC}_{\land\lor}$, this implies $\emptyset \vdash_{{\bf CPC}_{\land\lor}} x$, a contradiction. Consequently, ${\bf CPC}_{\land \lor}$ lacks any standard equational completeness theorem. 

Nonetheless, it admits a nonstandard one. For consider the three-element algebra $\A = \langle \{ 0^{+}, 0^{-}, 1 \}; \land, \lor \rangle$ whose binary commutative operations are defined by the following tables:

\begin{center}
\begin{tabular}{ c | c | c | c |}
$\land$ & $0^{-}$ & $0^{+}$ & $1$ \\
\hline
$0^{-}$ & $0^{+}$ & $0^{+}$ & $0^{+}$  \\ 
\hline
$0^{+}$ &  & $0^{-}$ & $0^{+}$  \\ 
\hline
$1$ &  &  & $1$ \\ 
\hline
\end{tabular}
\quad
\begin{tabular}{ c | c | c | c |}
$\lor$ & $0^{-}$ & $0^{+}$ & $1$ \\
\hline
$0^{-}$ & $0^{+}$ & $0^{+}$ & $1$  \\ 
\hline
$0^{+}$ &  & $0^{-}$ & $1$  \\ 
\hline
$1$ &  &  & $1$ \\ 
\hline
\end{tabular}
\end{center}

\noindent Let also $\D_{2}$ be the two-element distributive lattice with universe $\{ 0, 1 \}$. Notice that the map $f \colon \A \to \D_{2}$ such that $f(1) = 1$ and $f(0^{+}) = f(0^{-}) = 0$ is a surjective homomorphism. Bearing in mind that ${\bf CPC}_{\land\lor}$ is complete with respect to the logical matrix $\langle \D_{2}, \{ 1 \} \rangle$, this implies that it is also complete with respect to the matrix $\langle \A, f^{-1}(\{ 1 \})\rangle$, that is $\langle \A, \{ 1 \} \rangle$. Finally, since $\{ 1 \}$ is the set of solutions of the equations in $\btau(x)$ in $\A$, we conclude that the class $\{ \A \}$ constitutes a $\btau$-algebraic semantics for ${\bf CPC}_{\land \lor}$.
\end{proof}

The nonstandard equational completeness theorem for ${\bf CPC}_{\land \lor}$ described above is a special instance of a general construction, devised in \cite[Thm.\ 3.1]{BlRe03}, that produces a (possibly nonstandard) algebraic semantics for every logic possessing an idempotent connective. While most familiar logics possess such a connective and, therefore, an algebraic semantics, the existence of logics that do not possess any algebraic semantics is known since \cite{BP89}, see also \cite{BlRe03,Font13,JMF14opL,Ra06a}. It is therefore sensible to wonder whether an intelligible characterization of logics with an algebraic semantics could possibly be obtained \cite{JGRa11}. In this paper we provide a positive answer to this question for a wide family of logics.


%


To this end, it is convenient to isolate some limits cases:\ a logic is said to be \textit{graph-based} when its language comprises only constant symbols and, possibly, a single unary connective. The algebraic models of graph-based logics are suitable directed graphs with distinguished elements. This makes them amenable to a series of combinatorial observations, culminating in a classification of graph-based logics with an algebraic semantics (Theorem \ref{Thm:graph-based-solved} and Section \ref{Sec:strongly-finite}).

To tackle the case of logics that are not graph-based, we first introduce a new method for constructing nonstandard algebraic semantics based on a universal algebraic trick known as \textit{Maltsev's Lemma}, which provides a description of congruence generation (Theorem \ref{Thm:sufficient}). As a consequence, we derive a characterization of logics with an algebraic semantics that, moreover, have at least one tautology in which a variable occurs (Theorem \ref{Thm:with-thms}). 

This result is subsequently specialized to the case of \textit{protoalgebraic} logics, i.e., logics $\vdash$ for which there exists a set of formulas $\Delta(x, y)$ that globally behaves as a very weak implication, in the sense that $\emptyset \vdash \Delta(x, x)$ and $x, \Delta(x, y) \vdash y$. More precisely, we prove that a nontrivial protoalgebraic logic has an algebraic semantics if and only if syntactic equality differs from logical equivalence, that is, there are two distinct logically equivalent formulas (Theorem \ref{Thm:protoalgebraic-logics}).\ It follows that most familiar protoalgebraic logics have an algebraic semantics which, however, can be nonstandard. For instance, while the local consequences of the modal systems ${\bf K}$, ${\bf K4}$, and ${\bf S4}$ possess a nonstandard algebraic semantics, the classes of algebras naturally associated with them, namely the varieties of modal, K4, and interior algebras, are not an algebraic semantics for them (Corollary \ref{Cor:modal-not-an-algebraic-semantics}).

Lastly, a logic is called \textit{locally tabular} when, up to logical equivalence, it has only finitely many $n$-ary formulas for every nonnegative integer $n$. A detailed characterization of locally tabular logics with an algebraic semantics is established in Section \ref{Sec:strongly-finite}. Accordingly, every locally tabular logic, whose language comprises at least an $n$-ary operation with $n \geq 2$, has an algebraic semantics (Proposition \ref{Prop:SF-no-graph-based}). 

The paper ends by considering the computational aspects of the problem of determining whether a logic has an algebraic semantics. We show that this problem is decidable both for logics presented by a finite set of finite logical matrices and for locally tabular logics presented by a finite Hilbert calculus. It becomes undecidable, however, when formulated for arbitrary (i.e., not necessarily locally tabular) logics presented by a finite Hilbert calculus (Theorem \ref{Thm:comput-thm}).

\section{Matrix semantics}

For a systematic presentation of matrix semantics, we refer the reader to \cite{BP86,BP89,BP92,CN20xx-book,Cz01,AAL-AIT-f,FJa09,FJaP03b}, while for a basic introduction to universal algebra we recommend \cite{Be11g,BuSa81}. Given an algebraic language $\mathscr{L}$, we denote by $Fm_{\mathscr{L}}$ the set of its formulas build up with a denumerable set $\Var$ of variables and by $\Fm_{\LL}$ the corresponding algebra. When $\mathscr{L}$ is clear from the context, we shall write $Fm$ and $\Fm$ instead of $Fm_{\mathscr{L}}$ and  $\Fm_{\LL}$. Generic elements of $\Var$ will be denoted by $x, y, z \dots$

A \textit{logic} $\vdash$ is then a consequence relation on the set $Fm_{\mathscr{L}}$, for some algebraic language $\LL$, that is \textit{substitution invariant} in the sense that for every substitution $\sigma$ on $Fm_{\mathscr{L}}$ (a.k.a.\ endomorphism of $\Fm_{\LL}$) and every $\Gamma \cup \{ \varphi \} \subseteq Fm_{\LL}$,
\[
\text{if }\Gamma \vdash \varphi \text{, then }\sigma[\Gamma] \vdash \sigma(\varphi).
\]
\noindent Given $\Gamma, \Delta \subseteq Fm$, we write $\Gamma \vdash \Delta$ when $\Gamma \vdash \delta$ for all $\delta \in \Delta$. Accordingly, $\Gamma \sineq \Delta$ stands for $\Gamma \vdash \Delta$ and $\Delta \vdash \Gamma$. For $\varphi, \psi \in Fm$, we abbreviate $\{ \varphi \} \sineq \{ \psi \}$ by $\varphi \sineq \psi$.

The language in which a logic $\vdash$ is formulated is denoted by $\LL_{\vdash}$. If $\LL \subseteq Fm_{\LL_{\vdash}}$, the $\LL$-\textit{fragment} of $\vdash$ is the restriction of $\vdash$ to formulas in $Fm_{\LL}$. Moreover, a logic $\vdash'$ is said to be an \textit{extension} of $\vdash$ if it is formulated in the same language as $\vdash$ and ${\vdash}  \subseteq {\vdash'}$.

Given a logic $\vdash$, we denote by $\textup{Cn}_{\vdash} \colon \mathcal{P}(Fm) \to \mathcal{P}(Fm)$ the closure operator naturally associated with $\vdash$, i.e., the map defined by the rule
\[
\textup{Cn}_{\vdash}(\Gamma) \coloneqq \{ \gamma \in Fm \colon \Gamma \vdash \gamma \}.
\]
The elements of $\textup{Cn}_{\vdash}(\emptyset)$ are said to be the \textit{theorems} of $\vdash$. By a \textit{rule} we understand an expression of the form $\Gamma \rhd \varphi$ where $\Gamma \cup \{ \varphi \} \subseteq Fm$. A rule $\Gamma \rhd \varphi$ is \textit{valid} in $\vdash$ if $\Gamma \vdash \varphi$. 

Given an algebra $\A$, a map $p \colon A^{n} \to A$ is said to be a \textit{polynomial function} of $\A$ if there are a formula $\varphi(x_{1},\dots, x_{n}, \vec{y})$ and a tuple $\vec{c} \in A$ such that 
\[
p(a_{1}, \dots, a_{n}) = \varphi^{\A}(a_{1}, \dots, a_{n}, \vec{c})
\]
for every $a_{1}, \dots, a_{n} \in A$. 

Algebras whose language is $\LL$ are called $\LL$-\textit{algebras}. Given a logic $\vdash$ and an $\LL_{\vdash}$-algebra $\A$, a set $F \subseteq A$ is said to be a \textit{deductive filter} of $\vdash$ on $\A$ when for every $\Gamma \cup \{ \varphi\} \subseteq Fm$ such that $\Gamma \vdash \varphi$ and every homomorphism $h \colon \Fm \to \A$, if $h[\Gamma] \subseteq F$, then $h(\varphi) \in F$. The set of deductive filters of $\vdash$ on $\A$ is a closure system, whose closure operator is denoted by $\textup{Fg}_{\vdash}^{\A}(\cdot) \colon \mathcal{P}(A) \to \mathcal{P}(A)$. When $a \in A$, we shall write $\textup{Fg}_{\vdash}^{\A}(a)$ as a shorthand for $\textup{Fg}_{\vdash}^{\A}(\{ a \})$.

Every logic $\vdash$ can be associated with a distinguished class of algebras, as we proceed to explain. Given an $\LL_{\vdash}$-algebra $\A$, let $\equiv_{\vdash}^{\A}$ be the congruence\footnote{This congruence is named after Tarski in \cite{AAL-AIT-f,FJa09,FJaP03b}.} of $\A$ defined for every $a, c \in A$ as
\begin{align}\label{Eq:Tarski}
    \begin{split}
a \equiv_{\vdash}^{\A}c \Longleftrightarrow& \text{ for every unary polynomial functions $p$ of $\A$,}\\
& \text{ }\textup{Fg}_{\vdash}^{\A}(p(a)) = \textup{Fg}_{\vdash}^{\A}(p(c)).
\end{split}
\end{align}
The \textit{algebraic counterpart} of $\vdash$ is the class of algebras
\[
\Alg(\vdash) \coloneqq \{ \A \colon \! \A \text{ is an $\LL_{\vdash}$-algebra and}\equiv_{\vdash}^{\A} \text{is the identity relation on }A\}.
\]
If $\vdash$ is classical propositional logic, $\Alg(\vdash)$ is the variety (a.k.a.\ equational class) of Boolean algebras.


A canonical way to present logics is by means of classes of (logical) matrices. By a \textit{matrix} we understand a pair $\langle \A, F \rangle$ where $\A$ is an algebra and $F \subseteq A$. In this case, $\A$ is called the \textit{algebraic reduct} of the matrix. The logic \textit{induced} by a class of similar matrices $\class{M}$ is the consequence relation $\vdash_{\class{M}}$ on $Fm$ defined for every $\Gamma \cup \{ \varphi \} \subseteq Fm$ as
\begin{align*}
\Gamma \vdash_{\class{M}} \varphi \Longleftrightarrow&\text{ for every $\langle \A, F \rangle \in \class{M}$ and homomorphism }h \colon \Fm \to \A, \\
 &\text{ if }h[\Gamma] \subseteq F\text{, then }h(\varphi)\in F.
\end{align*}
Accordingly, a class of matrices $\class{M}$ is said to be a \textit{matrix semantics} for a logic $\vdash$ when the latter coincides with $\vdash_{\class{M}}$. By the classical  Lindenbaum-Tarski method, every logic $\vdash$ is induced by the class of matrices $\{ \langle \Fm, \Gamma \rangle \colon \Gamma \subseteq Fm \text{ and }\textup{Cn}_{\vdash}(\Gamma) = \Gamma\}$, whence

\begin{Theorem}\label{Thm:matrix-semantics}
Every logic has a matrix semantics.
\end{Theorem}




A congruence $\theta$ of an algebra $\A$ is \textit{compatible} with a set $F \subseteq A$ when for every $a, c \in A$,
\[
\text{if }\langle a, c \rangle \in \theta \text{ and }a \in F\text{, then }c \in F.
\]
The \textit{Leibniz congruence} $\leibniz^{\A}F$ is the largest congruence on $\A$ compatible with $F$ and can be described as follows:

\begin{Proposition}[\protect{\cite[Thm.\ 4.23]{AAL-AIT-f}}]\label{Prop:Polynomial}
Let  $\A$ be an algebra, $F \subseteq A$, and $a, c\in A$. 
\begin{align*}
\langle a, c \rangle \in \leibniz^{\A}F\Longleftrightarrow& \text{ for every unary polynomial function $p$ of $\A$,}\\
&  \text{ }p(a)\in F\textrm{ if and only if }p(c)\in F.
\end{align*}
\end{Proposition}
\noindent This description of $\leibniz^{\A}F$ can be often simplified. For instance, if $\A$ is a Boolean algebra and $F$ one of its filters, then
\[
\langle a, c \rangle \in \leibniz^{\A}F \Longleftrightarrow a \to^{\A} c, c \to^{\A} a \in F.
\]

Given a matrix $\langle \A, F \rangle$, we denote by $\HHSinv(\langle \A, F \rangle)$ the class of all matrices $\langle \B, G \rangle$ for which there exists a surjective homomorphism $h \colon \B \to \A$ such that $G = h^{-1}[F]$.

\begin{Lemma}[\protect{\cite[Prop.\ 4.35.1]{AAL-AIT-f}}]\label{Lem:reduction}
If $\langle \A, F \rangle$ is a matrix and $\langle \B, G \rangle \in \HHSinv(\langle \A, F \rangle)$, then $\langle \A, F \rangle$ and $\langle \B, G \rangle$ induce the same logic.
\end{Lemma}

\noindent A matrix $\langle \A, F \rangle$ is said to be \textit{reduced} when $\leibniz^{\A}F$ is the identity relation. The matrix
\[
\langle \A, F \rangle^{\ast} \coloneqq \langle \A / \leibniz^{\A}F, F / \leibniz^{\A}F\rangle
\]
is always reduced. Similarly, given a class of matrices $\class{M}$, set
\[
\class{M}^{\ast} \coloneqq \{ \langle \A, F \rangle^{\ast} \colon \langle \A, F \rangle \in \class{M} \}.
\]
As $\langle \A, F \rangle \in \HHSinv(\langle \A, F \rangle^{\ast})$ for every matrix $\langle \A, F \rangle$, from Lemma \ref{Lem:reduction} we obtain:
\begin{Corollary}\label{Cor:reduction-same-logic}
The logics induced by $\class{M}$ and $\class{M}^{\ast}$ coincide, for every class of matrices $\class{M}$.
\end{Corollary}

By Theorem \ref{Thm:matrix-semantics} and Corollary \ref{Cor:reduction-same-logic}, every logic $\vdash$ has a matrix semantics of the form $\class{M}^{\ast}$. Moreover, if $\langle \A, F \rangle \in \class{M}^{\ast}$, then $F$ is a deductive filter of $\vdash$ on $\A$ and $\leibniz^{\A}F$ is the identity relation. In this case, by Proposition \ref{Prop:Polynomial}, for every pair of distinct $a, c \in A$, there is a unary polynomial function $p$ of $\A$ such that
\[
\textup{Fg}_{\vdash}^{\A}(F \cup \{ p(a) \}) \ne \textup{Fg}_{\vdash}^{\A}(F \cup \{ p(c) \}).
\]
In particular, this implies
\[
\textup{Fg}_{\vdash}^{\A}(p(a)) \ne \textup{Fg}_{\vdash}^{\A}(p(c)).
\]
By (\ref{Eq:Tarski}), we conclude that $\equiv_{\vdash}^{\A}$ is the identity relation, whence $\A \in \Alg(\vdash)$. Consequently,

\begin{Corollary}\label{Cor:the-place-of-reductions}
Every logic has a matrix semantics whose algebraic reducts belong to $\Alg(\vdash)$.
\end{Corollary}

When $\A = \Fm$, we shall drop the superscript $\Fm$ and write $\equiv_{\vdash}$ instead of $\equiv_{\vdash}^{\Fm}$. Two formulas $\varphi$ and $\psi$ are said to be \textit{logically equivalent} in $\vdash$ when $\varphi \equiv_{\vdash} \psi$, i.e., if
\[
\delta(\varphi, \vec{z}) \sineq \delta(\psi, \vec{z}), \text{ for all }\delta(x, \vec{z}) \in Fm.
\]
In the case of propositional classical logic, two formulas $\varphi$ and $\psi$ are logically equivalent in the above sense if and only if $\varphi \leftrightarrow \psi$ is a classical tautology.

Lastly, the \textit{free Lindenbaum-Tarski algebra} of $\vdash$ is $\Fm(\vdash) \coloneqq \Fm / {\equiv}_{\vdash}$. In the case of classical propositional logic, $\Fm(\vdash)$ is the free Boolean algebra with a denumerable set of free generators. More in general,

\begin{Lemma}[\protect{\cite[Prop.\ 5.75.2]{AAL-AIT-f}}]\label{Lem:free-algebra}
For every logic $\vdash$, $\Fm(\vdash)$ is the free algebra of $\Alg(\vdash)$ with a denumerable set of free generators.
\end{Lemma}

\noindent Validity of equations in the algebraic counterpart of a logic and logical equivalence are related as follows:

\begin{Lemma}\label{Lem:matrix-to-rules-1}
The following conditions are equivalent for a logic $\vdash$ and $\epsilon, \delta \in Fm$:
\benroman
\item\label{item:matrix-to-rules-0} $\Alg(\vdash) \vDash \epsilon \thickapprox \delta$;
\item\label{item:matrix-to-rules-1} $\vdash$ is induced by a class of matrices $\class{M}$ such that $\class{M} \vDash \epsilon \thickapprox \delta$;
\item\label{item:matrix-to-rules-2} $\epsilon$ and $\delta$ are logically equivalent;
\item\label{item:matrix-to-rules-4} $\Fm(\vdash) \vDash \epsilon \thickapprox \delta$.
\eroman
\end{Lemma}

\noindent In the above result, (\ref{item:matrix-to-rules-0})$\Rightarrow$(\ref{item:matrix-to-rules-1}) follows from Corollary \ref{Cor:the-place-of-reductions}, (\ref{item:matrix-to-rules-1})$\Rightarrow$(\ref{item:matrix-to-rules-2}) is obvious, (\ref{item:matrix-to-rules-2})$\Rightarrow$(\ref{item:matrix-to-rules-4}) is a consequence of the definition of $\equiv_{\vdash}$, and  (\ref{item:matrix-to-rules-4})$\Rightarrow$(\ref{item:matrix-to-rules-0}) follows from Lemma \ref{Lem:free-algebra}.

\section{Algebraic semantics}


A logic is said to have an algebraic semantics when it admits an equational completeness theorem. This concept originated in Blok and Pigozzi's monograph on algebraizable logics \cite[Sec.\ 2]{BP89} and was further investigated by Blok and Rebagliato \cite{BlRe03} and Raftery \cite{Ra06a}. In order to review it, let $\class{K}$ be a class of similar algebras. 

The \textit{equational consequence relative to $\class{K}$} is the consequence relation $\vDash_{\class{K}}$ on the set of equations $Eq$ in variables $\Var$ defined for every $\Theta \cup \{ \epsilon \thickapprox \delta \} \subseteq Eq$ as follows:
\begin{align*}
\Theta \vDash_{\class{K}} \epsilon \thickapprox \delta \Longleftrightarrow&\text{ for every $\A \in \class{K}$ and homomorphism }h \colon \Fm \to \A, \\
 &\text{ if }h(\varphi) = h(\psi) \text{ for all }\varphi \thickapprox \psi \in \Theta\text{, then }h(\epsilon) = h(\delta).
\end{align*}
For every set of equations $\btau(x)$ and set of formulas $\Gamma \cup \{\varphi \}$, we shall abbreviate
\[
\{ \epsilon(\varphi) \thickapprox \delta(\varphi) \colon \epsilon \thickapprox \delta \in \btau \} \text{ as }\btau(\varphi), \text{ and } \bigcup_{\gamma \in \Gamma} \btau(\gamma) \text{ as }\btau[\Gamma].
\]
\begin{law}
A logic $\vdash$ is said to have an \textit{algebraic semantics} if there are a set of equations $\btau(x)$ and a class $\class{K}$ of $\LL_{\vdash}$-algebras such that for all $\Gamma \cup \{ \varphi \} \subseteq Fm$,
\[
\Gamma \vdash \varphi \Longleftrightarrow \btau[\Gamma] \vDash_{\class{K}} \btau(\varphi).
\]
In this case $\class{K}$ is said to be a $\btau$-\textit{algebraic semantics} (or simply an algebraic semantics) for $\vdash$. An algebraic semantics $\class{K}$ for a logic $\vdash$ is called \textit{standard} when $\class{K} \subseteq \Alg(\vdash)$. It is said to be \textit{nonstandard} otherwise.
\end{law}

The completeness theorem of classical propositional logic ${\bf CPC}$ with respect to Boolean algebras states that these form a $\btau$-algebraic semantics for ${\bf CPC}$ where $\btau = \{ x \thickapprox 1 \}$. A logic, however, can have nonstandard algebraic semantics. For instance, Glivenko's Theorem \cite{Gli29} implies that Heyting algebras form a $\btau$-algebraic semantics for ${\bf CPC}$ where $\btau = \{ \lnot \lnot x \thickapprox 1 \}$, as noticed in the introduction.

Given a set of equations $\btau(x)$ and an algebra $\A$, we set
\[
\btau(\A) \coloneqq \{ a \in A \colon \A \vDash \btau(a) \}.
\]
The following observations are immediate consequences of the definition of an algebraic semantics.

\begin{Proposition}[\protect{\cite[Thm.\ 2.4]{BP89}}]\label{Prop:matrix-sem-alg-sem}
A class of algebras $\class{K}$ is a $\btau$-algebraic semantics for a logic $\vdash$ if and only if $\{ \langle \A, \btau(\A) \rangle \colon \A \in \class{K} \}$ is a matrix semantics for $\vdash$.
\end{Proposition}

\begin{Proposition}[\protect{\cite[Thm.\ 2.16]{BlRe03}}]\label{Prop:Suszko-congruence}
If $\vdash$ has a $\btau$-algebraic semantics, then
\[
x, \varphi(\epsilon, \vec{z}) \sineq \varphi(\delta, \vec{z}), x
\]
for all $\epsilon \thickapprox \delta \in \btau$ and $\varphi(v, \vec{z}) \in Fm$.
\end{Proposition}


In the remaining part of this section, we shall rephrase in purely algebraic parlance the problem of determining whether a logic has a $\btau$-algebraic semantics. To this end, given an algebra $\A$ and $X \subseteq A^{2}$, we denote by $\textup{Cg}^{\A}(X)$ the congruence of $\A$ generated by $X$.
\begin{law}
For every $\Gamma \subseteq Fm$ and set of equations $\btau(x)$, define
\[
\theta(\Gamma, \btau) \coloneqq \textup{Cg}^{\Fm}(\{ \langle \epsilon(\gamma), \delta(\gamma)\rangle \colon \epsilon \thickapprox \delta \in \btau \text{ and }\gamma \in \Gamma \}).
\]
\end{law}
In the statement of the next result we identify equations with pairs of formulas, whence the notation $\epsilon \thickapprox \delta$ stands for $\langle \epsilon, \delta\rangle$. Notice that, under this convention, expressions of the form $\btau(\varphi) \subseteq \theta(\Gamma, \btau)$ make sense.

\begin{Proposition}\label{Prop:alg-sem-meaning}
A logic $\vdash$ has a $\btau$-algebraic semantics if and only if $\Gamma \vdash \varphi$ for all $\Gamma \cup \{ \varphi \} \subseteq Fm$ such that $\btau(\varphi) \subseteq \theta(\Gamma, \btau)$.
\end{Proposition}

\begin{proof}
We begin by proving the ``if'' part. Define
\[
\class{K} \coloneqq \{ \Fm/ \theta(\Gamma, \btau) \colon \Gamma \subseteq Fm \text{ and }\textup{Cn}_{\vdash}(\Gamma) = \Gamma \}.
\]
By assumption, for every $\varphi \in Fm$ and $\Gamma \subseteq Fm$ such that $\textup{Cn}_{\vdash}(\Gamma) = \Gamma$,
\[
\varphi / \theta(\Gamma, \btau) \in \btau(\Fm/ \theta(\Gamma, \btau)) \Longleftrightarrow \btau(\varphi) \subseteq \theta(\Gamma, \btau) \Longleftrightarrow \varphi \in \Gamma.
\]
Consequently, for every $\Gamma \subseteq Fm$ such that $\textup{Cn}_{\vdash}(\Gamma) = \Gamma$,
\[
\langle \Fm, \Gamma \rangle \in \HHSinv(\langle \Fm / \theta(\Gamma, \btau),  \btau(\Fm/ \theta(\Gamma, \btau)) \rangle).
\]
Hence, by Lemma \ref{Lem:reduction}, the following classes of matrices induce the same logic:
\begin{align*}
\class{M}_{1} &\coloneqq \{ \langle \Fm, \Gamma \rangle \colon \Gamma \subseteq Fm \text{ and }\textup{Cn}_{\vdash}(\Gamma) = \Gamma \}\\
\class{M}_{2} & \coloneqq \{ \langle \A,  \btau(\A) \rangle \colon \A \in \class{K} \}.
\end{align*}
As $\class{M}_{1}$ is a matrix semantics for $\vdash$, so is $\class{M}_{2}$. By Proposition \ref{Prop:matrix-sem-alg-sem} we conclude that $\class{K}$ is a $\btau$-algebraic semantics for $\vdash$.

To prove the ``only if'' part, let $\class{K}$ be a $\btau$-algebraic semantics for $\vdash$. Then consider $\Gamma \cup \{ \varphi \} \subseteq Fm$ such that $\btau(\varphi) \subseteq \theta(\Gamma, \btau)$. In order to prove that $\Gamma \vdash \varphi$, it suffices to show that for every $\A \in \class{K}$ and homomorphism $h \colon \Fm \to \A$, if $h[\Gamma] \subseteq \btau(\A)$, then $h(\varphi) \in \btau(\A)$. To this end, consider $\A \in \class{K}$ and a homomorphism $h \colon \Fm \to \A$ such that $h[\Gamma] \subseteq \btau(\A)$. Observe that the kernel $\textup{Ker}(h)$ of $h$ contains the generators of $\theta(\Gamma, \btau)$, whence
\[
\btau(\varphi) \subseteq \theta(\Gamma, \btau) \subseteq \textup{Ker}(h).
\]
But this amounts to $h(\varphi) \in \btau(\A)$. Thus, we conclude that $\Gamma \vdash \varphi$, as desired.
\end{proof}

\begin{Corollary}[\protect{\cite[Thm.\ 2.15]{BlRe03}}]\label{Cor:persist-extensions}
The property of having an algebraic semantics persists in extensions of a logic.
\end{Corollary}

\section{Tame examples: assertional logics}

The notion of an assertional logic was introduced by Pigozzi in \cite{Pi91} and further investigated by Blok and Raftery among others \cite{BlRa08}. For the present purpose, the interest of assertional logics comes from the fact that they are typical examples of logics with an algebraic semantics. To review these facts, it is convenient to recall some basic definitions:

\begin{law}
A class of matrices $\class{M}$ is said to be \textit{unital} if $F$ is either empty or a singleton, for every $\langle \A, F \rangle \in \class{M}$.
\end{law}

\begin{Proposition}\label{Prop:unital-rules}
A logic $\vdash$ has a unital matrix semantics if and only if for every $\varphi(v, \vec{z}) \in Fm$,
\[
x, y, \varphi(x, \vec{z}) \vdash \varphi(y, \vec{z}).
\]
\end{Proposition}

\begin{proof}
This is attributed to Suszko in \cite[Pag.\ 67]{Rt93}. For a detailed argument, however, the reader may consult the proof of \cite[Thm.\ 10]{AFRM15}.
\end{proof}

Let $\A$ be an algebra and $\varphi(x)$ a formula. If the term-function $\varphi^{\A} \colon A \to A$ is a constant function, we denote by $\varphi^{\A}$ both the term-function and its unique value.

\begin{law}
A logic $\vdash$ is said to be \textit{assertional} if it is induced by a class of matrices $\class{M}$ for which there is a formula $\varphi(x)$ such that $\varphi^{\A} \colon A \to A$ is a constant function and $F = \{ \varphi^{\A} \}$, for all $\langle \A, F \rangle \in \class{M}$ .
\end{law}
\noindent Observe that, in this case, $\class{M}$ is unital and $\varphi$ a theorem of $\vdash$. Moreover, the class $\class{K}$ of algebraic reducts of matrices in $\class{M}$ is a $\btau$-algebraic semantics for $\vdash$ where $\btau = \{ x \thickapprox \varphi(x) \}$. Consequently, we obtain the following:


\begin{Proposition}\label{Prop:assertional} Assertional logics have an algebraic semantics.
\end{Proposition}

\begin{exa}[\textsf{Intermediate logics}]
A \textit{Heyting algebra} is a structure $\A = \langle A; \land, \lor, \to, 0, 1 \rangle$ comprising a bounded lattice $\langle A; \land, \lor, 0, 1 \rangle$ such that for every $a, b, c \in A$,
\[
a \land b \leq c \Longleftrightarrow a \leq b \to c,
\]
where $x \leq y$ is a shorthand for $x = x \land y$. The assertional logic $\vdash_{\class{K}}$, associated with a variety $\class{K}$ of Heyting algebras, is defined by the rule
\begin{equation}\label{Eq:intermediate-def}
\Gamma \vdash_{\class{K}} \varphi \Longleftrightarrow \btau[\Gamma] \vDash_{\class{K}} \btau(\varphi),
\end{equation}
where $\btau = \{ x \thickapprox 1 \}$. Logics arising in this way are called \textit{intermediate} or \textit{superintuitionistic} \cite{ChZa97}. Notice that if $\class{K}$ is the variety of Boolean algebras (resp.\ of all Heyting algebras), then $\vdash_{\class{K}}$ is classical (resp.\ intuitionistic) propositional logic. 
\qed
\end{exa}

\begin{exa}[\textsf{Global modal logics}]
A \textit{modal algebra} is a structure $\A = \langle A; \land, \lor, \lnot, \Box, 0, 1 \rangle$ comprising a Boolean algebra $\langle A; \land, \lor, \lnot, 0, 1 \rangle$ and a unary operation $\Box$ such that for every $a, c \in A$,
\[
\Box 1 = 1 \text{ and }\Box(a \land c) = \Box a \land \Box c.
\]
The assertional logic $\vdash_{\class{K}}$, associated with a variety $\class{K}$ of modal algebras, is defined by the rule (\ref{Eq:intermediate-def}). Logics arising in this way have been called \textit{global modal logics} \cite{MK07c}.

Notice that assertionality persists in fragments of intermediate and global modal logics in which $1$ is term-definable. 
\qed
\end{exa}


As we mentioned, assertional logics necessarily have theorems. The next definition extends this concept beyond logics with theorems.

\begin{law}
A logic lacking theorems is said to be \textit{almost assertional} if it has a matrix semantics $\class{M}$ for which there is a formula $\varphi(x)$ such that $\varphi^{\A} \colon A \to A$ is a constant function and $F = \{ \varphi^{\A} \}$, for all $\langle \A, F \rangle \in \class{M}$ with $F \ne \emptyset$.
\end{law}

In general, almost assertional logics need not have an algebraic semantics. In order to explain when this is the case, we rely on the following (see \cite[Lem.\ 4.3]{AAL-AIT-f}, if necessary):

\begin{Lemma}\label{Lem:no-thms-matrix}
Let $\vdash$ be a logic lacking theorems. If $\class{M}$ is a matrix semantics for $\vdash$, then so is
\[
\class{N} \coloneqq \{ \langle \A, F \rangle \in \class{M} \colon F \ne \emptyset \} \cup \{ \langle \Fm, \emptyset \rangle \}.
\]
\end{Lemma}

\noindent Almost assertional logics with an algebraic semantics can be characterized as follows:

\begin{Proposition}\label{Prop:almost-assertional}
An almost assertional logic $\vdash$ has an algebraic semantics if and only if there exists a set of equations $\btau(x)$ for which the following conditions hold:
\benroman
\item\label{Eq:almost-assertional-1} There is no substitution $\sigma$ such that $\sigma(\epsilon) = \sigma(\delta)$ for all $\epsilon \thickapprox \delta \in \btau$;
\item\label{Eq:almost-assertional-2} For every $\epsilon \thickapprox \delta \in \btau$ and $\varphi(v, \vec{z}) \in Fm$,
\[
x, \varphi(\epsilon, \vec{z}) \sineq \varphi(\delta, \vec{z}), x.
\]
\eroman
\end{Proposition}

\begin{proof}
Let $\vdash$ be an almost assertional logic. Then $\vdash$ has a matrix semantics $\class{M}$ for which there is a formula $\varphi(x)$ such that $\varphi^{\A} \colon A \to A$ is a constant function and $F = \{ \varphi^{\A} \}$, for all $\langle \A, F \rangle \in \class{M}$ with $F \ne \emptyset$. Since $\class{M}$ is a matrix semantics for $\vdash$, we get
\begin{equation}\label{Eq:zzz-almost-assertional-zzz--2}
y \vdash \varphi(x).
\end{equation}

To prove the ``only if'' part, suppose that $\vdash$ has a $\btau$-algebraic semantics $\class{K}$. As $\vdash$ lacks theorems, there are $\A \in \class{K}$ and $a \in A$ such that $\varphi^{\A}(a) \notin \btau(\A)$. Together with (\ref{Eq:zzz-almost-assertional-zzz--2}), this implies $\btau(\A) = \emptyset$. Consequently, $\btau$ satisfies condition (\ref{Eq:almost-assertional-1}). Moreover, $\btau$ satisfies condition (\ref{Eq:almost-assertional-2}) by Proposition \ref{Prop:Suszko-congruence}.

Then we turn to prove the ``if'' part. As $\class{M}$ is a matrix semantics for $\vdash$, the same holds for $\class{M}^{\ast}$ by Corollary \ref{Cor:reduction-same-logic}. Furthermore, for every $\langle \A, F \rangle \in \class{M}^{\ast}$,
\begin{equation}\label{Eq:zzz-almost-assertional-zzz}
\text{if }a \in F \text{, then }a \in \btau(\A).
\end{equation}
To prove this, consider $\langle \A, F \rangle \in \class{M}^{\ast}$, $a \in F$, and $\epsilon \thickapprox \delta \in \btau$. From condition (\ref{Eq:almost-assertional-2}) it follows that for every formula $\varphi(v,\vec{z})$ and $\vec{c} \in A$,
\[
\varphi^{\A}(\epsilon^{\A}(a), \vec{c}) \in F \Longleftrightarrow \varphi^{\A}(\delta^{\A}(a), \vec{c}) \in F.
\]
This means that for every unary polynomial function $p$ of $\A$,
\[
p(\epsilon^{\A}(a)) \in F \Longleftrightarrow p(\delta^{\A}(a)) \in F.
\]
By Proposition \ref{Prop:Polynomial}, we conclude that $\langle \epsilon^{\A}(a), \delta^{\A}(a) \rangle \in \leibniz^{\A}F$. Since $\langle \A, F \rangle \in \class{M}^{\ast}$, the matrix $\langle \A, F \rangle$ is reduced, i.e., $\leibniz^{\A}F$ is the identity relation on $A$. Consequently, $\epsilon^{\A}(a) = \delta^{\A}(a)$, establishing (\ref{Eq:zzz-almost-assertional-zzz}).

Then set
\[
\brho \coloneqq \btau \cup \{ x \thickapprox \varphi(x) \}
\]
and observe that for every $\langle \A, F \rangle \in \class{M}^{\ast}$,
\begin{equation}\label{Eq:zzz-almost-assertional-zzz-2}
\text{if }F \ne \emptyset, \text{ then }\brho(\A) = F.
\end{equation}
To prove this, consider $\langle \A, F \rangle \in \class{M}^{\ast}$ such that $F \ne \emptyset$. Then $\langle \A, F \rangle = \langle \B / \leibniz^{\B}G, G / \leibniz^{\B}G\rangle$ for some $\langle \B, G \rangle \in \class{M}$ such that $G = \emptyset$. By assumption, $\varphi^{\B} \colon B \to B$ is a constant function and $G = \{ \varphi^{\B} \}$. This easily implies 
\[
F = \{ a \in A \colon a = \varphi^{\A}(a) \}.
\]
Together with (\ref{Eq:zzz-almost-assertional-zzz}), the above display implies (\ref{Eq:zzz-almost-assertional-zzz-2}).

Recall that $\class{M}^{\ast}$ is a matrix semantics for $\vdash$ and that $\vdash$ lacks theorems. Then by Lemma \ref{Lem:no-thms-matrix} the following is also a matrix semantics for $\vdash$:
\[
\class{N} \coloneqq \{ \langle \A, F \rangle \in \class{M}^{\ast} \colon F \ne \emptyset \} \cup \{ \langle \Fm, \emptyset \rangle \}.
\]
Moreover,
\begin{equation}\label{Eq:zzz-almost-assertional-zzz-3}
\text{if }\langle \A, F \rangle \in \class{N}\text{, then }F = \brho(\A).
\end{equation}
To prove this, consider $\langle \A, F \rangle \in \class{N}$. If $F \ne \emptyset$, then $F = \brho(\A)$ by (\ref{Eq:zzz-almost-assertional-zzz-2}). Then we consider the case where $F = \emptyset$. By definition of $\class{N}$, in this case $\langle \A, F \rangle = \langle \Fm, \emptyset \rangle$. By condition (\ref{Eq:almost-assertional-1}) we get $\brho(\Fm) \subseteq \btau(\Fm) = \emptyset$, whence $\brho(\Fm) = \emptyset$ as desired. This establishes (\ref{Eq:zzz-almost-assertional-zzz-3}). 

As a consequence,
\[
\class{N} = \{ \langle \A, \brho(\A) \rangle \colon \A \in \class{K} \},
\]
where $\class{K}$ is the class of algebraic reducts of matrices in $\class{N}$. Since $\class{N}$ is a matrix semantics for $\vdash$, we can apply Proposition \ref{Prop:matrix-sem-alg-sem} obtaining that $\class{K}$ is a $\brho$-algebraic semantics for $\vdash$.
\end{proof}

The following result will be used later on.

\begin{Proposition}\label{Prop:rules-almost-assertional}
A logic $\vdash$ is either assertional or almost assertional if and only if there is a formula $\psi(x)$ such that $y \vdash \psi(x)$ and for every $\varphi(v, \vec{z}) \in Fm$,
\[
x, y, \varphi(x, \vec{z}) \vdash \varphi(y, \vec{z}).
\]
\end{Proposition}

\begin{proof}
Observe that a logic $\vdash$ is either assertional or almost assertional if and only if it has a matrix semantics $\class{M}$ for which there is a formula $\psi(x)$ such that $\psi^{\A} \colon A \to A$ is a constant function and $F = \{ \psi^{\A} \}$, for all $\langle \A, F \rangle \in \class{M}$ with $F \ne \emptyset$. 

To prove the ``only if' part, consider a logic $\vdash$ that is either assertional or almost assertional. By the above remarks, $\vdash$ has a unital matrix semantics $\class{M}$ and a formula $\psi(x)$ such that $y \vdash \psi(x)$. Hence, with an application of Proposition \ref{Prop:unital-rules}, we obtain $x, y, \varphi(x, \vec{z}) \vdash \varphi(y, \vec{z})$ for every $\varphi(v, \vec{z}) \in Fm$, as desired.

Then we turn to prove the ``if' part. Consider a logic $\vdash$ for which there is a formula $\psi(x)$ such that $y \vdash \psi(x)$ and $x, y, \varphi(x, \vec{z}) \vdash \varphi(y, \vec{z})$ for every $\varphi(v, \vec{z}) \in Fm$. By Proposition \ref{Prop:unital-rules}, $\vdash$ has a unital matrix semantics $\class{M}$. The fact that $\class{M}$ is unital and $y \vdash \psi(x)$ imply that $\psi^{\A} \colon A \to A$ is a constant function and $F = \{ \psi^{\A} \}$, for all $\langle \A, F \rangle \in \class{M}$ with $F \ne \emptyset$. Thus, by the remark at the beginning of the proof, $\vdash$ is either assertional or almost assertional.
\end{proof}

\section{Getting rid of limit cases}

While characterizing logics with an algebraic semantics, it is convenient to treat separately some limit cases. These do not present special difficulties and will be examined briefly is in this section.

\begin{law}
Let $\vdash$ be a logic.
\benormal
\item $\vdash$ is said to be \textit{inconsistent} if $\Gamma \vdash \varphi$ for all $\Gamma \cup \{ \varphi \} \subseteq Fm$. 
\item $\vdash$ is said to be \textit{almost inconsistent} if it lacks theorems and $\Gamma \vdash \varphi$ for all $\Gamma \cup \{ \varphi \} \subseteq Fm$ such that $\Gamma \ne \emptyset$. 
\item $\vdash$ is said to be \textit{trivial} if it is either inconsistent or almost inconsistent.
\enormal
\end{law}
\noindent Notice that a logic $\vdash$ is trivial precisely when $x \vdash y$.

\begin{Proposition}\label{Prop:trivial}
The following conditions hold for a logic $\vdash$.
\benroman
\item\label{trivial-item1} If $\vdash$ is inconsistent, then it has an algebraic semantics.
\item\label{trivial-item2} If $\vdash$ is almost inconsistent, then it has an algebraic semantics if and only if $\LL_{\vdash}$ comprises either two distinct constants or a non-constant connective.
\eroman
\end{Proposition}

\begin{proof}
(\ref{trivial-item1}): Observe that if $\vdash$ is inconsistent, then any class of algebras $\class{K}$ is an $\btau$-algebraic semantics for $\vdash$ where $\btau = \emptyset$.

(\ref{trivial-item2}): Let $\vdash$ be almost inconsistent. To prove the ``only if'' part, suppose, with a view to contradiction, that $\vdash$ has a $\btau$-algebraic semantics $\class{K}$, but that $\LL_{\vdash}$ is either empty or comprises only a constant symbol $1$. Due to the poor language in which $\vdash$ is formulated,
\[
\btau \subseteq \{ x \thickapprox x, x \thickapprox 1, 1 \thickapprox x, 1 \thickapprox 1 \}.
\]
We can assume without loss of generality that $\btau$ does not contain any trivially valid equation, whence $\btau \subseteq \{ x \thickapprox 1, 1 \thickapprox x \}$. Furthermore, by symmetry we can assume that $\btau \subseteq \{ x \thickapprox 1 \}$. Thus, either $\btau = \emptyset$ or $\btau = \{ x \thickapprox 1 \}$. In both cases, the fact that $\class{K}$ is a $\btau$-algebraic semantics for $\vdash$ implies $\emptyset \vdash 1$, contradicting the fact that $\vdash$ lacks theorems.

To prove the ``if'' part, suppose first that $\LL_{\vdash}$ comprises a non-constant connective $f(x_{1}, \dots, x_{n})$. Then define
\[
\btau \coloneqq \{ f(x, \dots, x) \thickapprox f(f(x, \dots, x), \dots, f(x, \dots, x)) \}.
\]
Clearly, $\btau(\Fm) = \emptyset$. As $\vdash$ is almost inconsistent, this implies that $\class{K} \coloneqq \{ \Fm \}$ is a $\btau$-algebraic semantics for $\vdash$. 

The case where $\LL_{\vdash}$ comprises two constant symbols $1$ and $0$ is treated similarly, taking $\btau \coloneqq \{ 1 \thickapprox 0 \}$.
\end{proof}

\section{A sufficient condition for algebraic semantics}

\begin{law}
An algebraic language $\LL$ is said to be \textit{graph-based} if the arities of its connectives are $\leq 1$ and, moreover, $\LL$ comprises at most one unary connective.\footnote{This terminology comes from the fact that algebras in graph-based languages can be naturally identified with certain directed graphs with a set of distinguished elements.} Similarly, a logic $\vdash$ is said to be \textit{graph-based} when $\LL_{\vdash}$ is. 
\end{law}

Given a formula $\varphi$, we denote by $\Var(\varphi)$ the set of variables really occurring in $\varphi$. The aim of this section is to establish the following result:

\begin{Theorem}\label{Thm:sufficient}
Let $\vdash$ be a logic that is not graph-based. If there are two distinct logically equivalent formulas $\varphi$ and $\psi$ such that $\Var(\varphi) \cup \Var(\psi) = \{ x \}$, then $\vdash$ has an algebraic semantics.
\end{Theorem}

One of the consequences of the above result is that every logic is ``essentially'' equivalent to one with an algebraic semantics. In order to make this observation mathematically precise, we shall recall some basic concepts. 

Given two algebraic languages $\LL$ and $\LL'$, a \textit{translation} of $\LL$ into $\LL'$ is a map that associates a (possibly complex) formula $f_{\btau}(x_{1}, \dots, x_{n})$ of $\LL'$ with every $n$-ary connective $f$ of $\LL$. Furthermore, given an algebraic language $\LL$, let $\LL^{\dagger}$ be the algebraic language obtained by replacing every constant symbol $\textbf{c}$ of $\LL$ with a new unary connective $\textbf{c}(x)$. Every $\LL$-algebra $\A$ can then be transformed into an $\LL^{\dagger}$-algebra $\A^{\dagger}$ by interpreting each new operation $\textbf{c}(x)$ as constant the function with value $\textbf{c}^{\A}$.

Lastly, two logics $\vdash$ and $\vdash'$ are said to be \textit{term-equivalent} if there are matrix semantics $\class{M}$ and $\class{M}'$ for $\vdash$ and $\vdash'$, respectively, and translations $\btau$ of $\LL_{\vdash}^{\dagger}$ into $\LL_{\vdash'}^{\dagger}$ and $\brho$ of $\LL_{\vdash'}^{\dagger}$ into $\LL_{\vdash}^{\dagger}$ such that
\[
\langle \A^{\dagger}, F \rangle =\langle \A^{\dagger\btau\brho}, F \rangle \; \text{ and } \; \langle \B^{\dagger}, G \rangle = \langle \B^{\dagger\brho\btau}, G \rangle
\]
for every $\langle \A, F \rangle \in \class{M}'$ and $\langle \B, G \rangle \in \class{M}$.\footnote{This notion of term-equivalence appeared first in \cite{JaMor19-1,JaMor19lh,JaMor19-3}. Even if in these papers term-equivalence was defined by means of the so-called \textit{Suszko reduced models} of a logic \cite{Cze03}, the present definition is easily seen to be equivalent to the original one.}

\begin{Corollary}
Every logic is term-equivalent to one with an algebraic semantics.
\end{Corollary}

\begin{proof}
Consider an arbitrary logic $\vdash$. Furthermore, let $\class{M}$ be a matrix semantics for $\vdash$ (it exists by Theorem \ref{Thm:matrix-semantics}). For every $\langle \A, F \rangle \in \class{M}$, let also $\A^{+}$ be the expansion of $\A$ with a new binary operation $+ \colon A^{2} \to A$ that is interpreted as the projection of the first coordinate, i.e.,
\[
a + c \coloneqq a \text{, for all }a, c \in A.
\]
Then let $\vdash^{+}$ be the logic induced by the class of matrices
\[
\class{M}^{+} \coloneqq \{ \langle \A^{+}, F \rangle : \langle \A, F \rangle \in \class{M} \}.
\]
From the definition of $\class{M}^{+}$ it follows that $\vdash$ and $\vdash^{+}$ are term-equivalent. 

Now, from the definition of $\vdash^{+}$ it follows that the formulas $x$ and $x + x$ are logically equivalent in $\vdash^{+}$. Since $\vdash^{+}$ is not graph-based, we can apply Theorem \ref{Thm:sufficient} obtaining that it has an algebraic semantics.
\end{proof}

The above result implies that ``having an algebraic semantics'' is not a property of clones or free algebras, thus confirming the fragility of the notion.

The remaining part of the section is devoted to proving Theorem \ref{Thm:sufficient}. We rely on the description of congruence generation given by \textit{Maltsev's Lemma} \cite[Thm.\ 4.17]{Be11g}, which we proceed to recall:

\begin{Theorem}\label{Thm:Maltsev-lemma}
Let $\A$ be an algebra, $X \subseteq A \times A$, and $a, c \in A$. Then $\langle a, c \rangle \in \textup{Cg}^{\A}(X)$ if and only if there are $e_{0}, \dots, e_{n} \in A$, $\langle b_{0}, d_{0}\rangle, \dots, \langle b_{n-1}, d_{n-1}\rangle \in X$, and unary polynomial functions $p_{0}, \dots, p_{n-1}$ of $\A$ such that
\[
a = e_{0}, c = e_{n} \text{, and }\{ e_{i}, e_{i+1} \} = \{ p_{i}(b_{i}), p_{i}(d_{i}) \}\text{, for every }i < n.
\]
\end{Theorem}

Lastly, given a formula $\gamma$, the \textit{subformula tree} $T(\gamma)$ of $\gamma$ is defined by induction on the construction of $\gamma$ as follows. If $\gamma$ is a variable or a constant, then $T(\gamma)$ is the trivial tree whose unique element is labeled by $\gamma$. Moreover, if $\gamma = f(\delta_{1}, \dots, \delta_{n})$ for some $n$-ary connective $f$ and formulas $\delta_{1}, \dots, \delta_{n}$, then $T(\gamma)$ is obtained by adding a common bottom element labeled by $f$ to the disjoint union of $T(\delta_{1}), \dots, T(\delta_{n})$. 

\begin{Notation}\emph{Given a formula $\gamma$, we denote by $T(\gamma)^{-x}$ the tree obtained by removing from $T(\gamma)$ the leaves labeled by the variable $x$.}
\end{Notation}

We are now ready to prove the main result of this section:

\begin{proof}[Proof of Theorem \ref{Thm:sufficient}.]
Since $\vdash$ is not graph-based, either $\LL_{\vdash}$ comprises two distinct unary connectives $\Box$ and $\Diamond$ or it comprises at least an $n$-ary connective with $n \geq 2$. We shall detail the case where $\LL_{\vdash}$ comprises $\Box x$ and $\Diamond x$, as the other one is analogous. 

Suppose that $\vdash$ has two distinct logically equivalent formulas $\varphi$ and $\psi$ such that $\Var(\varphi) \cup \Var(\psi) = \{ x \}$. By Lemma \ref{Lem:matrix-to-rules-1}, the logic $\vdash$ has a matrix semantics $\class{M}$ validating the equation $\varphi \thickapprox \psi$. Notice that we can safely assume that
\begin{equation}\label{Eq:x-appear}
\Var(\varphi) = \Var(\psi) = \{ x \}.
\end{equation}
For suppose that (\ref{Eq:x-appear}) does not hold. Then by symmetry we can assume that $\Var(\varphi) = \{ x \}$ and $\Var(\psi) = \emptyset$. Now, from $\class{M} \vDash \varphi(x) \thickapprox \psi(x)$ it follows $\class{M}\vDash \varphi(\varphi(x)) \thickapprox \psi(\varphi(x))$. Since $\Var(\psi) = \emptyset$, the formula $\psi$ is closed, whence $\psi(\varphi(x)) = \psi(x)$. Thus,  $\class{M}\vDash \varphi(\varphi(x)) \thickapprox \psi(x)$. Together with $\class{M} \vDash \varphi(x) \thickapprox \psi(x)$, this implies 
\[
\class{M} \vDash \varphi(x) \thickapprox \varphi(\varphi(x)).
\] 
We have two cases: either $\varphi = x$ or $\varphi \ne x$. First suppose that $\varphi = x$. As $\psi$ is a closed formula, $\class{M} \vDash x \thickapprox \psi$ implies that algebraic reducts of the matrices in $\class{M}$ are trivial. Because $\class{M}$ is a matrix semantics for $\vdash$, this yields $x \vdash y$ which, in turn, means that $\vdash$ is trivial. As the language of $\vdash$ comprises a non-constant symbol, from Proposition \ref{Prop:trivial} it follows that $\vdash$ has an algebraic semantics, thus concluding the proof of the theorem. Then it only remains to consider the case where $\varphi \ne x$. In this case, $\varphi(\varphi(x)) \ne \varphi(x)$. Thus we can replace $\psi$ by $\varphi(\varphi(x))$, validating both the assumption of the theorem and (\ref{Eq:x-appear}). Accordingly, from now we shall assume that (\ref{Eq:x-appear}) holds.

Let $k \in \omega$ be greater than the length of all branches in $T(\varphi)$ and $T(\psi)$. As $\varphi = \varphi(x)$ and $\psi = \psi(x)$, it makes sense to define
\[
\varphi' \coloneqq \Box^{2k}\Diamond \varphi(\Box^{k} \Diamond x) \text{ and }\psi'\coloneqq \Box^{2k}\Diamond \psi(\Box^{k} \Diamond x).
\]
Furthermore, since $\class{M} \vDash \varphi \thickapprox \psi$,
\begin{equation}\label{Eq:M-proves-equation}
\class{M} \vDash \varphi' \thickapprox \psi'.
\end{equation}

Now, let $m$ be the number of occurrences of $x$ in $\varphi'(x)$. Then, consider distinct variables $x_{1}, \dots, x_{m}$ and let $\hat{\varphi}'(x_{1}, \dots, x_{m})$ be the formula obtained by replacing the $j$-th occurrence of $x$ in $\varphi'$ by $x_{j}$ for all $j \leq m$. Clearly,
\[
\varphi' = \hat{\varphi}'(x, \dots, x).
\]

\begin{Claim}\label{Claim:first}
There are no $\alpha_{1}, \dots, \alpha_{m}, \beta \in Fm$ such that
\[
\hat{\varphi}'(\alpha_{1}, \dots, \alpha_{m}) = \psi'(\beta).
\]
\end{Claim}

\begin{proof}
Given a tree $T$ and two nodes $w < u$, we denote by $[w, u]$ the set of nodes of $T$ in between $w$ and $u$ (including $w$ and $u$). Bearing this in mind, suppose, with a view to contradiction, that $\hat{\varphi}'(\alpha_{1}, \dots, \alpha_{m}) = \psi'(\beta)$ for some $\alpha_{1}, \dots, \alpha_{m}, \beta \in Fm$. Consequently,
\begin{equation}\label{Eq:tree-same}
T(\hat{\varphi}'(\alpha_{1}, \dots, \alpha_{m})) = T(\psi'(\beta)).
\end{equation}

Notice that, because of $\varphi' = \Box^{2k}\Diamond \varphi(\Box^{k} \Diamond x)$, the tree $T(\varphi)^{-x}$ can be identified with a unique downset of the tree obtained by removing the smallest $2k+1$ nodes from $T(\hat{\varphi}'(\alpha_{1}, \dots, \alpha_{m}))$.\ Working under this identification, a node $w$ of $T(\hat{\varphi}'(\alpha_{1}, \dots, \alpha_{m}))$ belongs to $T(\varphi)^{-x}$ if and only it it satisfies the following conditions:
\benroman
\item\label{Eq:nodes-good1} The downset of $w$ in $T(\hat{\varphi}'(\alpha_{1}, \dots, \alpha_{m}))$ has at least $2k+2$ and at most $3k$ elements;
\item\label{Eq:nodes-good2} There are no nodes $u < v$ in $T(\hat{\varphi}'(\alpha_{1}, \dots, \alpha_{m}))$ such that $w \in [u, v]$ and $[u, v] = \{ b_{1}, \dots, b_{k}, d \}$, where
\[
b_{1} < \dots < b_{k} < d
\]
and the nodes $b_{i}$ and $d$ are labeled, respectively, with $\Box$ and $\Diamond$.
\eroman
This observation is a consequence of the fact that $k$ is greater than the length of branches in $T(\varphi)$ and $\varphi' = \Box^{2k}\Diamond \varphi(\Box^{k}\Diamond x)$.

A similar argument shows that $T(\psi)^{-x}$ can be identified with the subtree of $T(\psi'(\beta))$ comprising the nodes $w$ of $T(\psi'(\beta))$ satisfying conditions (\ref{Eq:nodes-good1}, \ref{Eq:nodes-good2}). Together with (\ref{Eq:tree-same}), this implies $T(\varphi)^{-x} = T(\psi)^{-x}$ and, therefore, $\varphi = \psi$. But this contradicts the assumption that the formulas $\varphi$ and $\psi$ are distinct.
\end{proof}

Our aim is to show that $\vdash$ has a $\btau$-algebraic semantics, where
\[
\btau(x) \coloneqq \{ \varphi' \thickapprox \psi' \}.
\]
Suppose the contrary, with a view to contradiction. By Proposition \ref{Prop:alg-sem-meaning} there is $\Gamma \cup \{ \gamma \} \subseteq Fm$ such that $\btau(\gamma) \subseteq \theta(\Gamma, \btau)$ and $\Gamma \nvdash \gamma$. We can assume without loss of generality $\Gamma = \textup{Cn}_{\vdash}(\Gamma)$.

\begin{Claim}\label{Claim:compatibility-tree-trick}
The congruence $\theta(\Gamma, \btau)$ is compatible with $\Gamma$.
\end{Claim}

\begin{proof}
Consider $\alpha, \beta \in Fm$ such that $\alpha \in \Gamma$ and $\langle \alpha, \beta \rangle \in \theta(\Gamma, \btau)$. We need to prove that $\beta \in \Gamma$. As $\Gamma = \textup{Cn}_{\vdash}(\Gamma)$, it suffices to show $\Gamma \vdash \beta$. Bearing in mind that $\class{M}$ is a matrix semantics for $\vdash$, consider $\langle \A, F \rangle \in \class{M}$ and a homomorphism $h \colon \Fm \to \A$ such that $h[\Gamma] \subseteq F$. From (\ref{Eq:M-proves-equation}) it follows that the generators of $\theta(\Gamma, \btau)$ belong to the kernel $\textup{Ker}(h)$. Consequently, $\langle \alpha, \beta \rangle \in \theta(\Gamma, \btau) \subseteq \textup{Ker}(h)$ and, therefore, $h(\beta) = h(\alpha) \in h[\Gamma] \subseteq F$.
\end{proof}

Since $\btau(\gamma) \subseteq \theta(\Gamma, \btau)$, by Theorem \ref{Thm:Maltsev-lemma} there are 
\begin{align*}
\alpha_{0}, \dots, \alpha_{n} &\in Fm\\
\beta_{0}, \dots \beta_{n-1} &\in \Gamma\\
 \delta_{0}(x, \vec{z}), \dots, \delta_{n-1}(x, \vec{z})  &\in Fm
\end{align*}
such that $\varphi'(\gamma) = \alpha_{0}$, $\psi'(\gamma) = \alpha_{n}$, and
\begin{equation}\label{Eq:Maltsev}
\{ \alpha_{i}, \alpha_{i+1} \} = \{ \delta_{i}(\varphi'(\beta_{i}), \vec{z}), \delta_{i}(\psi'(\beta_{i}), \vec{z}) \}\text{, for every }i < n.
\end{equation}

We shall prove by induction on $i$ that for all $\alpha_{i}$ there are $\gamma_{1}^{i}, \dots, \gamma_{m}^{i} \in Fm$ such that
\begin{equation}\label{Eq:induction}
\alpha_{i} = \hat{\varphi}'(\gamma_{1}^{i}, \dots, \gamma_{m}^{i}) \text{ and }\gamma \equiv \gamma_{1}^{i} \equiv \cdots \equiv \gamma_{m}^{i} \mod \theta(\Gamma, \btau).
\end{equation}
For the base case, it suffices to take
\[
\gamma = \gamma_{1}^{0} = \cdots =  \gamma_{m}^{0}.
\] 
For the induction step, suppose that $\alpha_{i} = \hat{\varphi}'(\gamma_{1}^{i}, \dots, \gamma_{m}^{i})$ for some $\gamma_{1}^{i}, \dots, \gamma_{m}^{i} \in Fm$ such that $\gamma \equiv \gamma_{1}^{i}\equiv  \cdots \equiv \gamma_{m}^{i} \mod \theta(\Gamma, \btau)$. Moreover, from (\ref{Eq:Maltsev}) it follows $\{ \alpha_{i}, \alpha_{i+1} \} = \{ \delta_{i}(\varphi'(\beta_{i}), \vec{z}), \delta_{i}(\psi'(\beta_{i}), \vec{z}) \}$. There are two cases:
\benroman
\item\label{item:Malt1} either $\hat{\varphi}'(\gamma_{1}^{i}, \dots, \gamma_{m}^{i}) = \alpha_{i} = \delta_{i}(\varphi'(\beta_{i}), \vec{z})$ and $\alpha_{i+1} = \delta_{i}(\psi'(\beta_{i}), \vec{z})$, or
\item\label{item:Malt2} $\hat{\varphi}'(\gamma_{1}^{i}, \dots, \gamma_{m}^{i}) = \alpha_{i} = \delta_{i}(\psi'(\beta_{i}), \vec{z})$ and $\alpha_{i+1} = \delta_{i}(\varphi'(\beta_{i}), \vec{z})$.
\eroman
(\ref{item:Malt1}): We shall prove that there exist $\epsilon_{1}(x, \vec{z}), \dots, \epsilon_{m}(x, \vec{z}) \in Fm$ such that
\begin{equation}\label{Eq:what-happened-to-delta}
\delta_{i}(x, \vec{z}) = \hat{\varphi}'(\epsilon_{1}(x, \vec{z}), \dots, \epsilon_{m}(x, \vec{z})).
\end{equation}
Suppose the contrary, with a view to contradiction.\ By (\ref{item:Malt1}), $\hat{\varphi}'(\gamma_{1}^{i}, \dots, \gamma_{m}^{i}) = \delta_{i}(\varphi'(\beta_{i}), \vec{z})$. Thus, $T(\delta_{i}(x, \vec{z}))^{-x}$ and $T(\varphi')^{-x}$ can be identified with two unique (possibly empty) downsets of $T(\hat{\varphi}'(\gamma_{1}^{i}, \dots, \gamma_{m}^{i}) )$. Accordingly, from now on we shall work under this identification. 

As $\delta_{i}(x, \vec{z}) \ne \hat{\varphi}'(\epsilon_{1}, \dots, \epsilon_{m})$ for any $\epsilon_{1}, \dots, \epsilon_{m} \in Fm$, the tree $T(\varphi')^{-x}$ is not a downset of $T(\delta_{i}(x, \vec{z}))^{-x}$. Let then $w$ be a node of $T(\varphi')^{-x}$ that is minimal in $T(\varphi')^{-x} \smallsetminus T(\delta_{i}(x, \vec{z}))^{-x}$. Let also $\chi$ be the unique formula such that $T(\chi)$ is the upset of $w$ in $T(\hat{\varphi}'(\gamma_{1}^{i}, \dots, \gamma_{m}^{i}))$. As $\delta_{i}(\varphi'(\beta_{i}), \vec{z}) = \hat{\varphi}'(\gamma_{1}^{i}, \dots, \gamma_{m}^{i})$, the minimality of $w$ implies that
\begin{equation}\label{Eq:number-of-boxes}
\chi = \varphi'(\beta_{i}) = \Box^{2k}\Diamond \varphi(\Box^{k} \Diamond \beta_{i}).
\end{equation}
As $k$ is greater than the length of all branches in $T(\varphi)$ and $\varphi' = \Box^{2k}\Diamond \varphi(\Box^{k} \Diamond x)$, there cannot be any point  $u$ in $T(\varphi')^{-x}$ other than the root whose upset in $T(\hat{\varphi}'(\gamma_{1}^{i}, \dots, \gamma_{m}^{i}))$ is the subformula tree of a formula starting with the prefix $\Box^{2k}$. Thus we conclude that $w$ is the root of $T(\varphi')^{-x}$ and, therefore, of $T(\hat{\varphi}'(\gamma_{1}^{i}, \dots, \gamma_{m}^{i}))$. As $T(\chi)$ is the upset of the root $w$ of $T(\hat{\varphi}'(\gamma_{1}^{i}, \dots, \gamma_{m}^{i}))$, together with (\ref{Eq:number-of-boxes}), this implies
\begin{equation}\label{Eq:first-trick}
\hat{\varphi}'(\gamma_{1}^{i}, \dots, \gamma_{m}^{i}) = \chi = \varphi'(\beta_{i}).
\end{equation}
Since $\Var(\varphi') = \{ x \}$ by (\ref{Eq:x-appear}), this yields
\[
\beta_{i} = \gamma_{1}^{i} = \cdots = \gamma_{m}^{i}.
\]
As by inductive hypothesis $\gamma \equiv \gamma_{1}^{i}\equiv  \cdots \equiv \gamma_{m}^{i} \mod \theta(\Gamma, \btau)$, the above display guarantees $\langle \gamma, \beta_{i} \rangle \in \theta(\Gamma, \btau)$. Lastly, since $\beta_{i} \in \Gamma$ and $\theta(\Gamma, \btau)$ is compatible with $\Gamma$ by Claim \ref{Claim:compatibility-tree-trick}, this implies $\gamma \in \Gamma$, a contradiction. This establishes that there exist $\epsilon_{1}(x, \vec{z}), \dots, \epsilon_{m}(x, \vec{z}) \in Fm$ validating (\ref{Eq:what-happened-to-delta}).

From (\ref{Eq:what-happened-to-delta}) it follows
\[
\hat{\varphi}'(\gamma_{1}^{i}, \dots, \gamma_{m}^{i}) = \delta_{i}(\varphi'(\beta_{i}), \vec{z}) = \hat{\varphi}'(\epsilon_{1}(\varphi'(\beta_{i}), \vec{z}), \dots, \epsilon_{m}(\varphi'(\beta_{i}), \vec{z})).
\]
Consequently,
\[
\gamma_{1}^{i} = \epsilon_{1}(\varphi'(\beta_{i}), \vec{z}), \dots, \gamma_{m}^{i} = \epsilon_{m}(\varphi'(\beta_{i}), \vec{z}).
\]
Thus, defining
\[
\gamma_{1}^{i+1} \coloneqq \epsilon_{1}(\psi'(\beta_{i}), \vec{z}), \dots, \gamma_{m}^{i+1} \coloneqq \epsilon_{m}(\psi'(\beta_{i}), \vec{z}),
\]
we obtain
\begin{align*}
\alpha_{i+1} = \delta_{i}(\psi'(\beta_{i}), \vec{z}) &= \hat{\varphi}'(\epsilon_{1}(\psi'(\beta_{i}), \vec{z}), \dots, \epsilon_{m}(\psi'(\beta_{i}), \vec{z}))\\
&= \hat{\varphi}'(\gamma_{1}^{i+1}, \dots, \gamma_{m}^{i+1}).
\end{align*}
Furthermore, for every $m \geq j \in \omega$,
\[
\gamma \equiv \gamma_{j}^{i} = \epsilon_{j}(\varphi'(\beta_{i}), \vec{z}) \equiv \epsilon_{j}(\psi'(\beta_{i}), \vec{z}) = \gamma_{j}^{i+1} \mod \theta(\Gamma, \btau).
\]
The first equivalence above follows from the inductive hypothesis, while the third from $\beta_{i} \in \Gamma$. This concludes the proof of the induction step in case (\ref{item:Malt1}).

(\ref{item:Malt2}): Replicating the proof described for case (\ref{item:Malt1}) up to (\ref{Eq:first-trick}), we obtain
\[
\hat{\varphi}'(\gamma_{1}^{i}, \dots, \gamma_{m}^{i}) = \psi'(\beta_{i}),
\]
a contradiction with Claim \ref{Claim:first}. This concludes the inductive proof.

Finally, taking $i=n$ in (\ref{Eq:induction}), we obtain
\[
\hat{\varphi}'(\gamma_{1}^{n}, \dots, \gamma_{m}^{n}) = \alpha_{n} =  \psi'(\gamma).
\]
But this is impossible by Claim \ref{Claim:first}. Hence we reached the desired contradiction. We conclude that $\vdash$ has a $\btau$-algebraic semantics.

This establishes the theorem for logics whose language comprises two distinct unary symbols. To conclude the proof, suppose that $\LL_{\vdash}$ comprises an $n$-ary connective $f$ with $ n \geq 2$. Then define
\[
\Box x \coloneqq f(f(x, \dots, x), x, \dots, x) \text{ and }\Diamond x \coloneqq f(x, f(x, \dots, x), x \dots, x).
\]
Using this definition of $\Box$ and $\Diamond$, it is possible to reproduce the argument detailed in the case where $\LL_{\vdash}$ comprises two distinct unary connectives.
\end{proof}

\begin{Remark}
Interesting conditions implying the existence of an algebraic semantics were first discovered in \cite[Thms.\ 3.1 \& 3.8]{BlRe03}. In particular, \cite[Thm.\ 3.8]{BlRe03} is a consequence of Theorem \ref{Thm:sufficient}, while \cite[Thm.\ 3.1]{BlRe03} can be derived from a simple variant of the argument of detailed above.
\qed
\end{Remark}

\section{Examples: modal and substructural logics}

From Theorem \ref{Thm:sufficient} it follows that most well-known logics have an algebraic semantics, as we proceed to explain.

\begin{exa}[\textsf{Local modal logics}]
Given a variety $\class{K}$ of modal algebras, let $\vdash^{\ell}_{\class{K}}$ be the logic defined by rule
\begin{align*}
\Gamma \vdash^{\ell}_{\class{K}}\varphi \Longleftrightarrow& \text{ for all }\A \in \class{K}, a \in A, \text{and homomorphism }h \colon \Fm \to \A,\\
&\text{ if }a \leq h(\gamma)\text{ for all }\gamma \in \Gamma \text{, then }a \leq h(\varphi).
\end{align*}
Logics arising in this way have been called \textit{local modal logics} \cite{MK07c}.
\qed
\end{exa}

\begin{exa}[\textsf{Substructural logics}]
A \textit{commutative FL-algebra} is a structure $\A = \langle A; \land, \lor, \cdot, \to, 1, 0 \rangle$ comprising a lattice $\langle A; \land, \lor \rangle$, an Abelian monoid $\langle A; \cdot, 1 \rangle$, a constant $0$, and a binary operation $\to$ such that for every $a, b, c \in A$,
\[
a \cdot b \leq c \Longleftrightarrow a \leq b \to c.
\]
Notice that if $\cdot$ coincides with $\land$ and $0$ is the minimum of $\leq$, then $\A$ is a Heyting algebra.

The logic $\vdash_{\class{K}}$ associated with a variety $\class{K}$ of residuated lattices is defined by the rule
\[
\Gamma \vdash_{\class{K}} \varphi \Longleftrightarrow \btau[\Gamma] \vDash_{\class{K}} \btau(\varphi),
\]
where $\btau = \{ x \land 1 \thickapprox 1 \}$. Logics arising in this way have been called \textit{substructural logics with exchange} \cite{GaJiKoOn07}.\ Notice that when $\class{K}$ is a variety of Heyting algebras, $\vdash_{\class{K}}$ is the intermediate logic associated with $\class{K}$.
\qed
\end{exa}



An equation $\varphi \thickapprox \psi$ is said to be \textit{nontrivial} if $\varphi \ne \psi$.

\begin{Theorem}\label{Thm:exa} The following holds:
\benroman
\item\label{Eq:exa-1} Extensions of fragments of local and global modal logics in which a connective among $\{ \land, \lor, \to \}$ is term-definable have an algebraic semantics.
\item\label{Eq:exa-2} Extensions of fragments of substructural logics with exchange in which a connective among  $\{ \land, \lor, \cdot, \to \}$ is term-definable have an algebraic semantics.
\eroman
\end{Theorem}

\begin{proof}
We detail the proof of (\ref{Eq:exa-2}) only, as that of (\ref{Eq:exa-1}) is analogous. By Corollary \ref{Cor:persist-extensions} it suffices to show that fragments of substructural logics with exchange in which a connective among $\{ \land, \lor, \cdot, \to \}$ is term-definable have an algebraic semantics. Accordingly, let $\vdash$ be one of these fragments. Then there is a variety $\class{K}$ of commutative FL-algebras such that $\vdash$ is the logic induced by the class of matrices
\[
\class{M} \coloneqq \{ \langle \A, F \rangle \colon \A \text{ is the $\LL_{\vdash}$-reduct of an FL-algebra $\A^{+}$ and }F = \btau(\A^{+}) \},
\]
where $\btau = \{ x \land 1 \thickapprox 1 \}$. Thus, by Lemma \ref{Lem:matrix-to-rules-1} and Theorem \ref{Thm:sufficient}, it only remains to prove that $\class{M}$ validates a nontrivial equation $\epsilon \thickapprox \delta$ such that $\Var(\epsilon) \cup \Var(\delta) = \{ x \}$. To this end, recall that in $\vdash$ a connective in $\ast \in \{ \land, \lor, \cdot, \to \}$ is term-definable by means of a, possibly complex, formula $x + y$.
\benroman
\item If $\ast \in \{ \land, \lor \}$, take $\epsilon \coloneqq x$ and $\delta \coloneqq x + x$;
\item If $\ast = \cdot$, take $\epsilon \coloneqq x + (x + x)$ and $\delta \coloneqq (x + x) + x$;
\item If $\ast = {\to}$, take $\epsilon \coloneqq (x + x) + (x + x)$ and $\delta \coloneqq x + ((x + x) + x)$.
\eroman
In all cases, $\epsilon \thickapprox \delta$ is nontrivial, valid in $\class{M}$, and such that $\Var(\epsilon) \cup \Var(\delta) = \{ x \}$.
\end{proof}

\begin{Remark}
Theorem \ref{Thm:exa} was essentially discovered \cite[Cor.\ 3.4, 3.7, 3.9]{BlRe03}. With respect to its first incarnation, the current formulation provides a small improvement by admitting the presence of the operation $\cdot$ in condition (\ref{Eq:exa-2}) and, more in general, by allowing the operations $\land, \lor, \to, \cdot$ to be term-definable, as opposed to basic connectives. 
\qed
\end{Remark}

\section{Logics with theorems}

The aim of this section is to establish the following characterization of logics with theorems possessing an algebraic semantics.

\begin{Theorem}\label{Thm:with-thms}
Let $\vdash$ be a nontrivial logic with a theorem $\varphi$ such that $\Var(\varphi) \ne \emptyset$. The following conditions are equivalent:
\benroman
\item\label{Eq:with-thms-1} $\vdash$ has an algebraic semantics;
\item\label{Eq:with-thms-3} Either $\vdash$ is assertional and graph-based or it is not graph-based and there are two distinct logically equivalent formulas $\epsilon$ and $\delta$ such that $\Var(\epsilon) \cup \Var(\delta) = \{ x \}$;
\item\label{Eq:with-thms-1b} Either $\vdash$ is assertional and graph-based or it is not graph-based and $\Alg(\vdash) \vDash \epsilon \thickapprox \delta$ for some nontrivial equation $\epsilon \thickapprox \delta$ such that $\Var(\epsilon) \cup \Var(\delta) = \{ x \}$;
\item\label{Eq:with-thms-2} Either $\vdash$ is assertional and graph-based or it is not graph-based and has a matrix semantics validating a nontrivial equation $\epsilon \thickapprox \delta$ such that $\Var(\epsilon) \cup \Var(\delta) = \{ x \}$.
\eroman
\end{Theorem}

 We restricted the statement of Theorem \ref{Thm:with-thms} to the case of nontrivial logics for the sake of readability. However, the reader may consult Proposition \ref{Prop:trivial} for a catalogue of trivial logics with an algebraic semantics. 


\begin{Notation}\emph{From now on we shall denote by $\Box$ the unary connective (if any) of an arbitrary graph-based logic and its set of constants by $\{ \textbf{c}_{i} \colon i < \alpha \}$ where $\alpha$ is a suitable ordinal.}
\end{Notation}

\begin{Proposition}\label{Prop:unital}
Graph-based logics with an algebraic semantics have a unital matrix semantics.
\end{Proposition}

\begin{proof}
Let $\vdash$ be a graph-based logic with an algebraic semantics. By Proposition \ref{Prop:unital-rules}, to prove that $\vdash$ has a unital matrix semantics, it suffices to show that $x, y, \varphi(x, \vec{z}) \vdash \varphi(y, \vec{z})$ for all $\varphi(v, \vec{z}) \in Fm$. Since $\vdash$ is graph-based, if it does not have unary connectives, this condition holds vacuously. Then we can assume that $\vdash$ has a unary connective $\Box$. Since $\vdash$ is graph-based, it will be enough to show that for every $n \in \omega$,
\[
x, y, \Box^{n}x \vdash \Box^{n}y.
\]

To this end, consider $n \in \omega$. By assumption $\vdash$ has an algebraic semantics given by a class of algebras $\class{K}$ and a set of equations $\btau(x)$. Consider $\A \in \class{K}$ and $a, b \in A$ such that
\[
a, b, \underbrace{\Box^{\A} \dots \Box^{\A}}_{n\text{-times}}a \in \btau(\A).
\]
Then consider $\epsilon \thickapprox \delta \in \btau$. We need to prove
\begin{equation}\label{Eq:zsdrftgyhujikolp}
\epsilon^{\A}(\underbrace{\Box^{\A} \dots \Box^{\A}}_{n\text{-times}}b) = \delta^{\A}(\underbrace{\Box^{\A} \dots \Box^{\A}}_{n\text{-times}}b). 
\end{equation}
As $\vdash$ is graph-based, one of the following conditions holds:
\benroman
\item \label{item-equ-1}$\epsilon = \Box^{m}x$ and $\delta = \Box^{k}x$ for some $m, k \in \omega$;
\item \label{item-equ-2}$\epsilon = \Box^{m}\textbf{c}_{i}$ and $\delta = \Box^{k}x$ for some $m, k \in \omega$ and $i < \alpha$;
\item \label{item-equ-3}$\epsilon = \Box^{m}x$ and $\delta = \Box^{k}\textbf{c}_{i}$ for some $m, k \in \omega$ and $i < \alpha$;
\item \label{item-equ-4}$\epsilon = \Box^{m} \textbf{c}_{i}$ and $\delta = \Box^{k}\textbf{c}_{j}$ for some $m, k \in \omega$ and $i, j < \alpha$.
\eroman

(\ref{item-equ-1}): From $b \in \btau(\A)$ and $\epsilon \thickapprox \delta \in \btau$, it follows
\[
\underbrace{\Box^{\A} \dots \Box^{\A}}_{m\text{-times}}b = \underbrace{\Box^{\A} \dots \Box^{\A}}_{k\text{-times}}b.
\]
In particular, this implies
\[
\epsilon^{\A}(\underbrace{\Box^{\A} \dots \Box^{\A}}_{n\text{-times}}b) = \underbrace{\Box^{\A} \dots \Box^{\A}}_{m+n\text{-times}}b = \underbrace{\Box^{\A} \dots \Box^{\A}}_{k+n\text{-times}}b = \delta^{\A}(\underbrace{\Box^{\A} \dots \Box^{\A}}_{n\text{-times}}b).
\]

(\ref{item-equ-2}): From $a, \underbrace{\Box^{\A} \dots \Box^{\A}}_{n\text{-times}}a \in \btau(\A)$ and $\epsilon \thickapprox \delta \in \btau$ it follows
\[
\underbrace{\Box^{\A} \dots \Box^{\A}}_{n+m\text{-times}}\textbf{c}_{i}^{\A} = \underbrace{\Box^{\A} \dots \Box^{\A}}_{n+k\text{-times}}a = \underbrace{\Box^{\A} \dots \Box^{\A}}_{k+n\text{-times}}a = \underbrace{\Box^{\A} \dots \Box^{\A}}_{m\text{-times}}\textbf{c}_{i}^{\A}.
\]
Moreover, from $b \in \btau(\A)$ and $\epsilon \thickapprox \delta \in \btau$ it follows
\[
\underbrace{\Box^{\A} \dots \Box^{\A}}_{m\text{-times}} \textbf{c}_{i}^{\A} = \underbrace{\Box^{\A} \dots \Box^{\A}}_{k\text{-times}}b.
\]
From the two displays above we obtain
\[
\epsilon^{\A}(\underbrace{\Box^{\A} \dots \Box^{\A}}_{n\text{-times}}b) = \underbrace{\Box^{\A} \dots \Box^{\A}}_{m\text{-times}} \textbf{c}_{i}^{\A} = \underbrace{\Box^{\A} \dots \Box^{\A}}_{n+m\text{-times}}\textbf{c}_{i}^{\A} = \underbrace{\Box^{\A} \dots \Box^{\A}}_{k+n\text{-times}}b = \delta^{\A}(\underbrace{\Box^{\A} \dots \Box^{\A}}_{n\text{-times}}b).
\]
Case (\ref{item-equ-3}) is handled analogously to case (\ref{item-equ-2}). 

(\ref{item-equ-4}): As $\epsilon$ and $\delta$ are closed formulas, the fact that $\epsilon^{\A}(a) = \delta^{\A}(a)$ implies 
\[
\epsilon^{\A}(\underbrace{\Box^{\A} \dots \Box^{\A}}_{n\text{-times}}b) = \delta^{\A}(\underbrace{\Box^{\A} \dots \Box^{\A}}_{n\text{-times}}b).
\]

This establishes (\ref{Eq:zsdrftgyhujikolp}) and, therefore, that $\underbrace{\Box^{\A} \dots \Box^{\A}}_{n\text{-times}}b \in \btau(\A)$. As $\class{K}$ is a $\btau$-algebraic semantics for $\vdash$, we conclude that $x, y, \Box^{n}x \vdash \Box^{n}y$, as desired.
\end{proof}

\begin{Corollary}\label{Cor:graph-based-with-theorems}
A graph-based logic with theorems has an algebraic semantics if and only if it is assertional.
\end{Corollary}

\begin{proof}
Assertional logics have an algebraic semantics by Proposition \ref{Prop:assertional}. Then consider a logic $\vdash$ that is graph-based and has an algebraic semantics. By Proposition \ref{Prop:unital} $\vdash$ has a unital matrix semantics. Therefore, if $\vdash$ has theorems, it is also assertional.
\end{proof}

We shall also rely on the following technical observation:

\begin{Lemma}\label{Lem:technical-x-in-tau}
Let $\vdash$ be a logic with a $\btau$-algebraic semantics. If $\vdash$ is nontrivial, then $\btau$ contains a nontrivial equation $\epsilon \thickapprox \delta$ such that $\Var(\epsilon) \cup \Var(\delta) = \{ x \}$.
\end{Lemma}

\begin{proof}
We reason by contraposition. Suppose that $\vdash$ has a $\btau$-algebraic semantics $\class{K}$ and that $\btau$ contains no nontrivial equation $\epsilon \thickapprox \delta$ such that $\Var(\epsilon) \cup \Var(\delta) = \{ x \}$. We shall prove that $x \vdash y$, i.e., that $\vdash$ is trivial. As $\class{K}$ is a $\btau$-algebraic semantics for $\vdash$, it suffices to show that for all $\A \in \class{K}$ and $a, c \in A$, if $a \in \btau(\A)$, then $c \in \btau(\A)$. 

To this end, consider $\A \in \class{K}$ and $a, c \in A$ such that $a \in \btau(\A)$. We need to show that $\epsilon^{\A}(c) = \delta^{\A}(c)$ for all $\epsilon \thickapprox \delta \in \btau$. Then consider an equation $\epsilon \thickapprox \delta$ in $\btau$. By the assumption, we have two cases: either $\epsilon \thickapprox \delta$ is trivial or $\Var(\epsilon) \cup \Var(\delta) = \emptyset$. If $\epsilon \thickapprox \delta$ is trivial, then $\epsilon = \delta$, whence $\epsilon^{\A}(c) = \delta^{\A}(c)$. Then we consider the case where $\Var(\epsilon) \cup \Var(\delta) = \emptyset$, i.e., $\epsilon$ and $\delta$ are closed formulas. In this case, $a \in \btau(\A)$ implies $\epsilon^{\A}(a) = \delta^{\A}(a)$. But, since $\epsilon$ and $\delta$ are closed formulas, $\epsilon^{\A}(a) = \epsilon^{\A}(c)$ and $\delta^{\A}(a) = \delta^{\A}(c)$, whence $\epsilon^{\A}(c) = \delta^{\A}(c)$. We conclude that $c \in \btau(\A)$, as desired.
\end{proof}

We are now ready to prove the main result of this section:

\begin{proof}[Proof of Theorem \ref{Thm:with-thms}]
Let $\vdash$ be a nontrivial logic with a theorem $\varphi$ such that $\Var(\varphi) \ne \emptyset$. As $\vdash$ is substitution invariant, we can assume without loss of generality $\Var(\varphi) = \{ x \}$.

(\ref{Eq:with-thms-1})$\Rightarrow$(\ref{Eq:with-thms-2}): Suppose that $\vdash$ has a $\btau$-algebraic semantics $\class{K}$. First if $\vdash$ is graph-based, then it is assertional by Corollary \ref{Cor:graph-based-with-theorems}. Then we consider the case where $\vdash$ is not graph-based. 

As $\vdash$ is nontrivial, we can apply Lemma \ref{Lem:technical-x-in-tau}, obtaining that $\btau$ contains a nontrivial equation $\epsilon \thickapprox \delta$ such that $\Var(\epsilon) \cup \Var(\delta) = \{ x \}$. Recall from Proposition \ref{Prop:matrix-sem-alg-sem} that
\[
\class{M} \coloneqq \{ \langle \A, \btau(\A) \rangle \colon \A \in \class{K} \}
\]
is a matrix semantics for $\vdash$. As $\epsilon \thickapprox \delta \in \btau$ and $\varphi$ is a theorem of $\vdash$,
\[
\class{M} \vDash \epsilon(\varphi) \thickapprox \delta(\varphi).
\]
Moreover, as $\Var(\epsilon) \cup \Var(\delta) = \{ x \} = \Var(\varphi)$,
\[
\Var(\epsilon(\varphi)) \cup \Var(\delta(\varphi)) = \{ x \}. 
\]

In view of the two displays above, since $\class{M}$ is a matrix semantics for $\vdash$, it only remains to prove that the equation $\epsilon(\varphi) \thickapprox \delta(\varphi)$ is nontrivial. Suppose the contrary, with a view to contradiction. As $\Var(\varphi) = \{ x \}$, the tree $T(\epsilon)^{-x}$ can be identified with the subtree of $T(\epsilon(\varphi))$ obtained by removing from $T(\epsilon(\varphi))$ the upsets equal to $T(\varphi)$. Similarly, $T(\delta)^{-x}$ can be identified with the subtree of $T(\delta(\varphi))$ obtained by removing from $T(\delta(\varphi))$ the upsets equal to $T(\varphi)$. Since $\epsilon(\varphi) = \delta(\varphi)$, this implies $T(\epsilon)^{-x} = T(\delta)^{-x}$ and, therefore, $\epsilon = \delta$. But this contradicts the fact that the equation $\epsilon \thickapprox \delta$ is nontrivial.

(\ref{Eq:with-thms-3})$\Rightarrow$(\ref{Eq:with-thms-1}): If $\vdash$ is graph-based, then by assumption it is also assertional. Thus, in this case, $\vdash$ has an algebraic semantics by Proposition \ref{Prop:assertional}. Then we consider the case where $\vdash$ is not graph-based. Together with the assumption, this implies that $\vdash$ has an algebraic semantics by Theorem \ref{Thm:sufficient}.

Lastly, conditions (\ref{Eq:with-thms-3}), (\ref{Eq:with-thms-1b}), and (\ref{Eq:with-thms-2}) are equivalent by Lemma \ref{Lem:matrix-to-rules-1}.
\end{proof}

\begin{problem}
Extend the characterization of logics with an algebraic semantics given in Theorem \ref{Thm:with-thms} beyond the setting of logics with theorems.
\end{problem}

\section{Protoalgebraic logics}

In this section, Theorem \ref{Thm:with-thms} is specialized to the class of protoalgebraic logics \cite{Cz01}, introduced by Czelakowski \cite{Cz85,Cz86} and refined by Blok and Pigozzi \cite{BP86}.

\begin{law}
A logic $\vdash$ is \textit{protoalgebraic} if there exists a set of formulas $\Delta(x, y)$ such that
\begin{equation}\label{Eq:def-protoalgebraic}
\emptyset \vdash \Delta(x, x) \text{ and } x, \Delta(x, y) \vdash y.
\end{equation}
\end{law}


The class of protoalgebraic logics embraces all logics $\vdash$ possessing a term-definable operation $x \to y$ such that
\[
\emptyset \vdash x \to x\text{ and }x, x \to y \vdash y
\]
as, in this case, the set $\Delta(x, y) \coloneqq \{ x \to y \}$ satisfies condition (\ref{Eq:def-protoalgebraic}).\ In particular, protoalgebraic logics comprise all fragments of substructural logics and of local and global modal logics containing $\to$. On the other hand, examples of nonprotoalgebraic logics include Visser's logic \cite{Vi81}, as shown in  \cite[Thm.\ 14]{SuWoZa98} (see also \cite{CeJa01}), as well as many implication-less fragments of familiar logics such as positive modal logic \cite{Du95}, as proved in \cite[Thm.\ 9]{Ja02}.

\begin{Lemma}[\protect{\cite[Prop.\ 6.11.(5, 6, 7)]{AAL-AIT-f}}]\label{Lem:proto-facts}
If $\vdash$ is a nontrivial protoalgebraic logic, then $\vdash$ has an $n$-ary connective with $n \geq 2$ and a theorem $\varphi$ such that $\Var(\varphi) \ne \emptyset$.
\end{Lemma}



\begin{Theorem}\label{Thm:protoalgebraic-logics}
The following conditions are equivalent for a nontrivial protoalgebraic logic $\vdash$:
\benroman
\item\label{item:proto-1} $\vdash$ has an algebraic semantics;
\item\label{item:proto-3} There are two distinct logically equivalent formulas;
\item\label{item:proto-1b} $\Alg(\vdash)$ validates a nontrivial equation;
\item\label{item:proto-2} $\vdash$ has a matrix semantics validating a nontrivial equation;
\item\label{item:proto-4} The free Lindenbaum-Tarski algebra $\Fm(\vdash)$ is not isomorphic to $\Fm$.
\eroman
\end{Theorem}

%
%
%
%

\begin{proof}
Let $\vdash$ be a nontrivial protoalgebraic logic. By Lemma \ref{Lem:proto-facts}, $\vdash$ has a theorem $\varphi$ such that $\Var(\varphi) \ne \emptyset$ and an $n$-ary connective $f(x_{1}, \dots, x_{n})$ with $n \geq 2$.

(\ref{item:proto-1})$\Rightarrow$(\ref{item:proto-2}): Suppose that $\vdash$ has an algebraic semantics. As $\vdash$ has a theorem $\varphi$ such that $\Var(\varphi) \ne \emptyset$, it falls in the scope of Theorem \ref{Thm:with-thms}. Therefore, the implication (\ref{Eq:with-thms-1})$\Rightarrow$(\ref{Eq:with-thms-2}) in Theorem \ref{Thm:with-thms} yields that either $\vdash$ is graph-based or it has a matrix semantics validating a nontrivial equation.\ Thus, to conclude the proof, it suffices to show that $\vdash$ is not graph-based. But this follows from the fact that $\vdash$ has an $n$-ary connective $f(x_{1}, \dots, x_{n})$ with $n \geq 2$.

(\ref{item:proto-2})$\Rightarrow$(\ref{item:proto-1}): Suppose that $\vdash$ has a matrix semantics $\class{M}$ validating a nontrivial equation $\epsilon(x_{1}, \dots, x_{n}) \thickapprox \delta(x_{1}, \dots, x_{n})$. By the implication (\ref{Eq:with-thms-2})$\Rightarrow$(\ref{Eq:with-thms-1}) in Theorem \ref{Thm:with-thms}, to conclude the proof it suffices to show that there are formulas $\epsilon'$ and $\delta'$ such that $\Var(\epsilon') \cup \Var(\delta') = \{ x \}$ and the equation $\epsilon' \thickapprox \delta'$ is both nontrivial and valid in $\class{M}$.

To this end, recall that $f$ is $n$-ary with $n \geq 2$. Therefore, if $\epsilon$ and $\delta$ are closed formulas, we are done taking
\[
\epsilon' \coloneqq f(\epsilon, x, \dots, x) \text{ and }\delta' \coloneqq f(\delta, x, \dots, x).
\]

Then we can assume by symmetry that $x_{1} \in \Var(\epsilon)$. For every $n \in \omega$, we define recursively a formula $\varphi_{n}(x)$ by the rule
\[
\varphi_{0}\coloneqq f(x, \dots, x) \text{ and }\varphi_{m+1} \coloneqq f(\varphi_{m}, x, \dots, x).
\]
Clearly $\Var(\varphi_{n}) = \{ x \}$ for all $n \in \omega$. Let then $k \in \omega$ be greater than the length of all branches in the trees $T(\epsilon)$ and $T(\delta)$, and define
\[
\epsilon' \coloneqq \epsilon(\varphi_{k}, \varphi_{2k}, \dots, \varphi_{nk}) \text{ and }\delta' \coloneqq \epsilon(\varphi_{k}, \varphi_{2k}, \dots, \varphi_{nk}).
\]
Notice that $\Var(\epsilon') \cup \Var(\delta') = \{ x \}$ as $x_{1} \in \Var(\epsilon)$ and $\Var(\varphi_{m}) = \{ x \}$ for all $m \in \omega$. Furthermore, $\epsilon' \thickapprox \delta'$ is valid in $\class{M}$, because $\epsilon \thickapprox \delta$ is. 

Therefore, it only remains to prove that the equation $\epsilon' \thickapprox \delta'$ is nontrivial. To this end, observe that, since $k$ is greater than all branches in $T(\epsilon)$, the tree $T(\epsilon)$ can be obtained as follows. First, let $S_{n}$ be the tree obtained by replacing in $T(\epsilon)$ the principal upsets whose longest branch has length $nk + 1$ by a point labeled by $x_{n}$. Then let $S_{n-1}$ be the tree obtained by replacing in $S_{n}$ the principal upsets whose longest branch have length $(n-1)k + 1$ by a point labeled by $x_{n-1}$. Iterating this process, we eventually get a tree $S_{1}$ which coincides with $T(\epsilon)$. Notice that $T(\delta)$ can be obtained by pruning $T(\delta')$ in exactly the same way. Consequently, if $T(\epsilon') = T(\delta')$, then $T(\epsilon) = T(\delta)$ and, therefore, $\epsilon = \delta$. As by assumption $\epsilon \ne \delta$, we conclude that $T(\epsilon') \ne T(\delta')$, i.e., $\epsilon' \ne \delta'$.

(\ref{item:proto-3})$\Leftrightarrow$(\ref{item:proto-4}): It is easy to see that $\Fm(\vdash)$ is isomorphic to $\Fm$ if and only if $\equiv_{\vdash}$ is the identity relation on $Fm$. Bearing this in mind, the equivalence between conditions (\ref{item:proto-3}) and (\ref{item:proto-4}) follows from the definition of $\equiv_{\vdash}$.

Lastly, conditions (\ref{item:proto-3}), (\ref{item:proto-1b}), and (\ref{item:proto-2}) are equivalent by Lemma \ref{Lem:matrix-to-rules-1}. 
\end{proof}

\begin{Remark}
The above proof shows that Theorem \ref{Thm:protoalgebraic-logics} can be generalized to nontrivial logics with a theorem $\varphi$ such that $\Var(\varphi) \ne \emptyset$ and an $n$-ary connective with $n \geq 2$.
\qed
\end{Remark}

By Theorem \ref{Thm:protoalgebraic-logics}, a nontrivial protoalgebraic logic has an algebraic semantics if and only if its algebraic counterpart satisfies a nontrivial equation. It follows that almost all reasonable protoalgebraic logics have an algebraic semantics. Because of this, it makes sense to review some exceptions to this general rule.

\begin{exa}
To the best of our knowledge, the first example of a protoalgebraic logic lacking an algebraic semantics was discovered in \cite[Thm.\ 2.19]{BlRe03} (see also \cite[Prop.\ 5.2]{Font13}) and subsequently generalized as follows. Let $\mathscr{L}$ be an algebraic language comprising at least an $n$-ary connective with $n \geq 2$. A set $\Delta(x, y) \subseteq Fm_{\LL}$ is said to be \textit{coherent} if $\varphi(x, x) = \psi(x, x)$, for all $\varphi, \psi \in \Delta(x, y)$. Then, given a coherent set $\Delta(x, y) \subseteq Fm_{\LL}$, we denote by ${\bf I}(\Delta)$ the logic, formulated in $\mathscr{L}$, axiomatized by the rules 
\[
\emptyset \rhd \Delta(x, x) \text{ and }x, \Delta(x, y) \rhd y.
\]
The logic ${\bf I}(\Delta)$ is protoalgebraic and lacks an algebraic semantics \cite[Prop.\ 5.5]{JMF14opL}.

Another protoalgebraic logic without an algebraic semantics comes from the field of relevance logic. The logic ${\bf P}$--${\bf W}$ \cite{Komori94,MaMe82} is formulated in the language comprising a binary symbol $\to$ only and is axiomatized by the rules
\begin{align*}
\emptyset \rhd (x \to y) &\to ((z \to x) \to (z \to y))\\
\emptyset \rhd (x \to y) &\to (( y \to z) \to (x \to z))\\
\emptyset \rhd x &\to x\\
x, x &\to y \rhd y.
\end{align*}
Observe that the last two rules in the above display guarantee that ${\bf P}$--${\bf W}$ is protoalgebraic. The fact that ${\bf P}$--${\bf W}$ lacks an algebraic semantics was established in \cite[Prop.\ 38]{Ra06a}.
\qed
\end{exa}

As we mentioned, most familiar protoalgebraic logics have an algebraic semantics. However, these algebraic semantics were constructed \textit{ad hoc} and have little to do with the intuitive interpretation of these logics. For this reason, it is natural to wonder whether protoalgebraic logics $\vdash$ with an algebraic semantics have also a standard one. We close this section by showing that this is not the case, even in the well-behaved setting of local modal logics.

To this end, we assume that the reader is familiar with the basics of modal logic \cite{ChZa97,Kr99}. Given a class $\mathfrak{F}$ of Kripke frames, let $\vdash^{\ell}_{\mathfrak{F}}$ be the logic defined by the rule
\begin{align*}
\Gamma \vdash^{\ell}_{\mathfrak{F}} \varphi \Longleftrightarrow& \text{ for every }\langle W, R \rangle \in \mathfrak{F}, w \in W, \text{and evaluation }v  \text{ into }\langle W, R \rangle,\\
&\text{ if }w, v \Vdash \Gamma \text{, then }w, v \Vdash \varphi.
\end{align*}
As it happens, $\vdash^{\ell}_{\mathfrak{F}}$ is always a protoalgebraic logic.

Furthermore, given a Kripke frame $\mathcal{F} = \langle W, R \rangle$, we denote by $\mathcal{F}^{+}$ the Kripke frame obtained adding to $\mathcal{F}$ a new point $w^{+}$ that is related to everything in $W \cup \{ w^{+} \}$. 


\begin{Proposition}\label{Prop:local-moda-logic-not-as}
Let $\mathfrak{F}$ be a nonempty class of Kripke frames such that if $\mathcal{F} \in \mathfrak{F}$, then $\mathcal{F}^{+} \in \mathfrak{F}$. Then $\vdash^{\ell}_{\mathfrak{F}}$ does not have a standard algebraic semantics.
\end{Proposition}

\begin{proof}
We will make use of the following well-known fact: for all $\varphi, \psi \in Fm$,
\[
\Alg(\vdash^{\ell}_{\mathfrak{F}}) \vDash \varphi \thickapprox \psi \Longleftrightarrow \varphi \sineq^{\ell}_{\mathfrak{F}} \psi.
\]

Suppose, with a view to contradiction, that $\vdash^{\ell}_{\mathfrak{F}}$ has a $\btau$-algebraic semantics based on $\Alg(\vdash^{\ell}_{\mathfrak{F}})$. As $\mathfrak{F}$ is nonempty, the logic $\vdash^{\ell}_{\mathfrak{F}}$ is nontrivial. Together with the fact that $\Alg(\vdash^{\ell}_{\mathfrak{F}})$ is a $\btau$-algebraic semantics, this implies that there exists an equation $\epsilon \thickapprox \delta \in \btau$ such that $\Alg(\vdash^{\ell}_{\mathfrak{F}}) \nvDash \epsilon \thickapprox \delta$. Thus, in view of the above display, we can assume by symmetry that $\epsilon \nvdash^{\ell}_{\mathfrak{F}} \delta$. This means that there are $\mathcal{F} \in \mathfrak{F}$, a world $w \in \mathcal{F}$, and an evaluation $v$ in $\mathcal{F}$ such that $w, v \Vdash \epsilon$ and $w, v \nVdash \delta$.

Recall that, by the assumptions, $\mathcal{F}^{+} \in \mathfrak{F}$. Let $v^{+}$ be the unique evaluation on $\mathcal{F}^{+}$ such that for every $y \in \Var$ and $q \in \mathcal{F}^{+}$:
\[
q, v^{+} \Vdash y \Longleftrightarrow \text{ either }(q \in \mathcal{F} \text{ and }q, v \Vdash y)\text{ or }q = w^{+}.
\]
From the definition of $\mathcal{F}^{+}$ and $v^{+}$ it follows that
\[
q, v^{+} \Vdash \varphi \Longleftrightarrow q, v \Vdash \varphi
\]
for all $\varphi \in Fm$ and $q \in \mathcal{F}$. Consequently, as $w, v \Vdash \epsilon$ and $w, v \nVdash \delta$,
\[
w^{+}, v^{+} \Vdash x \text{ and } w^{+}, v^{+} \nVdash \Box(\epsilon \to \delta).
\]
This implies
\[
x \nvdash^{\ell}_{\mathfrak{F}} \Box(\epsilon \to \delta).
\]
On the other hand, clearly $\emptyset \vdash^{\ell}_{\mathfrak{F}} \Box(\delta \to \delta)$. Consequently,
\[
x, \Box(\delta \to \delta) \nvdash^{\ell}_{\mathfrak{F}} \Box(\epsilon \to \delta).
\]
Together with Proposition \ref{Prop:Suszko-congruence} and the fact that $\vdash$ has a $\btau$-algebraic semantics, this implies $\epsilon \thickapprox \delta \notin \btau$, a contradiction.
\end{proof}

\begin{Corollary}\label{Cor:modal-not-an-algebraic-semantics}
If $\mathfrak{F}$ is the class of all (resp.\ transitive, resp.\ reflexive and transitive) Kripke frames, then $\vdash^{\ell}_{\mathfrak{F}}$ does not have a standard algebraic semantics, namely one based on the variety of modal algebras (resp.\ of K4-algebras, resp.\ of interior algebras).
\end{Corollary}

\begin{proof}
We detail the case where $\mathfrak{F}$ is the class of all Kripke frames, as the remaining ones are analogous. First observe that if $\mathcal{F} \in \mathfrak{F}$, then $\mathcal{F}^{+} \in \mathfrak{F}$. Hence we can apply Proposition \ref{Prop:local-moda-logic-not-as}, obtaining that $\vdash^{\ell}_{\mathfrak{F}}$ does not have an algebraic semantics based on $\Alg(\vdash^{\ell}_{\mathfrak{F}})$. As $\Alg(\vdash^{\ell}_{\mathfrak{F}})$ is the variety of modal algebras, we are done.
\end{proof}

\begin{problem}
Characterize the local modal logics $\vdash$ possessing a standard algebraic semantics. More in general, while this paper focuses mostly on nonstandard equational completeness theorems, it would be interesting to obtain similar descriptions of logics admitting standard equational completeness theorems.
\end{problem}

\section{Intermezzo: graph-based logics}

This section contains a useful characterization of graph-based logics with an algebraic semantics. The reader may safely skip to the next section, after having absorbed the statement of Theorem \ref{Thm:graph-based-solved}.

In the next definition, given a nonempty set $X \subseteq \omega$, we denote by $\textup{gcd}(X)$ the greatest common divisor of the elements of $X$ and assume that $\textup{gcd}(\emptyset) \coloneqq 0$.

\begin{law}\label{Def:rules}
Let $\LL$ be a graph-based language comprising a unary symbol $\Box$ and a constant $\textbf{c}_{i}$, and let $k$ be a positive integer. 
\benormal
\item Let $\mathcal{U}$ be the set of rules of the following form, where $n \in \omega$,
\[
x, y, \Box^{n}x \rhd \Box^{n}y.
\]
\item Let $\mathcal{S}_{k, i}$ be the set of rules of the following form, where $n \in \omega$,
\[
x, \Box^{k+ n}x \rhd \Box^{n} \textbf{c}_{i} \text{ and }x, \Box^{n} \textbf{c}_{i} \rhd \Box^{k+n} x.
\]
\item Let $\mathcal{R}_{k, i}$ be the set of rules of the form
\[
\{ \Box^{t}x \} \cup \{ \Box^{u_{j}}x \colon s > j \in \omega\} \cup \{ \Box^{v_{j}}x \colon s > j \in \omega \} \rhd \Box^{t+g}x,
\]
where $t, s, g \in \omega$ and $\{ u_{j} \colon s > j \in \omega \} \cup \{ v_{j} \colon s > j \in \omega \} \subseteq \omega$ are such that $u_{j} < v_{j}$ for all $s > j \in \omega$, and $g$ is a multiple  of $\textup{gcd}(\{ v_{j} - u_{j} \colon s > j \in \omega \})$.
\item Let $\mathcal{I}_{k, i}$ be the set of rules of the form
\[
\{ \Box^{u_{j}}x \colon s > j \in \omega \} \cup \{ \Box^{v_{j}}x \colon s > j \in \omega \} \rhd \Box^{g}\textbf{c}_{i},
\]
where $s, g \in \omega$ and $\{ u_{j} \colon s > j \in \omega \} \cup \{ v_{j} \colon s > j \in \omega \} \subseteq \omega$ are such that $u_{j} < v_{j}$ for all $s >  j \in \omega$, and $g + k$ is a multiple  of $\textup{gcd}(\{ v_{j} - u_{j} \colon s > j \in \omega \})$.
\enormal
\end{law}

The aim of this section is to establish the following:

\begin{Theorem}\label{Thm:graph-based-solved}
A graph-based logic $\vdash$  has an algebraic semantics if and only if one of the following conditions holds:
\benroman
\item\label{graph-thm-1} $\vdash$ is assertional;
\item\label{graph-thm-2} $\vdash$ has a unary connective $\Box$ and $x \vdash \Box x$;
\item\label{graph-thm-3} $\vdash$ has a unary connective $\Box$  and the rules in $\mathcal{U} \cup \mathcal{S}_{k, i} \cup \mathcal{R}_{k, i} \cup \mathcal{I}_{k, i}$ are valid in $\vdash$  for some positive integer $k$ and some ordinal $i < \alpha$;
\item\label{graph-thm-4} $\vdash$ is almost assertional and there is a set of equations $\btau(x)$ satisfying conditions (\ref{Eq:almost-assertional-1}, \ref{Eq:almost-assertional-2}) in Proposition \ref{Prop:almost-assertional}.
\eroman
Moreover, when condition (\ref{graph-thm-3}) holds, $k$ can be taken to be the least positive integer $n$ for which the rules in $\mathcal{S}_{n, i}$ are valid in $\vdash$.
\end{Theorem}

The remaining part of the section is devoted to proving Theorem \ref{Thm:graph-based-solved}. We begin with the following observation:

\begin{Proposition}\label{Prop:x-proves-box-x}
Let $\vdash$ be a graph-based logic with a unary connective. The following conditions are equivalent:
\benroman
\item\label{item:x-proves-box-x-1} $\vdash$ has a $\btau$-algebraic semantics for some $\btau$ such that
\[
\btau \subseteq \{ 	\Box^{n}x \thickapprox \Box^{m}x \colon n, m \in \omega \} \cup \{ \Box^{n}\textbf{c}_{i} \thickapprox \Box^{m}\textbf{c}_{j} \colon i, j < \alpha \text{ and }n, m \in \omega \};
\]
\item\label{item:x-proves-box-x-2} $\vdash$ has a $\btau$-algebraic semantics where $\btau \coloneqq \{ x \thickapprox \Box x \}$;
\item\label{item:x-proves-box-x-3} $x \vdash \Box x$.
\eroman
\end{Proposition}

\begin{proof}
The implication (\ref{item:x-proves-box-x-2})$\Rightarrow$(\ref{item:x-proves-box-x-1}) is obvious. To prove the implication (\ref{item:x-proves-box-x-1})$\Rightarrow$(\ref{item:x-proves-box-x-3}), let $\class{K}$ be the  algebraic semantics given by condition (\ref{item:x-proves-box-x-1}). To show that $x \vdash \Box x$, it suffices to prove that $\Box^{\A}a \in \btau(\A)$ for all $\A \in \class{K}$ and $a \in \btau(\A)$. To this end, consider $\A \in \class{K}$, $a \in \btau(\A)$, and an equation $\epsilon \thickapprox \delta$ in $\btau$. We need to establish
\begin{equation}\label{Eq_x-proves-box-x-part-1}
\epsilon^{\A}(\Box^{\A}a) = \delta^{\A}(\Box^{\A}a).
\end{equation}
Now, recall from the assumption that either $\epsilon = \Box^{n}\textbf{c}_{i}$ and $\delta = \Box^{m}\textbf{c}_{j}$ for some $i, j < \alpha$ and $n, m \in \omega$, or $\epsilon = \Box^{n}x$ and $\delta = \Box^{m}x$ for some $n, m \in \omega$. Observe that if $\epsilon = \Box^{n}\textbf{c}_{i}$ and $\delta = \Box^{m}\textbf{c}_{j}$, then $a \in \btau(\A)$ implies $\underbrace{\Box^{\A} \dots \Box^{\A}}_{n\text{-times}}\textbf{c}_{i}^{\A} = \underbrace{\Box^{\A} \dots \Box^{\A}}_{m\text{-times}}\textbf{c}_{j}^{\A}$ and, therefore,
\[
\epsilon^{\A}(\Box^{\A} a) = \underbrace{\Box^{\A} \dots \Box^{\A}}_{n\text{-times}}\textbf{c}_{i}^{\A} = \underbrace{\Box^{\A} \dots \Box^{\A}}_{m\text{-times}}\textbf{c}_{j}^{\A} = \delta^{\A}(\Box^{\A} a).
\]
Then we consider the case where $\epsilon = \Box^{n}x$ and $\delta = \Box^{m}x$. From $a \in \btau(\A)$ it follows
\[
\underbrace{\Box^{\A} \dots \Box^{\A}}_{n\text{-times}}a = \epsilon^{\A}(a) = \delta^{\A}(a) = \underbrace{\Box^{\A} \dots \Box^{\A}}_{m\text{-times}}a.
\]
Consequently,
\[
\epsilon^{\A}(\Box^{\A}a) = \Box^{\A}\underbrace{\Box^{\A} \dots \Box^{\A}}_{n\text{-times}}a = \Box^{\A}\underbrace{\Box^{\A} \dots \Box^{\A}}_{m\text{-times}}a = \delta^{\A}(\Box^{\A}a).
\]
This establishes (\ref{Eq_x-proves-box-x-part-1}) and, therefore, $x \vdash \Box x$.

It only remains to prove the implication (\ref{item:x-proves-box-x-3})$\Rightarrow$(\ref{item:x-proves-box-x-2}). To this end, set $\btau \coloneqq \{ x \thickapprox \Box x \}$. By Proposition \ref{Prop:alg-sem-meaning}, in order to show that $\vdash$ has a $\btau$-algebraic semantics, it suffices to check that for every $\Gamma \cup \{ \varphi \} \subseteq Fm$,
\begin{equation}\label{Eq_x-proves-box-x-part-2}
\text{if }\btau(\varphi) \subseteq \theta(\Gamma, \btau)\text{, then }\Gamma \vdash \varphi.
\end{equation}
To this end, consider $\Gamma \cup \{ \varphi \} \subseteq Fm$ such that $\btau(\varphi) \subseteq \theta(\Gamma, \btau)$. Then define
\[
\phi \coloneqq \{ \langle \epsilon, \delta\rangle \in Fm \times Fm : \text{either }\epsilon, \delta \in \textup{Cn}_{\vdash}(\Gamma) \text{ or }\epsilon = \delta \}.
\]
We shall see that $\phi$ is a congruence of $\Fm$. As $\phi$ is clearly an equivalence relation on $Fm$, it suffices to show that $\phi$ preserves the operation $\Box$. To this end, consider $\langle \epsilon, \delta \rangle \in \phi$. There are two cases: either $\epsilon = \delta$ or $\epsilon \ne \delta$. If $\epsilon = \delta$, then $\Box \epsilon = \Box \delta$ and, therefore, $\langle \Box \epsilon, \Box \delta \rangle \in \phi$. Then we consider the case where $\epsilon \ne \delta$. As $\langle \epsilon, \delta \rangle \in \phi$, we get $\epsilon, \delta \in \textup{Cn}_{\vdash}(\Gamma)$. Now, recall from the assumption that $x \vdash \Box x$ and, therefore, $\epsilon \vdash \Box \epsilon$ and $\delta \vdash \Box \delta$. Together with $\epsilon, \delta \in \textup{Cn}_{\vdash}(\Gamma)$, this implies $\Box \epsilon, \Box \delta \in  \textup{Cn}_{\vdash}(\Gamma)$, whence $\langle \Box \epsilon, \Box \delta \rangle \in \phi$. We conclude that $\phi$ is a congruence of $\Fm$.

We shall prove that $\phi$ extends $\theta(\Gamma, \btau)$. As $\phi$ is a congruence of $\Fm$, it will suffice to show that the generators of $\theta(\Gamma, \btau)$ belong to $\phi$, i.e., that $\langle \gamma, \Box \gamma \rangle \in \phi$ for every $\gamma \in \Gamma$. To this end, consider $\gamma \in \Gamma$. By $x \vdash \Box x$, we obtain $\Box \gamma \in \textup{Cn}_{\vdash}(\Gamma)$, whence $\langle \gamma, \Box \gamma \rangle \in \phi$. As a consequence, $\theta(\Gamma, \btau) \subseteq \phi$. In particular, this implies
\[
\{ \langle \varphi, \Box \varphi \rangle \} = \btau(\varphi) \subseteq \theta(\Gamma, \btau) \subseteq \phi.
\]
As $\varphi \ne \Box \varphi$, the fact that $\langle \varphi, \Box \varphi \rangle \in \phi$ implies $\varphi, \Box \varphi \in \textup{Cn}_{\vdash}(\Gamma)$. Consequently, $\Gamma \vdash \varphi$, establishing (\ref{Eq_x-proves-box-x-part-2}). Hence, we conclude that $\vdash$ has a $\btau$-algebraic semantics.
\end{proof}

A logic is said to be \textit{mono-unary} if its language consists of a unary connective only. As a consequence of Proposition \ref{Prop:x-proves-box-x} we obtain a new proof of the following:
\begin{Corollary}[\protect{\cite[Thm.\ 2.20]{BlRe03}}]
A mono-unary logic $\vdash$ has an algebraic semantics if and only if $x \vdash \Box x$.
\end{Corollary}
%

\begin{Proposition}\label{Prop:dicothomy}
Let $\vdash$ be a graph-based logic with an algebraic semantics. One of the following conditions holds:
\benroman
\item\label{dicothomy-0} $\vdash$ has a unary connective and $x \vdash \Box x$;
\item\label{dicothomy-1} $\vdash$ is assertional;
\item\label{dicothomy-2} $\vdash$ is almost assertional;
\item\label{dicothomy-3} $\vdash$ has a unary connective and a $\btau$-algebraic semantics where $\btau = \{ \Box^{k} x \thickapprox \Box^{n} \textbf{c}_{i} \}$, $0 \leq n < k \in \omega$, and $i < \alpha$.
\eroman
\end{Proposition}

\begin{proof}
By assumption, $\vdash$ has a $\brho$-algebraic semantics $\class{K}$. There are two cases: either $\vdash$ has a unary connective or it does not. First, suppose that $\vdash$ has no unary connective, i.e., its language comprises constant symbols only. If $\vdash$ is trivial, then it is either assertional or almost assertional and we are done. Then, we consider the case where $\vdash$ is nontrivial. By Lemma \ref{Lem:technical-x-in-tau},  there is a nontrivial equation $\epsilon \thickapprox \delta \in \brho$ such that $\Var(\epsilon) \cup \Var(\delta) = \{ x \}$. As $\LL_{\vdash}$ comprises constant symbols only, we can assume by symmetry that there is an $i < \alpha$ such that $\epsilon = x$ and $\delta = \textbf{c}_{i}$. The fact that $\vdash$ has a $\brho$-algebraic semantics where $x \thickapprox \textbf{c}_{i} \in \brho$ easily implies that $y \vdash \textbf{c}_{i}$ and $x, y, \varphi(x, \vec{z}) \vdash \varphi(y, \vec{z})$ for all $\varphi(v, \vec{z}) \in Fm$. By Proposition \ref{Prop:rules-almost-assertional}, we conclude that $\vdash$ is either assertional or almost assertional, as desired.

Therefore, it only remains to consider the case where $\vdash$ has a unary connective $\Box$. First, if condition (\ref{dicothomy-0}) holds, we are done. Then we consider the case where condition (\ref{dicothomy-0}) fails. By the implication (\ref{item:x-proves-box-x-3})$\Rightarrow$(\ref{item:x-proves-box-x-1}) in Proposition \ref{Prop:x-proves-box-x}, this yields that there are $s, t \in \omega$ and an ordinal $j < \alpha$ such that
\[
\brho \cap \{ \Box^{s} x \thickapprox \Box^{t} \textbf{c}_{j}, \Box^{t}  \textbf{c}_{j} \thickapprox \Box^{s} x \} \ne \emptyset.
\] 
Consequently, for every $\A \in \class{K}$ and $a \in \brho(\A)$,
\[
\underbrace{\Box^{\A} \dots \Box^{\A}}_{s\text{-times}}a = \underbrace{\Box^{\A} \dots \Box^{\A}}_{t\text{-times}} \textbf{c}_{j}^{\A}.
\]
Accordingly, the set
\begin{align*}
S \coloneqq \{ m \in \omega \colon&\text{there are }n \in \omega \text{ and an ordinal }i < \alpha \text{ such that for all }\A \in \class{K},\\
&\text{if }a \in \brho(\A),\text{ then } \underbrace{\Box^{\A} \dots \Box^{\A}}_{m\text{-times}}a = \underbrace{\Box^{\A} \dots \Box^{\A}}_{n\text{-times}} \textbf{c}_{i}^{\A}  \}
\end{align*}
is nonempty. Let then $k$ be the minimum of $S$. There are $n \in \omega$ and an ordinal $i < \alpha$ such that 
\begin{equation}\label{Eq:almost-assertional-deciding-equations-zzz}
\underbrace{\Box^{\A} \dots \Box^{\A}}_{k\text{-times}}a = \underbrace{\Box^{\A} \dots \Box^{\A}}_{n\text{-times}} \textbf{c}_{i}^{\A}
\end{equation}
for all $\A \in \class{K}$ and $a \in \brho(\A)$. Furthermore, define
\[
\btau \coloneqq \{ \Box^{k}x \thickapprox \Box^{n} \textbf{c}_{i} \}.
\]

\begin{Claim}\label{Claim:dicothomy}
If $\A \in \class{K}$ and $\brho(\A) \ne \emptyset$, then $\btau(\A) = \brho(\A)$.
\end{Claim}

\begin{proof}
Consider $\A \in \class{K}$ such that there is some $b \in \brho(\A)$. From the definition of $\btau$ it follows immediately that $\brho(\A) \subseteq \btau(\A)$. To prove the other inclusion, consider an element $a \in \btau(\A)$ and an equation $\epsilon \thickapprox \delta$ in $\brho$. We need to show $\epsilon^{\A}(a) = \delta^{\A}(a)$.
As $\vdash$ is graph-based, one of the following conditions holds:
\benroman
\item \label{item-equ-1-tau}$\epsilon = \Box^{s}x$ and $\delta = \Box^{t}x$ for some $s, t \in \omega$;
\item \label{item-equ-2-tau}$\epsilon = \Box^{s}\textbf{c}_{j}$ and $\delta = \Box^{t}x$ for some $s, t \in \omega$ and $j < \alpha$;
\item \label{item-equ-3-tau}$\epsilon = \Box^{s}x$ and $\delta = \Box^{t}\textbf{c}_{j}$ for some $s, t \in \omega$ and $j < \alpha$;
\item \label{item-equ-4-tau}$\epsilon = \Box^{s} \textbf{c}_{h}$ and $\delta = \Box^{t}\textbf{c}_{j}$ for some $s, t \in \omega$ and $h, j < \alpha$.
\eroman

(\ref{item-equ-1-tau}): If $s = t$, then clearly $\epsilon^{\A}(a) = \delta^{\A}(a)$. Then by symmetry we can assume that $s < t$. There are two cases: either $k \leq t$ or $t < k$. First suppose that $k \leq t$. Observe that for every $\B \in \class{K}$ and $d \in \brho(\B)$,
\[
\underbrace{\Box^{\B} \dots \Box^{\B}}_{s\text{-times}}d = \underbrace{\Box^{\B} \dots \Box^{\B}}_{t\text{-times}}d = \underbrace{\Box^{\B} \dots \Box^{\B}}_{(t-k)\text{-times}}\underbrace{\Box^{\B} \dots \Box^{\B}}_{k\text{-times}}d = \underbrace{\Box^{\B} \dots \Box^{\B}}_{(t-k)\text{-times}}\underbrace{\Box^{\B} \dots \Box^{\B}}_{n\text{-times}}\textbf{c}_{i}^{\B}.
\]
The first equality in the above display follows from $\epsilon \thickapprox \delta \in \brho$ and the third one from (\ref{Eq:almost-assertional-deciding-equations-zzz}). Furthermore, together with the minimality of $k$, the above display implies $k \leq s$. Thus, $k \leq s, t$. Bearing this in mind, observe that from $b \in \brho(\A)$ and (\ref{Eq:almost-assertional-deciding-equations-zzz}) it follows
\[
\underbrace{\Box^{\A} \dots \Box^{\A}}_{(s-k+n)\text{-times}}\textbf{c}_{i}^{\A} = \underbrace{\Box^{\A} \dots \Box^{\A}}_{(s-k)\text{-times}}\underbrace{\Box^{\A} \dots \Box^{\A}}_{k\text{-times}}b = \underbrace{\Box^{\A} \dots \Box^{\A}}_{(t-k)\text{-times}}\underbrace{\Box^{\A} \dots \Box^{\A}}_{k\text{-times}}b = \underbrace{\Box^{\A} \dots \Box^{\A}}_{(t-k+n)\text{-times}}\textbf{c}_{i}^{\A}.
\]
Together with $a \in \btau(\A)$, this implies
\[
\underbrace{\Box^{\A} \dots \Box^{\A}}_{(s-k)\text{-times}}\underbrace{\Box^{\A} \dots \Box^{\A}}_{k\text{-times}}a = \underbrace{\Box^{\A} \dots \Box^{\A}}_{(s-k+n)\text{-times}}\textbf{c}_{i}^{\A} = \underbrace{\Box^{\A} \dots \Box^{\A}}_{(t-k+n)\text{-times}}\textbf{c}_{i}^{\A} = \underbrace{\Box^{\A} \dots \Box^{\A}}_{(t-k)\text{-times}}\underbrace{\Box^{\A} \dots \Box^{\A}}_{k\text{-times}}a.
\]
Hence, we conclude $\epsilon^{\A}(a) = \delta^{\A}(a)$, as desired.

Then we consider the case where $t < k$. We shall see that this case never happens, i.e., it is contradictory. To prove this, notice that for every $\B \in \class{K}$ and $d \in \brho(\B)$,
\[
\underbrace{\Box^{\B} \dots \Box^{\B}}_{(k-t+s)\text{-times}}d = \underbrace{\Box^{\B} \dots \Box^{\B}}_{k\text{-times}} d = \underbrace{\Box^{\B} \dots \Box^{\B}}_{n\text{-times}} \textbf{c}_{i}^{\B}.
\]
The first equality above follows from $\epsilon \thickapprox \delta \in \brho$ and the second one from (\ref{Eq:almost-assertional-deciding-equations-zzz}). Now, together with the minimality of $k$, the above display implies $k - t+ s \geq k$. But this contradicts the fact that $s < t$. Hence, the case where $t \leq k$ never happens.

(\ref{item-equ-2-tau}): As $\Box^{s} \textbf{c}_{j} \thickapprox \Box^{t}x \in \brho$, by minimality of $k$ we get $k \leq t$. Bearing this in mind, 
\[
\underbrace{\Box^{\A} \dots \Box^{\A}}_{(t-k+n)\text{-times}}\textbf{c}_{i} = \underbrace{\Box^{\A} \dots \Box^{\A}}_{t\text{-times}}b = \underbrace{\Box^{\A} \dots \Box^{\A}}_{s\text{-times}} \textbf{c}_{j}^{\A}.
\]
The first equality above follows from (\ref{Eq:almost-assertional-deciding-equations-zzz}) and $b \in \brho(\A)$ and the second one from $\epsilon \thickapprox \delta \in \brho$ and $b \in \brho(\A)$. From the above display and $a \in \btau(\A)$ it follows
\[
\epsilon^{\A}(a) = \underbrace{\Box^{\A} \dots \Box^{\A}}_{s\text{-times}} \textbf{c}_{j}^{\A} = \underbrace{\Box^{\A} \dots \Box^{\A}}_{(t-k+n)\text{-times}}\textbf{c}_{i} = \underbrace{\Box^{\A} \dots \Box^{\A}}_{t\text{-times}}a = \delta^{\A}(a).
\]

Case (\ref{item-equ-3-tau}) is analogous to case  (\ref{item-equ-2-tau}). Finally, in case (\ref{item-equ-4-tau}) the equality $\epsilon^{\A}(a) = \delta^{\A}(a)$ holds vacuously, because $\brho(\A) \ne \emptyset$.
\end{proof}

Now, if $\vdash$ has theorems, then it is assertional by Corollary \ref{Cor:graph-based-with-theorems}. In this case, condition (\ref{dicothomy-1}) of Proposition \ref{Prop:dicothomy} holds and we are done. Then we consider the case where $\vdash$ lacks theorems. Since $\class{K}$ is a $\brho$-algebraic semantics for $\vdash$, we can apply Proposition \ref{Prop:matrix-sem-alg-sem} obtaining that $\{ \langle \A, \brho(\A)\rangle \colon \A \in \class{K}\}$ is a matrix semantics for $\vdash$. As $\vdash$ lacks theorems, by Lemma \ref{Lem:no-thms-matrix} the following is also a matrix semantics for $\vdash$:
\[
\class{N} \coloneqq \{ \langle \A, \brho(\A)\rangle \colon \A \in \class{K} \text{ such that } \brho(\A) \ne \emptyset \} \cup \{ \langle \Fm, \emptyset \rangle \}.
\]
Lastly, with an application of Claim \ref{Claim:dicothomy}, 
\begin{equation}\label{Eq:last-corrections-to-the-paper-546}
\class{N} = \{ \langle \A, \btau(\A)\rangle \colon \A \in \class{K} \text{ such that }\rho(\A) \ne \emptyset \} \cup \{ \langle \Fm, \emptyset \rangle \}.
\end{equation}

There are two cases: either $n < k$ or $k \leq n$. If $n < k$, there is no formula $\varphi$ such that $\Box^{k}\varphi = \Box^{n}\textbf{c}_{i}$. Consequently, $\btau(\Fm) = \emptyset$. Together with (\ref{Eq:last-corrections-to-the-paper-546}), this yields
\[
\class{N} = \{ \langle \A, \btau(\A)\rangle \colon \A \in \class{K} \text{ such that }\rho(\A) \ne \emptyset \} \cup \{ \langle \Fm, \btau(\Fm) \rangle \}.
\]
Since $\class{N}$ is a matrix semantics for $\vdash$, this implies that the class of algebraic reducts of $\class{N}$ is a $\btau$-algebraic semantics for $\vdash$. Thus, condition (\ref{dicothomy-3}) of Proposition \ref{Prop:dicothomy} holds and we are done.

It only remains to consider the case where $k \leq n$. In this case, from the fact that $\class{N}$ is a matrix semantics for $\vdash$ and (\ref{Eq:last-corrections-to-the-paper-546}) it follows $x \vdash \Box^{n-k}\textbf{c}_{i}$. Recall from Proposition \ref{Prop:unital} that $\vdash$ has a unital matrix semantics $\class{M}$. Together with $x \vdash \Box^{n-k}\textbf{c}_{i}$, this implies
\[
F = \{ \underbrace{\Box^{\A} \dots \Box^{\A}}_{(n-k)\text{-times}} \textbf{c}_{i}^{\A} \}
\]
for all $\langle \A, F \rangle \in \class{M}$ such that $F \ne \emptyset$. As $\vdash$ lacks theorems, we conclude that $\vdash$ is almost assertional and, therefore, that condition (\ref{dicothomy-2}) of Proposition \ref{Prop:dicothomy} holds.
\end{proof}

The proof of the next result relies on two technical lemmas which, for the sake of readability, are established in the Appendix.

\begin{Proposition}\label{Prop:n-k-trick}
The following conditions are equivalent for a graph-based logic $\vdash$ with a unary connective:
\benroman
\item\label{the-nk-case-1} $\vdash$ has a $\btau$-algebraic semantics where $\btau = \{ \Box^{k}x \thickapprox \Box^{n} \textbf{c}_{i} \}$ for some nonnegative integers $n < k$ and ordinal $i < \alpha$;
\item\label{the-nk-case-2} $\vdash$ has a $\btau$-algebraic semantics where $\btau = \{ \Box^{k}x \thickapprox \textbf{c}_{i} \}$ for some positive integer $k$ and ordinal $i < \alpha$;
\item\label{the-nk-case-3} The rules in $\mathcal{U} \cup \mathcal{S}_{k, i} \cup \mathcal{R}_{k, i} \cup \mathcal{I}_{k, i}$ are valid in $\vdash$ for some positive integer $k$ and ordinal $i < \alpha$.
\eroman
Moreover, when condition (\ref{the-nk-case-3}) holds, $k$ can be taken to be the least positive integer $n$ for which the rules in $\mathcal{S}_{n, i}$ are valid in $\vdash$.
\end{Proposition}

\begin{proof}
The implication (\ref{the-nk-case-2})$\Rightarrow$(\ref{the-nk-case-1}) is obvious (just take $n \coloneqq 0$).

(\ref{the-nk-case-1})$\Rightarrow$(\ref{the-nk-case-3}): Assume that $\vdash$ has a $\btau$-algebraic semantics where $\btau = \{ \Box^{k}x \thickapprox \Box^{n} \textbf{c}_{i} \}$ for some nonnegative integers $n < k$ and ordinal $i < \alpha$. We begin with the following observation:

\begin{Claim}\label{Claim:the-n-k-trick-1}
The rules in $\mathcal{S}_{k-n, i}$ are valid in $\vdash$.
\end{Claim}

\begin{proof}
We need to prove that for every $m \in \omega$,
\[
x, \Box^{k-n+ m}x \sineq \Box^{m}\textbf{c}_{i}, x.
\]
To this end, consider $m \in \omega$. As $\vdash$ has a $\btau$-algebraic semantics, by Proposition \ref{Prop:alg-sem-meaning} it suffices to establish
\begin{align*}
\langle \Box^{k+k - n+m}x, \Box^{n}\textbf{c}_{i}\rangle &\in \theta(\{ x, \Box^{m}\textbf{c}_{i}\}, \btau)\\
\langle \Box^{k+m}\textbf{c}_{i}, \Box^{n}\textbf{c}_{i}\rangle &\in \theta(\{ x, \Box^{k-n+m}x\}, \btau).
\end{align*}
The above conditions follow, respectively, from Lemmas \ref{Lem:graph-based-combinatorics1} and \ref{Lem:graph-based-combinatorics2}, which are proved in the Appendix.
\end{proof}

Observe that $k -n$ is a positive integer, since by assumption $k > n$. Thus by Claim \ref{Claim:the-n-k-trick-1} there exists the least positive integer $m$ such that the rules $\mathcal{S}_{m, i}$ are valid in $\vdash$. We shall prove that the rules in $\mathcal{U} \cup \mathcal{R}_{m, i} \cup \mathcal{I}_{m, i}$ are valid in $\vdash$. This will establish both condition (\ref{the-nk-case-3}) and the remark at the end of the statement of Proposition \ref{Prop:n-k-trick}. 

First notice that, as $\vdash$ has an algebraic semantics, by Propositions \ref{Prop:unital} and \ref{Prop:unital-rules} the rules in $\mathcal{U}$ are valid in $\vdash$. Then, consider a typical rule in $\mathcal{R}_{m, i}$:
\begin{equation}\label{Eq:valid-rule-R}
\{ \Box^{t}x \} \cup \{ \Box^{u_{j}}x \colon s > j \in \omega\} \cup \{ \Box^{v_{j}}x \colon s > j \in \omega \} \rhd \Box^{t+g}x.
\end{equation}
Define
\[
\Gamma \coloneqq \{ \Box^{t}x \} \cup \{ \Box^{u_{j}} x \colon s > j \in \omega \} \cup \{ \Box^{v_{j}}x \colon s > j \in \omega \}.
\]
By Proposition \ref{Prop:alg-sem-meaning}, as $\vdash$ has a $\btau$-algebraic semantics for $\vdash$, in order to prove that the rule (\ref{Eq:valid-rule-R}) is valid in $\vdash$,  it suffices to show that
\begin{equation}\label{Eq:valid-rule-R-sem}
\langle \Box^{k}\Box^{t+g}x, \Box^{n}\textbf{c}_{i}\rangle \in \theta(\Gamma, \btau).
\end{equation}
But this is a consequence of the ``if'' part of Lemma \ref{Lem:graph-based-combinatorics1}. Hence we conclude that the rules in $\mathcal{R}_{m, i}$ are valid in $\vdash$. The fact that the rules in $\mathcal{I}_{m, i}$ are valid in $\vdash$ is proved similarly, by replacing the only application of Lemma \ref{Lem:graph-based-combinatorics1} with one of Lemma \ref{Lem:graph-based-combinatorics2}.

(\ref{the-nk-case-3})$\Rightarrow$(\ref{the-nk-case-2}): Suppose that the rules in $\mathcal{U} \cup \mathcal{S}_{k, i} \cup \mathcal{R}_{k, i} \cup \mathcal{I}_{k, i}$ are valid in $\vdash$ for some positive integer $k$ and some ordinal $i < \alpha$. Then let $\mathcal{R}_{k, i}^{+}$ be the set of rules of the form
\[
\{ \Box^{t}x \} \cup \{ \Box^{u_{j}}y_{j} \colon s > j \in \omega \} \cup \{ \Box^{u_{j}}y_{j} \colon s > j \in \omega \} \rhd \Box^{t+g}x,
\]
where $t, s, g \in \omega$ and $\{ u_{j} \colon s > j \in \omega \} \cup \{ v_{j} \colon s > j \in \omega \} \subseteq \omega$ are such that $u_{j} < v_{j}$ for all $s > j \in \omega$, and $g$ is a multiple  of $\textup{gcd}(\{ v_{j} - u_{j} \colon s > j \in \omega \})$. 

Let also $\mathcal{I}_{k, i}^{+}$ be the set of rules of the form
\[
\{ \Box^{u_{j}}y_{j}  \colon s > j \in \omega \} \cup \{ \Box^{v_{j}}y_{j}  \colon s > j \in \omega \} \rhd \Box^{g}\textbf{c}_{i},
\]
where $s, g \in \omega$ and $\{ u_{j} \colon s > j \in \omega \} \cup \{ v_{j} \colon s > j \in \omega \} \subseteq \omega$ are such that $u_{j} < v_{j}$ for all $s > j \in \omega$, and $g + k$ is a multiple  of $\textup{gcd}(\{ v_{j} - u_{j} \colon s > j \in \omega \})$. 

Notice that, as the rules in $\mathcal{U} \cup \mathcal{R}_{k, i} \cup \mathcal{I}_{k, i}$ are valid in $\vdash$, it is not hard to see that so are those in $\mathcal{R}_{k, i}^{+} \cup \mathcal{I}_{k, i}^{+}$. Furthermore, notice that
\begin{equation}\label{Eq:the-n-k-trick-2}
\Box^{w}\textbf{c}_{i} \vdash \Box^{2w+k}\textbf{c}_{i}, \text{ for all }w \in \omega.
\end{equation}
To justify the above display, consider $w \in \omega$. As the rules in $\mathcal{S}_{k, i}$ are valid in $\vdash$, we have $x, \Box^{w}\textbf{c}_{i} \vdash \Box^{w+k}x$. Hence, by substitution invariance, we obtain $\Box^{w}\textbf{c}_{i}, \Box^{w}\textbf{c}_{i} \vdash \Box^{2w+k}\textbf{c}_{i}$ and, therefore, $\Box^{w}\textbf{c}_{i} \vdash \Box^{2w+k}\textbf{c}_{i}$.

Now, we shall prove that $\vdash$ has a $\btau$-algebraic semantics for
\[
\btau \coloneqq \{ \Box^{k}x \thickapprox \textbf{c}_{i} \}.
\]
By Proposition \ref{Prop:alg-sem-meaning} it suffices to check that 
\begin{equation}\label{Eq:zdfgrwhxlqcq}
\text{if }\btau(\varphi) \subseteq \theta(\Gamma, \btau)\text{, then }\Gamma \vdash \varphi
\end{equation}
for all $\Gamma \cup \{ \varphi \} \subseteq Fm$. 

To this end, consider $\Gamma \cup \{ \varphi \} \subseteq Fm$ such that $\btau(\varphi) \subseteq \theta(\Gamma, \btau)$, i.e., such that
\begin{equation}\label{Eq:generation-congruence-n-k-i-rule1}
\langle \Box^{k}\varphi, \textbf{c}_{i}\rangle \in \theta(\Gamma, \btau).
\end{equation}
We denote by $At(\LL_{\vdash})$ the set of atomic formulas of $\vdash$, i.e., the union of $\Var$ and the set of constants of $\LL_{\vdash}$. As $\vdash$ is graph-based, there are $p \in At(\LL_{\vdash})$ and $h \in \omega$ such that $\varphi = \Box^{h}p$. There are two cases: either $p \ne \textbf{c}_{i}$ or $p = \textbf{c}_{i}$.

First we consider the case where $p \ne \textbf{c}_{i}$. By applying the ``only if'' part of Lemma \ref{Lem:graph-based-combinatorics1} to (\ref{Eq:generation-congruence-n-k-i-rule1}), we obtain
\begin{align*}
s, s^{\ast} \in \omega \text{ and }\{ q_{j} \colon &s > j \in \omega \} \subseteq At(\LL_{\vdash}) \text{ and}\\
\{ u_{j} \colon s > j \in \omega \} \cup \{ v_{j} \colon & s >  j \in \omega \} \cup \{ w_{j^{\ast}} \colon s^{\ast} > i \in \omega \} \subseteq \omega
\end{align*}
such that
\benormal
\item\label{Eq:cong-gen-last-trick-1} $u_{j} < v_{j}$, for all $s > j \in \omega$;
\item\label{Eq:cong-gen-last-trick-2} $\Box^{u_{j}}q_{j}, \Box^{v_{j}}q_{j}, \Box^{w_{j^{\ast}}}\textbf{c}_{i} \in \Gamma$, for all $s > j \in \omega$ and $s^{\ast} > j^{\ast} \in \omega$;
\item\label{Eq:cong-gen-last-trick-3} $h = t+g$ for some $t, g \in \omega$ such that $\Box^{t}p \in \Gamma$ and $g$ is a multiple of
\[
d \coloneqq \textup{gcd}(\{ v_{j} - u_{j} \colon s > j \in \omega \} \cup \{ k + w_{j^{\ast}} \colon s^{\ast} > j^{\ast} \in \omega \}).
\]
\enormal

Define
\begin{align*}
\Delta \coloneqq &\{ \Box^{t}p \} \cup \{ \Box^{u_{j}}q_{j} \colon s > j \in \omega \} \cup \{ \Box^{v_{j}}q_{j} \colon s > j \in \omega \} \cup \\
& \{ \Box^{w_{j^{\ast}}}\textbf{c}_{i} \colon s^{\ast} > j^{\ast} \in \omega \} \cup \{ \Box^{2w_{j^{\ast}}+k}\textbf{c}_{i} \colon s^{\ast} > j^{\ast} \in \omega \}.
\end{align*}
From conditions (\ref{Eq:cong-gen-last-trick-2}) and (\ref{Eq:cong-gen-last-trick-3}) above and (\ref{Eq:the-n-k-trick-2}) it follows that $\Gamma \vdash \Delta$. Hence, to prove $\Gamma \vdash \varphi$, it suffices to show $\Delta \vdash \varphi$. From the fact that $g$ is a multiple of $d$, it follows that it is also a multiple of
\[
\textup{gcd}(\{ r - r' \colon \Box^{r}q, \Box^{r'}q \in \Delta \text{ for some }q \in At(\LL_{\vdash}) \text{ such that } r > r'\}).
\]
Together with the fact that $2w_{j^{\ast}}+ k > w_{j^{\ast}}$ for all $s^{\ast} > j^{\ast} \in \omega$ (as $k$ is positive) and that the rules in $\mathcal{R}_{k, i}^{+}$ are valid in $\vdash$, this implies $\Delta \vdash \Box^{t+ g}p$. As, by condition (\ref{Eq:cong-gen-last-trick-3}), $h = t+g$, this amounts to $\Delta \vdash \Box^{h}p$. Since $\varphi = \Box^{h}p$, we conclude $\Delta \vdash \varphi$ and, therefore, $\Gamma \vdash \varphi$.

Then we consider the case where $p = \textbf{c}_{i}$. By applying the ``only if'' part of Lemma \ref{Lem:graph-based-combinatorics2} to (\ref{Eq:generation-congruence-n-k-i-rule1}), we obtain 
\begin{align*}
s, s^{\ast} \in \omega \text{ and }\{ q_{j} \colon &s > j \in \omega \} \subseteq At(\LL_{\vdash}) \text{ and}\\
\{ u_{j} \colon s > j \in \omega \} \cup \{ v_{j} \colon & s >  j \in \omega \} \cup \{ w_{j^{\ast}} \colon s^{\ast} > i \in \omega \} \subseteq \omega
\end{align*}
validating conditions (\ref{Eq:cong-gen-last-trick-1}) and (\ref{Eq:cong-gen-last-trick-2}) above and such that
\benormal
\item[4.]\label{Eq:cong-gen-last-trick-3-c} $h+k$ is a multiple of 
\[
d \coloneqq \textup{gcd}(\{ v_{j} - u_{j} \colon s > j \in \omega \} \cup \{ k + w_{j^{\ast}} \colon s^{\ast} > j^{\ast} \in \omega \}).
\]
\enormal
Defining $\Delta$ as in the previous case (but omitting $\Box^{t}p$), we get $\Gamma \vdash \Delta$. Thus, to prove $\Gamma \vdash \varphi$, it suffices to show $\Delta \vdash \varphi$. Now, as $h+k$ is a multiple of $d$, it is also a multiple of
\[
\textup{gcd}(\{ r - r' \colon \Box^{r}q, \Box^{r'}q \in \Delta \text{ for some }q \in At(\LL_{\vdash}) \text{ such that } r > r'\}).
\]
Together with the fact that $2w_{j^{\ast}}+ k > w_{j^{\ast}}$ for all $s^{\ast} > j^{\ast} \in \omega$ (as $k$ is positive) and that the rules in $\mathcal{I}_{k, i}^{+}$ are valid in $\vdash$, this implies $\Delta \vdash \Box^{h}p$. Hence, we conclude that $\Gamma \vdash \varphi$, establishing (\ref{Eq:zdfgrwhxlqcq}).
\end{proof}

We are now ready to prove the main result of this section:

\begin{proof}[Proof of Theorem \ref{Thm:graph-based-solved}.]
To prove the ``only if'' part, suppose that $\vdash$ has an algebraic semantics and that conditions (\ref{graph-thm-1}, \ref{graph-thm-2}, \ref{graph-thm-4}) fail. As $\vdash$ has an algebraic semantics, Proposition \ref{Prop:almost-assertional} and the assumption that condition (\ref{graph-thm-4}) fails imply that $\vdash$ is not almost assertional. Therefore, from Proposition \ref{Prop:dicothomy} it follows that $\vdash$ has a $\btau$-algebraic semantics where $\btau = \{ \Box^{k} x \thickapprox \Box^{n} \textbf{c}_{i} \}$ for some nonnegative integers $n < k$ and ordinal $i < \alpha$. By Proposition \ref{Prop:n-k-trick} this implies that the rules in $\mathcal{U} \cup \mathcal{S}_{m, i} \cup \mathcal{R}_{m, i} \cup \mathcal{I}_{m, i}$ are valid in $\vdash$ for some positive integer $m$, i.e., that condition (\ref{graph-thm-3}) holds, as desired.

Then we turn to prove the ``if'' part. The fact that each condition in the list (\ref{graph-thm-1}, \ref{graph-thm-2}, \ref{graph-thm-3}, \ref{graph-thm-4}) implies that $\vdash$ has an algebraic semantics follows respectively from Propositions \ref{Prop:assertional}, \ref{Prop:x-proves-box-x}, \ref{Prop:n-k-trick}, and \ref{Prop:almost-assertional}.

Finally, the last sentence in the statement of Theorem \ref{Thm:graph-based-solved} follows from the last part of Proposition \ref{Prop:n-k-trick}.
\end{proof}

\section{Locally tabular logics}\label{Sec:strongly-finite}

\begin{law}
A logic $\vdash$ is said to be \textit{locally tabular} if for every integer $n \geq 1$ there are only finitely many formulas in variables $x_1, \dots, x_n$ up to logical equivalence in $\vdash$.\footnote{The expression \textit{locally tabular} originates in the field of modal and intermediate logics.}
\end{law}

An algebra is said to be \textit{locally finite} when its finitely generated subalgebras are finite. Similarly, a class of algebras $\class{K}$ is \textit{locally finite} when its members are.


\begin{Proposition}\label{Prop:locally-finite}
The following conditions are equivalent for a logic $\vdash$:
\benroman
\item\label{item:loc-fin:new-def-1} $\vdash$ is locally tabular;
\item\label{item:loc-fin:new-def-2} $\Alg(\vdash)$ is locally finite;
\item\label{item:loc-fin:new-def-3} $\vdash$ has a matrix semantics whose algebraic reducts generate a locally finite variety.
\eroman
\end{Proposition}

\begin{proof}
The implications (\ref{item:loc-fin:new-def-1})$\Rightarrow$(\ref{item:loc-fin:new-def-2}) and (\ref{item:loc-fin:new-def-3})$\Rightarrow$(\ref{item:loc-fin:new-def-1}) follow immediately from Lemma \ref{Lem:matrix-to-rules-1}.

(\ref{item:loc-fin:new-def-2})$\Rightarrow$(\ref{item:loc-fin:new-def-3}): Suppose that $\Alg(\vdash)$ is locally finite.\ Since $\Alg(\vdash)$ contains a free algebra with a denumerable set of free generators by Lemma \ref{Lem:free-algebra}, the fact that $\Alg(\vdash)$ is locally finite implies that also variety generated by $\Alg(\vdash)$ is locally finite. Since $\vdash$ has a matrix semantics whose algebraic reducts belong to $\Alg(\vdash)$ by Corollary \ref{Cor:the-place-of-reductions}, we are done.
\end{proof}

In this section we present a characterization of locally tabular logics with an algebraic semantics. For the sake of readability, this result is split in the next three propositions.

\begin{Proposition}\label{Prop:SF-no-graph-based}
Every locally tabular logic that is not graph-based has an algebraic semantics.
\end{Proposition}

\begin{proof}
Let $\vdash$ be a locally tabular logic that is not graph-based. The latter assumption guarantees that there is a complex formula $\varphi$ such that $\Var(\varphi) = \{ x \}$. As $\vdash$ is locally tabular, there are two nonnegative integers $m \leq n$ such that
\[
\varphi^{m}(x) \equiv_{\vdash} \varphi^{n+1}(x).
\]
Since $\varphi$ is a complex formulas, the two formulas in the above display are different. Moreover, from $\Var(\varphi) = \{ x \}$ it follows
\[
\Var(\varphi^{m}(x)) \cup \Var(\varphi^{n+1}(x)) = \{ x \}.
\]
As $\vdash$ is not graph-based, we can apply Theorem \ref{Thm:sufficient} obtaining that $\vdash$ has an algebraic semantics.
\end{proof}

\begin{Corollary}
A locally tabular logic with theorems has an algebraic semantics if and only if either it is graph-based and assertional or it is not graph-based.
\end{Corollary}

\begin{proof}
The ``only if'' part follows from Corollary \ref{Cor:graph-based-with-theorems}, while the ``if' one is a consequence of Propositions \ref{Prop:assertional} and \ref{Prop:SF-no-graph-based}.
\end{proof}

In view of Proposition \ref{Prop:SF-no-graph-based}, it only remains to describe graph-based locally tabular logics with an algebraic semantics. To some extent, this has already been done in Theorem \ref{Thm:graph-based-solved}. However, in this section we aim for a more computational characterization, which will be instrumental to derive some decidability results (Theorem \ref{Thm:comput-thm}).

To this end, notice that for every locally tabular graph-based logic $\vdash$ with a unary connective there are two nonnegative integers $m \leq n$ such that $\Box^{m}x \equiv_{\vdash} \Box^{n+1}x$. The next result explains when such a logic has an algebraic semantics.

\begin{Proposition}\label{Prop:complicated-case-comput}
Let $\vdash$ be a graph-based logic with a unary connective for which there are two nonnegative integers $m \leq n$ such that $\Box^{m}x \equiv_{\vdash} \Box^{n+1}x$. Then $\vdash$ has an algebraic semantics if and only if $x, y, \Box^{t}x \vdash \Box^{t}y$ for all nonnegative integers $t \leq n$ and one of the following conditions holds:
\benroman
\item\label{comput-g-1} $x \vdash \Box x$;
\item\label{comput-g-2} $y \vdash \Box^{t}p$ for some $p \in \{ x \} \cup \{ \textbf{c}_{i} \colon i \text{ is an ordinal}< \alpha \}$ and nonnegative integer $t \leq n$;
\item\label{comput-g-3} The following conditions hold for some ordinal $i < \alpha$ and positive integer $k \leq n$:
\benormal
\item[\emph{(a)}]\label{comput-g-a} $x, \Box^{t+k}x \sineq \Box^{t}\textbf{c}_{i}, x$ for all nonnegative integers $t \leq n$;
\item[\emph{(b)}]\label{comput-g-b} For all $s, g, h, t \in \omega$ such that $s \leq (2n-m+1)^{2}$ and $g, h, t \leq 2n-m+1$, and for all $\{ u_{j} \colon s > j \in \omega \} \cup \{ v_{j} \colon s > j \in \omega \} \subseteq \omega$,
\begin{align*}
\{ \Box^{t}x \} \cup \{ \Box^{u_{j}}x \colon s > j \in \omega\} \cup \{ \Box^{v_{j}}x \colon s > j \in \omega \} &\vdash \Box^{t+g}x\\
\{ \Box^{u_{j}}x \colon s > j \in \omega \} \cup \{ \Box^{v_{j}}x \colon s > j \in \omega \} &\vdash \Box^{h}\textbf{c}_{i},
\end{align*}
provided that $u_{j} < v_{j} \leq n + n-m+1$ for all $s > j \in \omega$, and that $g$ and $h+k$ are a multiples  of $\textup{gcd}(\{ v_{j} - u_{j} \colon s > j \in \omega \})$.
\enormal
\eroman
\end{Proposition}

\begin{proof}
First observe that $\Box^{m}x \equiv_{\vdash} \Box^{n+1}x$ implies $\Box^{m}x \sineq \Box^{n+1}x$. Consequently, for every formula $\varphi$ there are some $p \in \Var \cup \{ \textbf{c}_{i} \colon  i \text{ is an ordinal}< \alpha \}$ and some nonnegative integer $t \leq n$ such that $\varphi \sineq \Box^{t}p$. In what follows, we will make repeated use of this observation.
 
To prove the ``only if'' part, suppose that $\vdash$ has an algebraic semantics. From Propositions \ref{Prop:unital-rules} and \ref{Prop:unital} it follows $x, y, \Box^{t}x \vdash \Box^{t}y$ for all nonnegative integer $t \leq n$. Furthermore, as $\vdash$ is a graph-based logic with an algebraic semantics, one of conditions (\ref{graph-thm-1})--(\ref{graph-thm-4}) in Theorem \ref{Thm:graph-based-solved} holds. 

If condition (\ref{graph-thm-1}) of Theorem \ref{Thm:graph-based-solved} holds, then $x \vdash \Box x$ and we are done. 

Then suppose that condition (\ref{graph-thm-2}) or (\ref{graph-thm-4}) in Theorem \ref{Thm:graph-based-solved} holds. In this case, $\vdash$ is either assertional or almost assertional. Thus, by Proposition \ref{Prop:rules-almost-assertional}, there is a formula $\psi(x)$ such that $y \vdash \psi(x)$. As $\psi$ can be chosen of the form $\Box^{t}p$ for some $p \in \{ x \} \cup \{ \textbf{c}_{i} \colon  i \text{ is an ordinal}< \alpha \}$ and some nonnegative integer $t \leq n$, condition (\ref{comput-g-2}) in the statement holds. 

It only remains to consider the case where condition (\ref{graph-thm-3}) in Theorem \ref{Thm:graph-based-solved} holds. In this case, the rules in $\mathcal{S}_{k, i} \cup \mathcal{R}_{k, i} \cup \mathcal{I}_{k, i}$ are valid in $\vdash$  for some positive integer $k$ and some ordinal $i < \alpha$. Accordingly, conditions (a) and (b) in the statement hold. Therefore, to establish condition (\ref{comput-g-3}) in the statement, it only remains to prove that $k \leq n$. To this end, in view of Theorem \ref{Thm:graph-based-solved}, we can safely assume that $k$ is the least positive integer $t$ for which the rules in $\mathcal{S}_{t, i}$ are valid in $\vdash$. Furthermore, let $n \geq s \in \omega$ be such that $\Box^{s}x \sineq \Box^{k}x$. By substitution invariance, for every $t \in \omega$,
\[
\Box^{s+t}x \sineq \Box^{k+t}x.
\]
Then consider $t \in \omega$. As the rules in $\mathcal{S}_{k, i}$ are valid in $\vdash$, we obtain $x, \Box^{k+t}x \sineq \Box^{t} \textbf{c}_{i}, x$. Together with the above display, this yields $x, \Box^{s+t}x \sineq \Box^{t} \textbf{c}_{i}, x$. Consequently, the rules in $\mathcal{S}_{s, i}$ are also valid in $\vdash$. By the minimality of $k$, we conclude that $k \leq s \leq n$. 

Then we turn to prove the ``if'' part. First, suppose that $x, y, \Box^{t}x \vdash \Box^{t}y$ for all nonnegative integers $t \leq n$. Since $m \leq n$ and $\Box^{m}x \sineq \Box^{n+1}x$, this implies that the rules in $\mathcal{U}$ are valid in $\vdash$. To conclude the proof, it suffices to show that each condition in the list (\ref{comput-g-1}, \ref{comput-g-2}, \ref{comput-g-3}) forces $\vdash$ to possess an algebraic semantics. 

By Theorem \ref{Thm:graph-based-solved} this is true for (\ref{comput-g-1}). Then suppose that condition (\ref{comput-g-2}) holds. As the rules in $\mathcal{U}$ are valid in $\vdash$, (\ref{comput-g-2}) and Proposition \ref{Prop:rules-almost-assertional} imply that $\vdash$ is either assertional or almost assertional. Recall that assertional logics have an algebraic semantics by Proposition \ref{Prop:assertional}. Therefore we only need to consider the case where $\vdash$ is almost assertional. Since $\Box^{m}x \equiv_{\vdash} \Box^{n+1}x$, for every $\varphi(v, \vec{z}) \in Fm$,
\[
x, \varphi(\Box^{m}x, \vec{z}) \sineq \varphi(\Box^{n+1}, \vec{z}), x.
\]
Moreover, as $m \leq n$, there is no formula $\varphi$ such that $\Box^{m}\varphi = \Box^{n+1}\varphi$. Thus the set $\btau \coloneqq \{ \Box^{m}x \thickapprox \Box^{n+1}x \}$ satisfies conditions (\ref{Eq:almost-assertional-1}, \ref{Eq:almost-assertional-2}) in Proposition \ref{Prop:almost-assertional}. We conclude that condition (\ref{graph-thm-4}) in Theorem \ref{Thm:graph-based-solved} holds, whence $\vdash$ has an algebraic semantics.

It only remains to consider the case where condition (\ref{comput-g-3}) in the statement holds. By Theorem \ref{Thm:graph-based-solved}, in order to prove that $\vdash$ has an algebraic semantics, it suffices to show that the rules in $\mathcal{U} \cup \mathcal{S}_{k, i} \cup \mathcal{R}_{k, i} \cup \mathcal{I}_{k, i}$ are valid in $\vdash$. We know that this is the case for the rules in $\mathcal{U}$. Moreover, from (a) and $\Box^{m}x \sineq \Box^{n+1}x$ it follows that also the rules in $\mathcal{S}_{k, i}$ are valid in $\vdash$. Therefore, it only remains to prove that the rules in $\mathcal{R}_{k, i} \cup \mathcal{I}_{k, i}$ are valid in $\vdash$. 

To this end, consider a generic rule in $\mathcal{R}_{k, i}$:
\[
\{ \Box^{t}x \} \cup \{ \Box^{u_{j}}x \colon s > j \in \omega\} \cup \{ \Box^{v_{j}}x \colon s > j \in \omega \} \rhd \Box^{t+g}x.
\]
From the definition $\mathcal{R}_{k, i}$ it follows that $g$ is a multiple of 
\[
d \coloneqq \textup{gcd}(\{ v_{j} - u_{j} \colon s > j \in \omega \}).
\]
Define
\begin{align*}
\Gamma &\coloneqq \{ \Box^{t}x \} \cup \{ \Box^{u_{j}}x \colon s > j \in \omega\} \cup \{ \Box^{v_{j}}x \colon s > j \in \omega \}\\
\Delta & \coloneqq \{ \Box^{j}x \colon n + n-m+1 \geq j \in \omega \text{ and }\Box^{j}x \sineq \gamma \text{ for some }\gamma \in \Gamma \}.
\end{align*}
Since $\Gamma \vdash \Delta$, to conclude the proof, it suffices to prove $\Delta \vdash \Box^{t+g}x$. To this end, define
\[
d^{\ast} \coloneqq \textup{gcd}(\{ v - u \colon \Box^{u}x, \Box^{v}x \in \Delta \text{ and }u < v\}).
\]

\begin{Claim}\label{Claim:numerical}
$d$ is a multiple of $d^{\ast}$.
\end{Claim}

\begin{proof}
By the definition of $d$ and $d^{\ast}$, it suffices to prove that $v_{j} - u_{j}$ is a multiple of $d^{\ast}$, for every $s > j \in \omega$. To this end, consider $s > j \in \omega$. If $v_{j} \leq n + n-m+1$, then $u_{j}, v_{j} \leq n + n-m+1$, since $u_{j} < v_{j}$. Consequently, $\Box^{u_{j}}x, \Box^{v_{j}} x \in \Delta$.  Thus, $v_{j} - u_{j}$ is a multiple of $d^{\ast}$ by definition of $d^{\ast}$.

Then suppose that $v_{j} > n + n - m + 1$. In this case, there are $1 \leq p_{v} \in \omega$ and $n - m \geq q_{v} \in \omega$ such that
\[
v_{j} = n + q_{v} + 1 + p_{v}( n- m+ 1 ).
\] 
Notice that $\Box^{m} x \sineq \Box^{n+1}x$ yields $\Box^{v_{j}}x \sineq \Box^{m+q_{v}}x \sineq \Box^{n+q_{v}+1}x$. Together with $\Box^{v_{j}}x \in \Gamma$ and $m+q_{v}, n+q_{v}+1 \leq n + n - m + 1$, this implies
\[
\Box^{m + q_{v}}x, \Box^{n+q_{v}+1}x \in \Delta.
\]
Notice that $m + q_{v} < n+q_{v}+1$, since $m \leq n$. Consequently,
\[
n-m+1 = (m+q_{v}+1)-(m+q_{v})
\]
is a multiple of $d^{\ast}$, by definition of $d^{\ast}$. Thus, there is an $r_{1} \in \omega$ such that
\begin{equation}\label{Eq:multiple-last-day-of-writing}
d^{\ast} r_{1} = n-m+1.
\end{equation}

Now, if $u_{j} \leq n$, then $\Box^{u_{j}}x \in \Delta$. Together with $\Box^{n+q_{v}+1}x \in \Delta$ and $u_{j} \leq n < n + q_{v} +1$, this implies that $d^{\ast}$ is a divisor of $n+ q_{v} + 1 - u_{j}$. Thus, there is an $r_{2} \in \omega$ such that $d^{\ast}r_{2} = n+ q_{v} + 1 - u_{j}$. Consequently,
\begin{align*}
v_{j} - u_{j} &= n + q_{v} + 1 + p_{v}( n- m+ 1 ) - u_{j} \\
&= (n + q_{v} + 1 - u_{j}) + p_{v}( n- m+ 1 )  \\
&= d^{\ast} (r_{2} + r_{1}p_{v}).
\end{align*}
Hence, we conclude that $v_{j} - u_{j}$ is a multiple of $d^{\ast}$, as desired.

Then we consider the case where $u_{j} > n$.  There are $1 \leq p_{u} \in \omega$ and $n - m \geq q_{u} \in \omega$ such that
\[
u_{j} = n + q_{u} + 1 + p_{u}( n- m+ 1 ).
\] 
Furthermore, from $\Box^{u_{j}} \in \Gamma$ it follows
\[
\Box^{m + q_{u}}x, \Box^{n+q_{u}+1}x \in \Delta.
\]
Since $u_{j} < v_{j}$, clearly $p_{u} \leq p_{v}$. There are two cases: either $q_{u} \leq q_{v}$ or $q_{v} < q_{u}$. 

If $q_{u} \leq q_{v}$, then 
\[
v_{j} - u_{j} = (p_{v} - p_{u}) (n - m + 1) + q_{v} - q_{u}.
\]
We shall see that $q_{v} - q_{u}$ is a multiple of $d^{\ast}$. If $q_{v} = q_{u}$, this is obvious. Then suppose that $q_{v} > q_{u}$. In this case, $n+q_{u}+1 < n+q_{v}+1$. Together with $\Box^{n+q_{v}+1}x, \Box^{n+q_{u}+1}x \in \Delta$ and $n+q_{v}+1 \leq n + n - m +1$, this implies that $q_{v} - q_{u}$ is a multiple of $d^{\ast}$ by definition of $d^{\ast}$. Together with (\ref{Eq:multiple-last-day-of-writing}) and the above display, this guarantees that $v_{j} - u_{j}$ is a multiple of $d^{\ast}$.

It only remains to consider the case where $q_{v} < q_{u}$. As $u_{j} < v_{j}$, in this case $p_{u} < p_{v}$. Consequently,
\begin{align*}
v_{j} - u_{j} = (p_{v}- p_{u}- 1)(n - m+1) + (n - m + 1 + q_{v} - q_{u}),
\end{align*}
where $(n - m + 1 + q_{v} - q_{u}) \geq 0$, since $q_{u} \leq n-m$. In view of (\ref{Eq:multiple-last-day-of-writing}) and the above display, to conclude the proof it suffices to show that $n - m + 1 + q_{v} - q_{u}$ is a multiple of $d^{\ast}$. To this end, recall that $\Box^{n+q_{v}+ 1}x, \Box^{m + q_{u}}x \in \Delta$. Together with $q_{u} \leq n-m \leq n < n +1$, this implies that $n - m + 1 + q_{v} - q_{u} = (n+q_{v}+ 1)- (m + q_{u})$ is a multiple of $d^{\ast}$.
\end{proof}

Now, let $n \geq t^{\ast} \in \omega$ be such that $\Box^{t} x\sineq \Box^{t^{\ast}}x$. As $\Box^{t}x \in \Gamma$, we get $\Box^{t^{\ast}}x \in \Delta$. Furthermore, from the definition of $\Delta$ it follows that the set
\[
Z \coloneqq \{ v - u \colon \Box^{u}x, \Box^{v}x \in \Delta \text{ and }u < v\}
\] 
has cardinality $\leq (2n-m+2)^{2}$. Lastly, $d^{\ast} \leq 2n - m +1$. To prove this, observe that if $Z = \emptyset$, then $d^{\ast} = 0$ and we are done. Otherwise, $d^{\ast}$ is a divisor of $v - u$ for some $u, v \leq 2n -m +1$. Consequently, $d^{\ast} \leq 2n - m +1$, as desired.

Then we can apply assumption (b), obtaining $\Delta \vdash \Box^{t^{\ast}+ d^{\ast}}x$. Recall that $g$ is a multiple of $d$. Thus, by Claim \ref{Claim:numerical} there is an $r \in \omega$ such that $rd^{\ast} = g$. Since the rules in $\mathcal{U}$ are valid in $\vdash$, we have $\Box^{t^{\ast}}x, \Box^{t^{\ast}+ d^{\ast}}x \vdash \Box^{t^{\ast}+rd^{\ast}}x$ and, therefore, $\Delta \vdash \Box^{t^{\ast}+g}x$. Together with $\Box^{t} x\sineq \Box^{t^{\ast}}x$, this implies $\Delta \vdash \Box^{t+g}x$. We conclude that the rules in $\mathcal{R}_{k, i}$ are valid in $\vdash$.

It only remains to prove that the same holds for the rules in $\mathcal{I}_{k, i}$. To this end, consider a generic rule in $\mathcal{I}_{k, i}$:
\[
\{ \Box^{u_{j}}x \colon s > j \in \omega\} \cup \{ \Box^{v_{j}}x \colon s > j \in \omega \} \rhd \Box^{h}\textbf{c}_{i}.
\]
From the definition $\mathcal{I}_{k, i}$ it follows that $h+k$ is a multiple of 
\[
d \coloneqq \textup{gcd}(\{ v_{j} - u_{j} \colon s > j \in \omega \}).
\]
As $k$ is positive, this implies that $d \geq 1$.

Define
\begin{align*}
\Gamma &\coloneqq \{ \Box^{u_{j}}x \colon s > j \in \omega\} \cup \{ \Box^{v_{j}}x \colon s > j \in \omega \}\\
\nabla & \coloneqq \{ \Box^{j}x \colon n + n-m+1 \geq j \in \omega \text{ and }\Box^{j}x \sineq \gamma \text{ for some }\gamma \in \Gamma \}\\
Z &\coloneqq \{ v - u \colon \Box^{u}x, \Box^{v}x \in \nabla \text{ and }u < v\}\\
d^{\ast} &\coloneqq \textup{gcd}(Z).
\end{align*}
Notice that the proof of Claim \ref{Claim:numerical} can be repeated with these definitions, thus $d$ is a multiple of $d^{\ast}$. Consequently, as $d$ is positive, $d^{\ast}$ is as well. In particular, this implies that $Z \ne \emptyset$ and, therefore, $\nabla \ne \emptyset$.

We shall define a set of formulas $\Delta$ as follows. If there is a $\Box^{u}x \in \nabla$ such that $u \geq m$, then just set $\Delta \coloneqq \nabla$. Otherwise, 
\[
\nabla \subseteq \{ \Box^{0}x, \Box^{1} x, \dots, \Box^{m-1}x \}. 
\]
Since $\nabla \ne \emptyset$, we can choose a formula $\Box^{t}x \in \nabla$. Because of $Z \ne \emptyset$ and of the above display,
\begin{equation}\label{last-display-1}
1 \leq d^{\ast} \leq m-1 \leq n \text{ and }t \leq m-1.
\end{equation}
As the rules in $\mathcal{R}_{k, i}$ are valid in $\vdash$,  we get $\nabla \vdash \Box^{d^{\ast}+ t}x$. Notice that $\Box^{t}x, \Box^{d^{\ast}+t}x \vdash \Box^{rd^{\ast}+t}x$ for all $r \in \omega$, since the rules in $\mathcal{U}$ are valid in $\vdash$. Consequently, 
\begin{equation}\label{last-display-2}
\nabla \vdash \Box^{rd^{\ast} +t}x\text{, for all }r \in \omega. 
\end{equation}
By (\ref{last-display-1}) and (\ref{last-display-2}) there is a $u \in \omega$ such that $m \leq u \leq n+ n - m +1$ and $\nabla \vdash \Box^{u}x$. Then set
\[
\Delta \coloneqq \nabla \cup \{ \Box^{r}x \colon n + n - m + 1 \geq r \in \omega \text{ and }\Box^{u}x \vdash \Box^{r}x \}.
\]
This completes the definition of $\Delta$.

Observe that  in both cases $\Gamma \vdash \nabla \vdash \Delta$. Thus, to conclude the proof, it suffices to show that $\Delta \vdash \Box^{h}\textbf{c}_{i}$. To this end, define
\[
d^{\ast\ast} \coloneqq \textup{gcd}(\{ v - u \colon \Box^{u}x, \Box^{v}x \in \Delta \text{ and }u < v\}).
\]
Recall that $d$ is a multiple of $d^{\ast}$. Furthermore, as $\nabla \subseteq \Delta$, $d^{\ast}$ is clearly a multiple of $d^{\ast\ast}$. Thus, $d$ is also a multiple of $d^{\ast\ast}$. Furthermore, $d^{\ast\ast}$ is positive, as $d$ is.

Since $h+k$ is a multiple of $d$ and $d^{\ast\ast}$ is a divisor of $d$, there is an $r_{1} \in \omega$ such that
\begin{equation}\label{last-display-3}
h + k = r_{1}d^{\ast\ast}.
\end{equation}
Now, if $h \leq n + n - m + 1$, we can apply condition (b) to $\Delta$, obtaining that $\Delta \vdash \Box^{h}\textbf{c}_{i}$, as desired.  The suppose that $n + n - m + 1 < h$. In this case, there are $1 \leq p \in \omega$ and $n- m\geq q \in \omega$ such that
\begin{equation}\label{last-display-4}
h = p(n-m+1) + n + q.
\end{equation}
By definition of $\Delta$ there is a $m \leq u \in \omega$ such that $\Box^{u}x \in \Delta$. Since $\Box^{m}x \sineq \Box^{n+1}x$ and $m \leq u$, there is an $n \geq  r \in \omega$ such that
\[
\Box^{u}x \sineq \Box^{r}x \sineq \Box^{r + n - m + 1}x. 
\]
By the definition of $\Delta$, we obtain $\Box^{r}x, \Box^{r + n - m + 1}x \in \Delta$. Consequently, $n - m +1$ is a multiple of $d^{\ast\ast}$, i.e., there is an $r_{2} \in \omega$ such that
\[
r_{2}d^{\ast\ast} = n- m +1.
\]
Together with (\ref{last-display-3}) and (\ref{last-display-4}) this implies
\[
r_{1}d^{\ast\ast} = pr_{2}d^{\ast\ast} + n + q + k.
\]
Consequently, $n+q+k$ is a multiple of $d^{\ast\ast}$. Since $q \leq n-m$, we can apply condition (b) to $\Delta$, obtaining that $\Delta \vdash \Box^{n+q}\textbf{c}_{i}$. Lastly, from $\Box^{m}x \sineq \Box^{n+1}x$ and (\ref{last-display-4}) it follows $\Box^{n+q}\textbf{c}_{i} \sineq \Box^{h}\textbf{c}_{i}$, whence $\Delta \vdash \Box^{h}\textbf{c}_{i}$, as desired.
\end{proof}

We conclude with a characterization of (locally tabular) logics, formulated in languages comprising constant symbols only, with an algebraic semantics.

\begin{Proposition}\label{Prop:no-language}
Let $\vdash$ be a logic such that $\LL_{\vdash}$ comprises constant symbols only. Then $\vdash$ has an algebraic semantics if and only if one of the following conditions holds:
\benroman
\item\label{item:no-language-1} Either $\emptyset \vdash x$ or $\emptyset \vdash \textbf{c}_{i}$ for some ordinal $i < \alpha$;
\item\label{item:no-language-2} $x \vdash \textbf{c}_{i}$ and $x \vdash \textbf{c}_{j}$ for some ordinals $i < j < \alpha$;
\item\label{item:no-language-3} $x \vdash \textbf{c}_{k}$ and $\textbf{c}_{i} \sineq \textbf{c}_{j}$ for some ordinals $i, j, k < \alpha$ such that $i \ne j$.
\eroman
\end{Proposition}

\begin{proof}
To prove the ``only if'' part, suppose that $\vdash$ has an algebraic semantics. Then $\vdash$ is either assertional or almost assertional by Theorem \ref{Thm:graph-based-solved}. If $\vdash$ is assertional, then it has theorems, whence condition (\ref{item:no-language-1}) holds. Then suppose that $\vdash$ is almost assertional. By Proposition \ref{Prop:almost-assertional} there is a set of equations $\btau(x)$ such that
\begin{equation}\label{Eq:constant-kvqlknlknvlkvnqlk}
x, \epsilon \sineq \delta, x\text{, for all }\epsilon \thickapprox \delta \in \btau, 
\end{equation}
and for which there is no substitution $\sigma$ such that $\sigma(\epsilon) = \sigma(\delta)$ for all $\epsilon \thickapprox \delta \in \btau$. Because $\LL_{\vdash}$ comprises constant symbols only, the latter condition implies that there are two ordinals $i < j < \alpha$ such that either $\textbf{c}_{i} \thickapprox \textbf{c}_{j} \in \btau$ or, by symmetry, $\{ x \thickapprox \textbf{c}_{i}, x \thickapprox \textbf{c}_{j} \} \subseteq \btau$. 

If $\{ x \thickapprox \textbf{c}_{i}, x \thickapprox \textbf{c}_{j} \} \subseteq \btau$, then (\ref{Eq:constant-kvqlknlknvlkvnqlk}) yields $x \vdash \textbf{c}_{i}$ and $x \vdash \textbf{c}_{j}$. Therefore, condition (\ref{item:no-language-2}) holds and we are done. Then we consider the case where $\textbf{c}_{i} \thickapprox \textbf{c}_{j} \in \btau$. By (\ref{Eq:constant-kvqlknlknvlkvnqlk}) we get $x, \textbf{c}_{i} \sineq \textbf{c}_{j}, x$. Notice that, by substitution invariance, this yields $\textbf{c}_{i} \sineq \textbf{c}_{j}$. Now, if $\vdash$ is trivial, then $x \vdash y$ and, therefore, $x \vdash \textbf{c}_{i}$. Accordingly, taking $k \coloneqq i$, condition (\ref{item:no-language-3}) holds and we are done. Then suppose that $\vdash$ is nontrivial, i.e., that $x \nvdash y$. Since $\vdash$ is almost assertional, by Proposition \ref{Prop:rules-almost-assertional} there is a formula $\psi(y)$ such that $x \vdash \psi(y)$. Since $x \nvdash y$, the formula $\psi(y)$ is not a variable. As $\LL_{\vdash}$ comprises constant symbols only, we conclude that $y = \textbf{c}_{k}$ for some ordinal $k < \alpha$. Therefore, $x \vdash \textbf{c}_{k}$, whence condition (\ref{item:no-language-3}) holds, as desired.

Then we turn to prove the ``if'' part. Suppose that a condition in (\ref{item:no-language-1}, \ref{item:no-language-2}, \ref{item:no-language-3}) holds. Because of the poorness of the language in which $\vdash$ if formulated, it is easy to see that $x, y, \varphi(x, \vec{z}) \vdash \varphi(y, \vec{z})$ for every formula $\varphi(v, \vec{z})$. Moreover, by assumption there is a formula $\psi(y)$ such that $x \vdash \psi(y)$. Consequently, we can apply Proposition \ref{Prop:rules-almost-assertional} obtaining that $\vdash$ is either assertional or almost assertional. If $\vdash$ is assertional, then it has an algebraic semantics by Proposition \ref{Prop:assertional}. Then suppose that $\vdash$ is almost assertional. As $\vdash$ lacks theorems, condition (\ref{item:no-language-1}) fails. Then either (\ref{item:no-language-2}) or (\ref{item:no-language-3}) holds. In both cases, define $\btau \coloneqq \{  \textbf{c}_{i} \thickapprox \textbf{c}_{j} \}$. It is easy to see that conditions (\ref{Eq:almost-assertional-1}) and (\ref{Eq:almost-assertional-2}) in Proposition \ref{Prop:almost-assertional} hold, whence $\vdash$ has an algebraic semantics.
\end{proof}

\section{Computational aspects}

Two of the most common ways to present a logic are by exhibiting either a class of matrices that induces the logic or a set of rules that axiomatizes it. In the latter case, the set of rules is sometimes called an \textit{Hilbert-style axiomatization} of the logic. We close our journey by studying the decidability of the problem of determining whether a logic presented by any of these two methods has an algebraic semantics, cf.\ \cite{Mor16a,Mor16compl}. To this end, a rule $\Gamma \rhd \varphi$ is said to be \textit{finite} when $\Gamma$ is finite. The aim of this section is to establish the following:

\begin{Theorem}\label{Thm:comput-thm}
\
\benroman
\item\label{item:comput-thm-1} The problem of determining whether a logic presented by a finite set of finite matrices in a finite language has an algebraic semantics is decidable.
\item\label{item:comput-thm-2} The problem of determining whether a locally tabular logic presented by a finite set of finite rules in a finite language has an algebraic semantics is decidable.
\item\label{item:comput-thm-3} The problem of determining whether a logic presented by a finite set of finite rules in a finite language has an algebraic semantics is undecidable.
\eroman
\end{Theorem}

\begin{proof}[Proof of Theorem \ref{Thm:comput-thm}(\ref{item:comput-thm-1})]
Let $\class{M}$ be a finite set of finite matrices in a finite language. As the variety generated by the algebraic reducts of $\class{M}$ is finitely generated, it is also locally finite \cite[Thm.\ II.10.16]{BuSa81}. Thus $\vdash_{\class{M}}$ is locally tabular by Proposition \ref{Prop:locally-finite}. 

In order to determine whether $\vdash_{\class{M}}$ has an algebraic semantics, our algorithm checks first if $\vdash_{\class{M}}$ is graph-based or not. If $\vdash_{\class{M}}$ is not graph-based, then it has an algebraic semantics in view of Proposition \ref{Prop:SF-no-graph-based} and we are done. Otherwise, $\vdash_{\class{M}}$ is graph-based and our algorithm checks whether its language comprises a unary symbol or not. If it does not, then $\vdash_{\class{M}}$ has an algebraic semantics if and only if one of conditions (\ref{item:no-language-1}, \ref{item:no-language-2}, \ref{item:no-language-3}) in Proposition \ref{Prop:no-language} holds. Notice that this can be checked mechanically, as $\class{M}$ is a finite set of finite matrices in a finite language. Therefore, it only remains to consider the case where the language of $\vdash_{\class{M}}$ comprises a unary symbol. As the $\class{M}$ is a finite set of finite matrices, our algorithm can find two nonnegative integers $m \leq n$ such that $\class{M} \vDash \Box^{m}x \thickapprox \Box^{n+1}x$. Consequently, $\Box^{m}x \equiv_{\vdash} \Box^{n+1}x$ and $\vdash_{\class{M}}$ has an algebraic semantics when one of conditions (\ref{comput-g-1}, \ref{comput-g-2}, \ref{comput-g-3}) in Proposition \ref{Prop:complicated-case-comput} holds. But, again, this can be checked mechanically, as $\class{M}$ is a finite set of finite matrices in a finite language.
\end{proof}

\begin{proof}[Proof of Theorem \ref{Thm:comput-thm}(\ref{item:comput-thm-2})]
In order to determine whether a locally tabular logic $\vdash$, presented by a finite set of finite rules in a finite language, has an algebraic semantics, our algorithm checks first whether $\vdash$ is graph-based. If $\vdash$ is not graph-based, then it has an algebraic semantics by Proposition \ref{Prop:SF-no-graph-based}. 

Suppose then that $\vdash$ is graph-based. In view of Theorem \ref{Thm:comput-thm}(\ref{item:comput-thm-1}), to determine whether $\vdash$ has an algebraic semantics, our algorithm needs only to produce a finite set of finite matrices that induces $\vdash$. To this end, we shall first detail the case where $\vdash$ has a unary connective $\Box$. In this case, the algorithm finds two nonnegative integers $m \leq n$ such that $\Box^{m} x \sineq \Box^{n+1}x$. Notice that $n$ and $m$ exist because $\vdash$ is locally tabular, and that we can find them mechanically by enumerating all proofs obtained from the Hilbert-style calculus axiomatizing $\vdash$.

Now, observe that, up to isomorphism, there are only finitely many $(2^{n+1}+1)$-generated $\LL_{\vdash}$-algebras $\A$ such that $\A \vDash \Box^{m}x \thickapprox \Box^{n+1}x$. Furthermore, these algebras $\A_{1}, \dots, \A_{m}$ can be produced mechanically. As $\vdash$ is presented by a finite set of finite rules, our algorithm can construct the following finite set of finite matrices:
\[
\class{M} \coloneqq \{ \langle \A_{i}, F \rangle \colon m \geq i \in \omega \text{ and the logic induced by }\langle \A_{i}, F \rangle \text{ extends }\vdash \}.
\]
Therefore, to conclude the proof, it only remains to show that $\vdash$ is induced by $\class{M}$. 

That $\vdash_{\class{M}}$ is an extension of $\vdash$ follows immediately from the definition of $\class{M}$. To prove the converse, consider $\Gamma \cup \{ \varphi \} \subseteq Fm$ such that $\Gamma \nvdash \varphi$. Observe that every formula $\delta \in Fm$ has the form $\Box^{t}p$ for some $t \in \omega$ and some $p$ that is either a constant or a variable. As $\Box^{m} x \sineq \Box^{n+1}x$, there is a $n \geq k \in \omega$ such that $\delta \sineq \Box^{k}p$. Accordingly, set $\delta^{\ast} \coloneqq \Box^{k}p$ where $k$ is the least natural number such that $\delta \sineq \Box^{k}p$. Clearly, $\delta \sineq \delta^{\ast}$. Similarly, given a set of formulas $\Gamma$, define
\[
\Gamma^{\ast} \coloneqq \{ \gamma^{\ast} \colon \gamma \in \Gamma \}.
\]
From $\Gamma \nvdash \varphi$ it follows $\Gamma^{\ast} \nvdash \varphi^{\ast}$. 

Now, we define an equivalence relation $R$ on the set
\[
\Var(\Gamma^{\ast}) \coloneqq \bigcup_{\gamma \in \Gamma^{\ast}} \Var(\gamma),
\]
by the rule
\[
\langle x, y \rangle \in R \Longleftrightarrow \{ k \in \omega \colon \Box^{k}x \in \Gamma^{\ast} \} = \{ k \in \omega \colon \Box^{k}y \in \Gamma^{\ast} \}.
\]
Observe that $\Var(\Gamma^{\ast}) / R$ has at most $2^{n+1}$ equivalence classes, because $\Gamma^{\ast}$ contains only formulas of the form $\Box^{t}p$ where $n \geq t \in \omega$ and $p$ is either a constant or a variable. Accordingly, we can choose representatives $v_{1}, \dots, v_{k}$ with $k \leq 2^{n+1}$ of the equivalence classes in $\Var(\Gamma^{\ast}) / R$ in such a way that if $\Var(\varphi) = \{ x \} \subseteq \Var(\Gamma^{\ast})$ for some variable $x$, then the representative of $x / R$ is $x$. Furthermore, we consider the unique substitution $\sigma$ defined for every $y \in \Var$ as
\[
\sigma(y) \coloneqq \begin{cases}
y & \text{if $y \notin \Var(\Gamma^{\ast})$}\\
v_{i} & \text{if $y \in \Var(\Gamma^{\ast})$ and $\langle y, v_{i}\rangle \in R$.}\\
\end{cases}
\]
Notice that $\sigma(\varphi^{\ast}) = \varphi^{\ast}$.

We shall prove that
\begin{equation}\label{Eq:sigma-Gamma-trick}
\sigma[\Gamma^{\ast}] \subseteq \Gamma^{\ast}.
\end{equation}
To this end, consider $\gamma \in \Gamma^{\ast}$. As $\vdash$ is graph-based, there are two cases: either $\Var(\gamma) = \emptyset$ or $\Var(\gamma)$ is a singleton. If $\Var(\gamma) = \emptyset$, then $\gamma$ is a closed formula, whence $\sigma(\gamma) = \gamma \in \Gamma^{\ast}$, as desired. Then we consider the case where $\Var(\gamma) = \{ x \}$ for some $x \in \Var$. By definition of $\Gamma^{\ast}$, there is a $n \geq t \in \omega$ such that $\gamma = \Box^{t}x$. Considering the unique $k \geq i \in \omega$ such that $\langle x, v_{i}\rangle \in R$, we obtain $\Box^{t}v_{i} \in \Gamma^{\ast}$. Consequently, $\sigma(\gamma) = \Box^{t}\sigma(x) = \Box^{t}v_{i} \in \Gamma^{\ast}$, establishing (\ref{Eq:sigma-Gamma-trick}).


Furthermore, as $\vdash$ is graph-based and substitution invariant, from $\Box^{m} x \sineq \Box^{n+1}x$ we obtain $\varphi(\Box^{m}x, \vec{z}) \sineq \varphi(\Box^{n+1}x, \vec{z})$ for all $\varphi(v, \vec{z}) \in Fm$ and, therefore, 
\[
\Box^{m}x \equiv_{\vdash}\Box^{n+1}x.
\]
Consequently, $\vdash$ is induced by a class of matrices $\class{N}$ validating the equation $\Box^{m}x \thickapprox \Box^{n+1}x$ by Lemma \ref{Lem:matrix-to-rules-1}. 

From $\Gamma^{\ast} \nvdash \varphi^{\ast}$ and (\ref{Eq:sigma-Gamma-trick}) it follows $\sigma[\Gamma^{\ast}] \nvdash \varphi^{\ast}$. Therefore, since $\class{N}$ is a matrix semantics for $\vdash$, there are a matrix $\langle \B, G \rangle \in \class{N}$ and a homomorphism $h \colon \Fm \to \B$ such that $h[\sigma[\Gamma^{\ast}]] \subseteq G$ and $h(\varphi^{\ast}) \notin G$. Let $\A$ be the subalgebra of $\B$ generated by $h(v_{1}), \dots, h(v_{k})$ if $\Var(\varphi^{\ast}) = \emptyset$ and by $h(v_{1}), \dots, h(v_{k}), h(x)$ if $\Var(\varphi^{\ast}) = \{ x \}$ for some variable $x$. As $k \leq 2^{n+1}$, in both cases $\A$ is a $(2^{n+1}+1)$-generated algebra. Moreover, as $\langle \B, G \rangle \in \class{N}$ and $\vdash$ is the logic induced by $\class{N}$, the logic induced by $\langle \A, G \cap A \rangle$ is an extension of $\vdash$. Lastly, $\A$ validates $\Box^{m}x \thickapprox \Box^{n+1}x$, as $\B$ does. Consequently, $\langle \A, G \cap A \rangle \cong \langle \A_{i}, F_{i}\rangle$ for some there is an $m \geq i \in \omega$. For the sake of simplicity, we shall assume that $\langle \A, G \cap A \rangle = \langle \A_{i}, F_{i}\rangle$ and, therefore, $\langle \A, G \cap A \rangle \in \class{M}$. By considering any homomorphism $g \colon \Fm \to \A$ that coincides with $h$ on the variables occurring in $\sigma[\Gamma^{\ast}]$ and $\varphi^{\ast}$, we obtain $g[\sigma[\Gamma^{\ast}]] \subseteq G \cap A$ and $g(\varphi^{\ast}) \notin G \cap A$. Since $\langle\A, G \cap A \rangle \in \class{M}$, this yields
\[
\sigma[\Gamma^{\ast}] \nvdash_{\class{M}} \varphi^{\ast}.
\]

Now, recall that $\sigma(\varphi^{\ast}) = \varphi^{\ast}$. Thus, from the above display it follows $\sigma[\Gamma^{\ast}] \nvdash_{\class{M}} \sigma(\varphi^{\ast})$. As $\vdash_{\class{M}}$ is substitution invariant, this implies $\Gamma^{\ast} \nvdash_{\class{M}} \varphi^{\ast}$. Finally, by $\class{M} \vDash \Box^{m}x \thickapprox \Box^{n+1}x$ we obtain $\Box^{m}x \sineq_{\class{M}} \Box^{n+1}x$, whence $\Gamma \sineq_{\class{M}} \Gamma^{\ast}$ and $\varphi \sineq_{\class{M}} \varphi^{\ast}$. Together with $\Gamma^{\ast} \nvdash_{\class{M}} \varphi^{\ast}$, this implies $\Gamma \nvdash_{\class{M}}\varphi$. Hence, $\vdash$ is an extension of $\vdash_{\class{M}}$, as desired. This concludes the proof that $\vdash$ is the logic induced by $\class{M}$.

It only remains to consider the case where $\vdash$ is a graph-based logic whose language comprises only constant symbols. But this case is handled similarly to the previous one (hint: repeat the proof by taking $n = m = 0$).
\end{proof}

\begin{problem}
Investigate the complexity of the decision problems mentioned in conditions (\ref{item:comput-thm-1}) and (\ref{item:comput-thm-2}) of Theorem \ref{Thm:comput-thm}.\footnote{The naive procedure for problem (\ref{item:comput-thm-1}) described above runs in exponential time in the length of the input.}
\end{problem}

The remaining part of this section is devoted to the proof of condition (\ref{item:comput-thm-3}) in Theorem \ref{Thm:comput-thm}. To this end, we assume the reader is familiar with computability theory, and sketch the basic definitions only to fix some terminology and conventions. 

By a \textit{Turing machine} ${\bf M}$ we understand a tuple $\langle P, Q, q_{0}, \delta \rangle$ where $P$ and $Q$ are sets of ``states'', $q_{0} \in Q$ is the \textit{initial state}, $Q$ the set of \textit{nonfinal states}, $P$ the set of \textit{final states}, and
\[
\delta \colon Q \times \{ 0, 1, \emptyset \} \to (Q \cup P) \times \{ 0, 1 \} \times \{ L, R \}.
\]
Instruction of the form $\delta(q, a) = \langle q', b, L\rangle$ (resp.\ $\delta(q, a)=\langle q', b, R\rangle$) should be understood as follows: if the machine ${\bf M}$ reads $a$ at state $q$, then it replaces $a$ with $b$, moves left (resp.\ right), and switches to state $q'$. 

Our Turing machines work on tapes that are infinite both to the left and to the right, can write only zeros and ones, but can read $0$, $1$, and the empty symbol $\emptyset$. At the beginning of the computation, the empty symbol $\emptyset$ occupies all cells in the tape.

A \textit{configuration} for a Turing machine ${\bf M}$ is a tuple $\langle q, \vec{w}, v, \vec{u}\rangle$ where $q \in Q \cup P$, $\vec{w}$ and $\vec{u}$ are either finite nonempty sequences of zeros and ones or the one-element sequence $\langle \emptyset \rangle$, and $v \in \{ \langle 0 \rangle, \langle 1 \rangle, \langle \emptyset \rangle \}$ which, in addition, satisfies the following requirement: if $\vec{w}$ and $\vec{u}$ are different from $\langle \emptyset \rangle$, then $v$ is also. Intuitively, the configuration $\langle q, \vec{w}, v, \vec{u}\rangle$ represents the instant in which ${\bf M}$ is in state $q$, it is reading the unique symbol in $v$, and the tape contains exactly the concatenation
\[
\dots \emptyset, \emptyset, \emptyset \rangle^{\frown}  \vec{w}^{\frown}v^{\frown}\vec{u}^{\frown} \langle \emptyset, \emptyset, \emptyset \dots
\]
Given two configurations $\texttt{c}$ and $\texttt{d}$ for ${\bf M}$, we say that $\texttt{c}$ \textit{yields} $\texttt{d}$ if ${\bf M}$ allows to move from $\texttt{c}$ to $\texttt{d}$ in a single step of computation. 

An \textit{input} is a finite sequence $\vec{t} = \langle t_{1}, \dots, t_{m} \rangle$ of zeros and ones such that $m \geq 2$.\footnote{The assumption that $m \geq 2$ is inessential and is intended to simplify the proof of Theorem \ref{Thm:comput-thm}(\ref{item:comput-thm-3}) only.} The \textit{initial configuration} for ${\bf M}$ under $\vec{t}$ is the tuple
\[
\textup{In}({\bf M}, \vec{t} \?\?) \coloneqq \langle q_{0}, \langle \emptyset \rangle, \langle t_{1}\rangle, \langle t_{2}, \dots, t_{m}\rangle \rangle.
\]
Then ${\bf M}$ is said to \textit{halt} on $\vec{t}$ if there is a finite sequence $\texttt{c}_{1}, \dots, \texttt{c}_{n}$ of configurations for ${\bf M}$ such that $\texttt{c}_{1} = \textup{In}({\bf M}, \vec{t} \?\?)$, the sate in $\texttt{c}_{n}$ belongs to the set $P$ of final states, and $\texttt{c}_{m}$ yields $\texttt{c}_{m+1}$ for every positive integer $m < n$. 

The \textit{halting problem} asks to determine, given a Turing machine ${\bf M}$ and an input $\vec{t}$, whether ${\bf M}$ halts on $\vec{t}$. It is well known that this problem is undecidable, as shown by Turing in \cite{Tu36}. In the remaining part of the section we shall see that the halting problem reduces to that of determining whether a logic presented by a finite set of finite rules in a finite language has an algebraic semantics. It follows that the latter problem is also undecidable.

To this end, we shall associate a logic with each Turing machine as follows:

\begin{law}
Let ${\bf M} = \langle P, Q, q_{0}, \delta \rangle$ be a Turing machine.
\benormal
\item Let $\LL({\bf M})$ be the algebraic language, whose set of constant symbols is
\[
P \cup Q \cup \{ 0, 1, \emptyset \},
\]
and which comprises a binary connective $x \cdot y$, and a ternary connective $\lambda(x, y, z)$.
\item Let $\vdash_{{\bf M}}$ be the logic in $\LL({\bf M})$ axiomatized by the following Hilbert-style calculus:
\begin{align*}
q \cdot \lambda(x \cdot y,  a, z) &\rhd q' \cdot \lambda(x,  y, b \cdot z) \tag{R1}\label{r1}\\
\hat{q} \cdot \lambda(x, \hat{a}, y \cdot z) &\rhd \hat{q}' \cdot \lambda( x \cdot \hat{b}, y, z)\tag{R2}\label{r2}\\
p \cdot \lambda(x, y, z) \lhd &\rhd p \cdot \lambda(\emptyset  \cdot x, y, z)\tag{R3}\label{r3}\\
p \cdot \lambda(x, y, z) \lhd &\rhd p \cdot \lambda(x, y, z \cdot \emptyset)\tag{R4}\label{r4}
\end{align*}
for every $p, q, q', \hat{q}, \hat{q}' \in P \cup Q$ and every $a, \hat{a}, b, \hat{b} \in \{ 0, 1, \emptyset \}$ such that 
\[
\delta(q, a) = \langle q', b, L\rangle \text{ and } \delta(\hat{q}, \hat{a}) = \langle \hat{q}', \hat{b}, R\rangle.
\]
\enormal
\end{law}

\noindent Let $\texttt{c} = \langle q, \vec{w}, v, \vec{u}\rangle$ be a configuration for a Turing machine ${\bf M}$, where
\[
\vec{w} = \langle w_{1}, \dots, w_{n}\rangle, v = \langle a \rangle, \text{ and }\vec{u} = \langle u_{1}, \dots, u_{m}\rangle.
\]
We associate a formula of $\mathscr{L}({\bf M})$ with $\texttt{c}$ as follows:
\[
\varphi_{\texttt{c}} \coloneqq q \cdot \lambda(   ( \cdots (( w_{1} \cdot w_{2}) \cdot w_{3})\cdots ) \cdot w_{n}, a, u_{1}\cdot ( u_{2} \cdot ( \dots (u_{m-1} \cdot  u_{m}) \dots ))).
\]
If $n =1$ (resp.\ $m =1$), the expression $( \cdots (( w_{1} \cdot w_{2}) \cdot w_{3})\cdots ) \cdot w_{n}$ (resp.\ $u_{1}\cdot ( u_{2} \cdot ( \dots (u_{m-1} \cdot  u_{m}) \dots ))$) in the above display should be understood as $w_{1}$ (resp.\ $u_{1}$).

%
%

\begin{Lemma}\label{Lem:configurations}
Let ${\bf M}$ be a Turing machine with configurations $\emph{\texttt{c}}$ and $\emph{\texttt{d}}$. If  $\emph{\texttt{c}}$ yields $\emph{\texttt{d}}$, then $\varphi_{\emph{\texttt{c}}} \vdash_{{\bf M}} \varphi_{\emph{\texttt{d}}}$.
\end{Lemma}

\begin{proof}
Consider two configurations $\texttt{c}$ and $\texttt{d}$ for ${\bf M}$ such that $\texttt{c}$ yields $\texttt{d}$. By definition of configuration, $\texttt{c}$ has the form $\langle q, \vec{w}, v, \vec{u}\rangle$ for some state $q$ and sequences
\[
\vec{w} = \langle w_{1}, \dots, w_{n}\rangle, v = \langle a \rangle, \text{ and }\vec{u} = \langle u_{1}, \dots, u_{m}\rangle.
\]
By symmetry, we can assume that $\delta(q, a) = \langle q', b, L\rangle$ for some state $q'$ and $b \in \{ 0, 1 \}$. Since $\texttt{c}$ yields $\texttt{d}$, the following holds:
\benroman
\item\label{item:config-1} If $n \geq 2$ and $u_{1} \ne \emptyset$, then $\texttt{d} = \langle q', \langle w_{1}, \dots, w_{n-1}\rangle, \langle w_{n}\rangle, \langle b, u_{1}, \dots, u_{m}\rangle\rangle$;
\item\label{item:config-2} If $n = 1$ and $u_{1} \ne \emptyset$, then $\texttt{d} = \langle q', \langle \emptyset\rangle, \langle w_{1}\rangle, \langle b, u_{1}, \dots, u_{m}\rangle\rangle$;
\item\label{item:config-3} If $n \geq 2$ and $u_{1} = \emptyset$, then $\texttt{d} = \langle q', \langle w_{1}, \dots, w_{n-1}\rangle, \langle w_{n}\rangle, \langle b\rangle\rangle$;
\item\label{item:config-4} If $n =1$ and $u_{1} = \emptyset$, then $\texttt{d} = \langle q', \langle \emptyset\rangle, \langle w_{1}\rangle, \langle b\rangle\rangle$.
\eroman

Accordingly, first suppose that $n \geq 2$ and $u_{1} \ne \emptyset$. By (\ref{item:config-1}) we can apply rule (\ref{r1}) to $\varphi_{\texttt{c}}$, obtaining $\varphi_{\texttt{c}} \vdash_{{\bf M}} \varphi_{\texttt{d}}$. 

Then we consider the case where $n = 1$ and $u_{1} \ne \emptyset$. In this case,
\[
\varphi_{\texttt{c}} = q \cdot \lambda(   w_{1}, a, u_{1}\cdot ( u_{2} \cdot ( \dots (u_{m-1} \cdot  u_{m}) \dots ))).
\]
Therefore, applying rule (\ref{r3}) to $\varphi_{\texttt{c}}$, we get
\[
\varphi_{\texttt{c}} \vdash_{{\bf M}}q \cdot \lambda(   \emptyset \cdot w_{1}, a, u_{1}\cdot ( u_{2} \cdot ( \dots (u_{m-1} \cdot  u_{m}) \dots ))).
\]
Furthermore, by (\ref{item:config-2}) we can apply rule (\ref{r1}) to the right-hand side of the above display, obtaining
\[
q \cdot \lambda(   \emptyset \cdot w_{1}, a, u_{1}\cdot ( u_{2} \cdot ( \dots (u_{m-1} \cdot  u_{m}) \dots ))) \vdash_{{\bf M}} \varphi_{\texttt{d}}.
\]
Consequently, $\varphi_{\texttt{c}} \vdash_{{\bf M}} \varphi_{\texttt{d}}$, as desired.

Consider then the case where $n \geq 2$ and $u_{1} = \emptyset$. Since $u_{1} = \emptyset$, the fact that $\texttt{c}$ is a configuration implies $m = 1$. Consequently,
\[
\varphi_{\texttt{c}} = q \cdot \lambda(  ( \cdots (( w_{1} \cdot w_{2}) \cdot w_{3}) \cdots ) \cdot w_{n}, a, \emptyset).
\]
By applying to the above formula the rule (\ref{r1}), we obtain
\[
\varphi_{\texttt{c}} \vdash_{{\bf M}} q' \cdot \lambda(  ( \cdots (( w_{1} \cdot w_{2}) \cdot w_{3})\cdots ) \cdot w_{n-1}, w_{n}, b \cdot \emptyset).
\]
By (\ref{item:config-3}), applying the rule (\ref{r4}) to the right-hand side of the above display, we get
\[
q' \cdot \lambda(  ( \cdots (( w_{1} \cdot w_{2}) \cdot w_{3}) \cdots ) \cdot w_{n-1}, w_{n}, b \cdot \emptyset) \vdash_{{\bf M}} \varphi_{\texttt{d}}.
\]
Consequently, $\varphi_{\texttt{c}} \vdash_{{\bf M}} \varphi_{\texttt{d}}$.

It only remains to consider the case where $n =1$ and $u_{1} = \emptyset$. As in the previous case, we get $m = 1$. Consequently,
\[
\varphi_{\texttt{c}} = q \cdot \lambda(w_{1}, a, \emptyset).
\]
By applying to the above formula the rule (\ref{r3}), we obtain
\[
\varphi_{\texttt{c}} \vdash_{{\bf M}}q \cdot \lambda(\emptyset \cdot w_{1}, a, \emptyset).
\]
Furthermore, applying the rule (\ref{r1}) to the right-hand side of the above display, we get
\[
q \cdot \lambda(\emptyset \cdot w_{1}, a, \emptyset) \vdash_{{\bf M}} q' \cdot  \lambda(  \emptyset, w_{1}, b \cdot \emptyset).
\]
By (\ref{item:config-4}), applying the rule (\ref{r4}) to the right-hand side of the above display, we have
\[
q' \cdot \lambda(  \emptyset, w_{1}, b \cdot \emptyset) \vdash_{{\bf M}}\varphi_{\texttt{d}}.
\]
Hence we conclude that, also in this case, $\varphi_{\texttt{c}} \vdash_{{\bf M}} \varphi_{\texttt{d}}$.
\end{proof}

We shall also associate a logic with every pair consisting of a Turing machine and an input.

\begin{law}
Let ${\bf M}$ be a Turing machine and $\vec{t}$ and input. 
\benormal
\item Let $\mathscr{L}({\bf M}, \vec{t})$ be $\mathscr{L}({\bf M})$ extended with a new binary connective $x \to y$.
\item Let also $\vdash_{{\bf M}}^{\vec{t}}$ be the logic axiomatized by the sets of rules (\ref{r1}), \dots, (\ref{r4}) plus
\begin{align*}
\emptyset &\rhd \varphi_{\text{In}({\bf M}, \vec{t}\?\? )}\tag{R5}\label{r5}\\
p \cdot y & \rhd x \to (x \cdot x)\tag{R6}\label{r6}\\
\emptyset & \rhd x \to x\tag{R7}\label{r7}\\
x, x \to y & \rhd y\tag{R8}\label{r8}\\
x_{1} \to y_{1}, \dots, x_{n} \to y_{n} & \rhd \ast(x_{1}, \dots, x_{n}) \to \ast(y_{1}, \dots, y_{n})\tag{R9}\label{r9}
\end{align*}
for every $p \in P$, every positive integer $n$, and every $n$-ary connective $\ast$.
\enormal
\end{law}

\begin{Lemma}\label{Lem:equivalential}
Let ${\bf M}$ be a Turing machine and $\vec{t}$ an input. For every pair of formulas $\epsilon, \delta \in Fm$,
\[
\emptyset \vdash_{{\bf M}}^{\vec{t}} \epsilon \to \delta \Longleftrightarrow \epsilon \text{ and }\delta \text{ are logically equivalent in }\vdash_{{\bf M}}^{\vec{t}}.
\]
\end{Lemma}

\begin{proof}
To prove the implication from right to left, observe that, taking $\varphi \coloneqq v \to \epsilon$, the assumption implies $\epsilon \to \epsilon \vdash_{{\bf M}}^{\vec{t}} \epsilon \to \delta$. By (\ref{r7}) this yields $\emptyset \vdash_{{\bf M}}^{\vec{t}} \epsilon \to \delta$, as desired.

Then we turn to prove the implication from left to right. First observe that from the rule (\ref{r9}) it follows
\[
x \to y, x \to x \vdash_{{\bf M}}^{\vec{t}} (x \to x) \to (y \to x).
\]
By rule (\ref{r7}) this simplifies to
\[
x \to y \vdash_{{\bf M}}^{\vec{t}} (x \to x) \to (y \to x).
\]
Finally, from (\ref{r8}) we get $x \to x, (x \to x) \to (y \to x) \vdash_{{\bf M}}^{\vec{t}} y \to x$, which simplifies to $(x \to x) \to (y \to x) \vdash_{{\bf M}}^{\vec{t}} y \to x$ by (\ref{r7}). Together with the above display, this implies
\begin{equation}\label{Eq:reverse:arrow}
x \to y \vdash_{{\bf M}}^{\vec{t}}  y \to x.
\end{equation}

Now, suppose that $\emptyset \vdash_{{\bf M}}^{\vec{t}} \epsilon \to \delta$ and consider a formula $\varphi(v, \vec{z})$. By applying repeatedly rule (\ref{r9}) and using the assumption $\emptyset \vdash_{{\bf M}}^{\vec{t}} \epsilon \to \delta$, it is not hard to see that $\emptyset \vdash_{{\bf M}}^{\vec{t}} \varphi(\epsilon, \vec{z}) \to \varphi(\delta, \vec{z})$.  By (\ref{Eq:reverse:arrow}) we get
\[
\emptyset \vdash_{{\bf M}}^{\vec{t}} \varphi(\epsilon, \vec{z}) \to \varphi(\delta, \vec{z}), \varphi(\delta, \vec{z}) \to \varphi(\epsilon, \vec{z}).
\]
Together with (\ref{r8}) this implies $\varphi(\epsilon, \vec{z})\sineq_{{\bf M}}^{\vec{t}} \varphi(\delta, \vec{z})$. It follows that $\epsilon$ and $\delta$ are logically equivalent in $\vdash_{{\bf M}}^{\vec{t}}$.
\end{proof}

\begin{Corollary}\label{Cor:undec-trick-different}
Let ${\bf M}$ be a Turing machine and $\vec{t}$ an input. The logic $\vdash_{{\bf M}}^{\vec{t}}$ has an algebraic semantics if and only if $\emptyset \vdash_{{\bf M}}^{\vec{t}} \epsilon \to \delta$ for some distinct $\epsilon,\delta \in Fm$.
\end{Corollary}

\begin{proof}
Suppose first that $\vdash_{{\bf M}}^{\vec{t}}$ is trivial. As $\vdash_{{\bf M}}^{\vec{t}}$ has theorems by rule (\ref{r7}), this means that it is inconsistent. Consequently, it has an algebraic semantics by Proposition \ref{Prop:trivial}. Moreover, $\emptyset \vdash_{{\bf M}}^{\vec{t}} \epsilon \to \delta$ for some (every) pair of distinct $\epsilon,\delta \in Fm$.

Then we consider the case where $\vdash_{{\bf M}}^{\vec{t}}$ is nontrivial. Observe that $\vdash_{{\bf M}}^{\vec{t}}$ is protoalgebraic by rules (\ref{r7}) and (\ref{r8}). Therefore, in view of Theorem \ref{Thm:protoalgebraic-logics}, $\vdash_{{\bf M}}^{\vec{t}}$ has an algebraic semantics if and only if there are distinct logically equivalent formulas $\epsilon$ and $\delta$. By Lemma \ref{Lem:equivalential} this happens precisely when $\emptyset \vdash_{{\bf M}}^{\vec{t}} \epsilon \to \delta$.
\end{proof}

The next result is the cornerstone of the argument.

\begin{Proposition}\label{Prop:main}
A Turing machine ${\bf M}$ halts on an input $\vec{t}$ if and only if the logic $\vdash_{\bf M}^{\vec{t}}$ has an algebraic semantics.
\end{Proposition}

\begin{proof}
First suppose that ${\bf M}$ halts on $\vec{t}$. Then there is a finite sequence of configurations $\texttt{c}_{1}, \dots, \texttt{c}_{n}$ such that $\texttt{c}_{1} = \textup{In}({\bf M}, \vec{t} \?\?)$, the state in $\texttt{c}_{n}$ is final, and $\texttt{c}_{i}$ yields $\texttt{c}_{i +1}$ for every positive integer $i < n$. By Lemma \ref{Lem:configurations} we get
\[
\varphi_{\texttt{c}_{1}} \vdash_{\bf M}^{\vec{t}} \varphi_{\texttt{c}_{2}} \vdash_{\bf M}^{\vec{t}} \dots \vdash_{\bf M}^{\vec{t}} \varphi_{\texttt{c}_{n}}.
\]
Together with (\ref{r5}) and $\texttt{c}_{1} = \textup{In}({\bf M}, \vec{t} \?\?)$, this implies $\emptyset \vdash_{\bf M}^{\vec{t}} \varphi_{\texttt{c}_{n}}$.
Thus, there is a formula $\psi$ and a final state $p \in P$ such that $\emptyset \vdash_{\bf M}^{\vec{t}} p \cdot \psi$. By rule (\ref{r6}) we obtain
\[
\emptyset \vdash_{\bf M}^{\vec{t}} x \to (x \cdot x).
\]
Hence, Corollary \ref{Cor:undec-trick-different} implies that $\vdash_{\bf M}^{\vec{t}}$ has an algebraic semantics.

To prove the ``if'' part, we reason by contraposition. Suppose that ${\bf M}$ does not halt on $\vec{t}$. Then there is an infinite sequence of configurations $\texttt{c}_{1}, \texttt{c}_{2}, \texttt{c}_{3}, \dots, \texttt{c}_{n}, \dots$ whose states are not final such that $\texttt{c}_{1} = \textup{In}({\bf M}, \vec{t}\?\?)$, and $\texttt{c}_{i}$ yields $\texttt{c}_{i+1}$ for every positive integer $i$. 

We shall associate a set of formulas $C_{n}$ to every $\texttt{c}_{n}$. To this end, consider a positive integer $n$ and recall that $\texttt{c}_{n}$ has the form $\langle q, \vec{w}, v, \vec{u}\rangle$ for some state $q$ and sequences
\[
\vec{w} = \langle w_{1}, \dots, w_{k}\rangle, v = \langle a \rangle, \text{ and }\vec{u} = \langle u_{1}, \dots, u_{m}\rangle.
\]
For every $i \in \omega$, we define two formulas
\begin{align*}
\alpha_{i} &\coloneqq \underbrace{\emptyset \cdot (\emptyset \cdot (\cdots (\emptyset}_{i\text{-times}} \cdot (w_{1} \cdot ( \cdots ( w_{k-1} \cdot w_{k}) \cdots ))) \cdots ))\\
\beta{i} &\coloneqq u_{1} \cdot (u_{2} \cdot (\cdots (u_{m} \cdot (\underbrace{\emptyset \cdot ( \cdots ( \emptyset \cdot \emptyset}_{i\text{-times}}) \cdots ))) \cdots )).
\end{align*}
Furthermore, let $[\alpha_{i}]$ and $[\beta_{i}]$ be the sets of formulas equivalent, respectively, to $\alpha_{i}$ and $\beta_{i}$ under the assumption that $\cdot$ is associative. In other words, $[\alpha_{i}]$ (resp.\ $[\beta{i}]$) is the set of formulas obtained reordering the parentheses in $\alpha_{i}$ (resp.\ $\beta_{i}$). Bearing this in mind, we set
\[
C_{n} \coloneqq \{ q \cdot \lambda(\epsilon, a, \delta ) \colon \text{there are }i, j \in \omega \text{ such that }\epsilon \in [\alpha_{i}] \text{ and }\delta \in [\beta_{i}] \}.
\]
Notice that $\varphi_{\texttt{c}_{n}} \in C_{n}$.

Now, consider the sets
\[
\Delta \coloneqq \{ q \cdot \lambda(\epsilon, \gamma, \delta) \colon q \in Q \text{ and } \epsilon, \gamma, \delta \in Fm \text{ are such that }\gamma \notin \{ 0, 1,  \emptyset \} \}
\]
and
\[
\Gamma \coloneqq \bigcup_{1 \leq n \in \omega} C_{n} \cup \Delta \cup \{ \gamma \to \gamma \colon \gamma \in Fm \}.
\]

\begin{Claim}\label{Claim:hopefully-last-itp}
The logic induced by $\langle \Fm, \Gamma \rangle$ extends $\vdash_{\bf M}^{\vec{t}}$. 
\end{Claim}

\begin{proof}
To prove the claim, it suffices to show that the rules axiomatizing $\vdash_{\bf M}^{\vec{t}}$ are valid in $\langle \Fm, \Gamma \rangle$, as we proceed to do.

(\ref{r1}-\ref{r2}): We detail only the case of (\ref{r1}) as the other one is analogous. Among the formulas in $\Gamma$, the rules (\ref{r1}) can be applied only to those in $\bigcup_{1 \leq n \in \omega}C_{n}$. Then suppose that a rule $q \cdot \lambda(x \cdot y, a, z) \rhd q' \cdot \lambda(x,  y, b \cdot z)$ in (\ref{r1}) is applied to a formula in $C_{n}$. In particular, this implies that $q$ is the state in $\texttt{c}_{n}$ and $q'$ the state in $\texttt{c}_{n+1}$. Consequently, $q' \in Q$, i.e., $q$ is not a final state. Moreover, observe the result of the application of the rule is a formula of the form $q' \cdot \lambda(\epsilon,  \gamma, b \cdot \delta)$ where $q \cdot \lambda(\epsilon \cdot \gamma, a, \delta) \in C_{n}$. Since $q'$ is not a final state, if $\gamma \notin \{ 0, 1, \emptyset \}$, then $q' \cdot \lambda(\epsilon,  \gamma, b \cdot \delta) \in \Delta \subseteq \Gamma$, and we are done. Then suppose that $\gamma \in \{ 0, 1, \emptyset \}$. But in this case, $q' \cdot \lambda(\epsilon,  \gamma, b \cdot \delta) \in C_{n+1} \subseteq \Gamma$, since $q' \cdot \lambda(\epsilon \cdot \gamma, a, \delta) \in C_{n}$. Thus, we conclude that the rules (\ref{r1}) are valid in $\langle \Fm, \Gamma \rangle$.

(\ref{r3}-\ref{r4}): Among formulas in $\Gamma$, these rules can be applied only to formulas in $\bigcup_{n \in \omega}C_{n} \cup \Delta$. But this set is closed under the application of the rules by construction.

(\ref{r5}): Recall that $\varphi_{\textup{In}({\bf M}, \vec{t} \?\? )} = \texttt{c}_{1} \in C_{1} \subseteq \Gamma$. Thus, the rule (\ref{r5}) is valid in $\langle \Fm, \Gamma \rangle$.

(\ref{r6}): Recall that the states in $\texttt{c}_{1}, \texttt{c}_{2}, \texttt{c}_{3}, \dots, \texttt{c}_{n}, \dots$ are not final. Therefore, this rule cannot be applied to formulas in $\Gamma$.

(\ref{r7}): This rule is valid in $\langle \Fm, \Gamma \rangle$ because $\{ \gamma \to \gamma \colon \gamma \in Fm \} \subseteq \Gamma$. 

(\ref{r8}-\ref{r9}): Among formulas in $\Gamma$, these rules can be applied only to formulas in $\{ \gamma \to \gamma \colon \gamma \in Fm \}$. This makes the validity of both rules in $\langle \Fm, \Gamma \rangle$ straightforward.
\end{proof} 

By Claim \ref{Claim:hopefully-last-itp}, if a formula is a theorem of $\vdash_{\bf M}^{\vec{t}}$, then it must belong to $\Gamma$. In particular, by definition of $\Gamma$, formulas of the form $\epsilon \to \delta$ can be theorems of $\vdash_{\bf M}^{\vec{t}}$ only if $\epsilon = \delta$. Hence, by Corollary \ref{Cor:undec-trick-different} we conclude that $\vdash_{\bf M}^{\vec{t}}$ does not have an algebraic semantics.
\end{proof}

\begin{proof}[Proof of Theorem \ref{Thm:comput-thm}(\ref{item:comput-thm-3})]
Observe that for every Turing machine ${\bf M}$ and input $\vec{t}$, the logic $\vdash_{\bf M}^{\vec{t}}$ is presented by a finite set of finite rules in a finite language. Thus, the result follows from the undecidability of the halting problem and Proposition \ref{Prop:main}.
\end{proof}

\begin{Remark}
The reader familiar with abstract algebraic logic might have noticed that the logic $\vdash_{\bf M}^{\vec{t}}$ is finitely equivalential \cite{Cz01,AAL-AIT-f}. Consequently, Theorem \ref{Thm:comput-thm}(\ref{item:comput-thm-3}) remains true when restricted to the family of finitely equivalential logics.
\qed
\end{Remark}

\paragraph{\bfseries Acknowledgements.}
Thanks are due to R. Jansana and M. Stronkowski for inspiring conversation on this topic and for many useful suggestions which helped to improve the presentation of the paper. Moreover, we gratefully thank the anonymous referee for the careful reading of the manuscript and I. Shillito for raising a question which inspired Corollary \ref{Cor:modal-not-an-algebraic-semantics}. The author was supported by the research grant $2017$ SGR $95$ of the AGAUR from the Generalitat de Catalunya, 
by the I+D+i research project PID2019-110843GA-I00 \textit{La geometria de las logicas no-clasicas} funded by the Ministry of Science and Innovation of Spain, and by the \textit{Beatriz Galindo} grant BEAGAL\-$18$/$00040$ funded by the Ministry of Science and Innovation of Spain.

\section{Appendix}

Given a graph-based language $\LL$, we denote by $At(\LL)$ the set of its atomic formulas, i.e., the union of $\Var$ and the set of constants of $\LL$. 

\begin{Lemma}\label{Lem:graph-based-combinatorics1}
Let $\LL$ be a graph-based language, $\Gamma \subseteq Fm_{\LL}$, and 
\[
\btau(x) = \{ \Box^{k}x \thickapprox \Box^{n} \textbf{c}_{i} \}
\]
for some nonnegative integers $n < k$ and ordinal $i < \alpha$. Moreover, let $p \in At(\LL) \smallsetminus \{ \textbf{c}_{i} \}$ and $h \in \omega$. Then $\langle \Box^{k}\Box^{h}p, \Box^{n}\textbf{c}_{i}\rangle \in \theta(\Gamma, \btau)$ if and only if there are $s, s^{\ast} \in \omega$, $\{ q_{j} \colon s > j \in \omega \} \subseteq At(\LL)$, and $\{ u_{j} \colon s > j \in \omega \} \cup \{ v_{j} \colon s > j \in \omega \} \cup \{ w_{j^{\ast}} \colon s^{\ast} > j^{\ast} \in \omega \} \subseteq \omega$ such that 
\benroman
\item\label{Eq:graph-combinatorics-1} $u_{j} < v_{j}$, for all $s > j \in \omega$;
\item\label{Eq:graph-combinatorics-2} $\Box^{u_{j}}q_{j}, \Box^{v_{j}}q_{j}, \Box^{w_{j^{\ast}}}\textbf{c}_{i} \in \Gamma$, for all $s > j \in \omega$ and $s^{\ast} > j^{\ast} \in \omega$;
\item\label{Eq:graph-combinatorics-3} $h = t+g$ for some $t, g \in \omega$ such that $\Box^{t}p \in \Gamma$ and $g$ is a multiple of
\[
\textup{gcd}(\{ v_{j} - u_{j} \colon s > j \in \omega \} \cup \{ k + w_{j^{\ast}} - n \colon s^{\ast} > j^{\ast} \in \omega \}).
\]
\eroman
\end{Lemma}

\begin{proof}
We begin by proving the ``only if'' part. To this end, suppose that $\langle \Box^{k}\Box^{h}p, \Box^{n}\textbf{c}_{i}\rangle \in \theta(\Gamma, \btau)$. There are two cases: either $\Box^{h}p \in \Gamma$ or $\Box^{h}p \notin \Gamma$. In the case where $\Box^{h}p \in \Gamma$, we define $s \coloneqq 0$ and $s^{\ast} \coloneqq 0$. Consequently, conditions (\ref{Eq:graph-combinatorics-1}, \ref{Eq:graph-combinatorics-2}) are vacuously satisfied. By the same token,
\[
d \coloneqq \textup{gcd}(\{ v_{j} - u_{j} \colon s > j \in \omega \} \cup \{ k + w_{j^{\ast}} - n \colon s^{\ast} > j^{\ast} \in \omega \}) = \textup{gdc}(\emptyset) = 0.
\]
Let also $g \coloneqq 0$ and $t \coloneqq h$. Clearly $h = t+g$ and $\Box^{t}p = \Box^{h} \in \Gamma$. Moreover, by the above display, $g$ is a multiple of $d$. This establishes condition (\ref{Eq:graph-combinatorics-3}), thus concluding the proof.

Then we consider the case where $\Box^{h}p \notin \Gamma$. By Theorem \ref{Thm:Maltsev-lemma} there are $\alpha_{0}, \dots, \alpha_{r} \in Fm$, $\gamma_{0}, \dots, \gamma_{r-1} \in \Gamma$, and $\delta_{0}(x, \vec{z}), \dots, \delta_{r-1}(x, \vec{z}) \in Fm$ such that
\[
\Box^{k}\Box^{h}p = \alpha_{0}, \Box^{n}\textbf{c}_{i} = \alpha_{r-1} \text{, and }\{ \alpha_{j}, \alpha_{j+1} \} = \{ \delta_{j}(\Box^{k}\gamma_{j}, \vec{z}), \delta_{j}(\Box^{n}\textbf{c}_{i}, \vec{z}) \}\text{, for every }r > j \in \omega.
\]
This implies 
\begin{equation}\label{Eq:infity-if-you-say-so}
\langle \Box^{k}\Box^{h}p, \Box^{n}\textbf{c}_{i}\rangle \in \theta(\Gamma^{-}, \btau),
\end{equation}
where
\[
\Gamma^{-} \coloneqq \{ \gamma_{0}, \dots, \gamma_{r-1} \}.
\]

Observe that we can assume without loss of generality that $x \in \Var(\delta_{0})$. As $\LL$ is graph-based, this implies the existence of a $g^{\ast} \in \omega$ such that $\delta_{0}(x, \vec{z}) = \Box^{g^{\ast}}x$. Bearing this in mind, we obtain
\[
\{\Box^{k}\Box^{h}p, \alpha_{1}\} = \{ \alpha_{0}, \alpha_{1} \} = \{ \delta_{0}(\Box^{k}\gamma_{0}, \vec{z}), \delta_{0}(\Box^{n}\textbf{c}_{i}, \vec{z}) \} = \{ \Box^{g^{\ast}}\Box^{k}\gamma_{0}, \Box^{g^{\ast}}\Box^{n} \textbf{c}_{i} \}.
\]
As $p \ne \textbf{c}_{i}$ by assumption, necessarily $\Box^{k}\Box^{h}p = \Box^{g^{\ast}}\Box^{k}\gamma_{0}$. Consequently, there is a $t^{\ast} \in \omega$ such that $\Box^{t^{\ast}}p = \gamma_{0} \in \Gamma^{-}$. Then let $t \coloneqq \min\{ s \in \omega \colon \Box^{s}p \in \Gamma^{-} \}$.  Clearly $t \leq t^{\ast}$. Consequently, there is a $g \in \omega$ such that
\[
\Box^{k}\Box^{h}p = \Box^{g^{\ast}}\Box^{k}\gamma_{0} = \Box^{g^{\ast}} \Box^{k} \Box^{t^{\ast}}p = \Box^{g} \Box^{k} \Box^{t}p.
\]
We conclude
\begin{equation}\label{Eq:to-prove-point-3}
h = t + g \text{ and }\Box^{t}p \in \Gamma^{-}.
\end{equation}

Now, consider the set
\begin{align*}
\Delta \coloneqq \{ \langle \gamma, \delta \rangle \colon & \gamma, \delta \in \Gamma^{-} \text{ and }\gamma = \Box^{u}q \text{ and } \delta = \Box^{v}q,\\
& \text{for some }q \in At(\LL) \text{ and nonnegative integers }u < v \}.
\end{align*}
Then let $s \in \omega$, $\{ q_{j} \colon s > j \in \omega \} \subseteq At(\LL)$, and $\{ u_{j} \colon s > j \in \omega \} \cup \{ v_{j} \colon s > j \in \omega \}\subseteq \omega$ be such that
\[
\Delta = \{ \langle \Box^{u_{j}}q_{j}, \Box^{v_{j}}q_{j}\rangle \colon s > j \in \omega \}.
\]
Similarly, consider
\[
\nabla \coloneqq \{ \Box^{w} \textbf{c}_{i} \in \Gamma^{-} \colon w \in \omega \}
\]
and let $s^{\ast} \in \omega$ and $\{ w_{j^{\ast}} \colon s^{\ast} > j^{\ast} \in \omega \} \subseteq \omega$ be such that
\[
\nabla = \{ \Box^{w_{j^{\ast}}}\textbf{c}_{i} \colon s^{\ast} > j^{\ast} \in \omega \}.
\]
Observe that conditions (\ref{Eq:graph-combinatorics-1}, \ref{Eq:graph-combinatorics-2}) hold for these definitions of $u_{j}, v_{j}, w_{j^{\ast}}$, and $q_{j}$.

Thus, to conclude the proof, it suffices to establish condition (\ref{Eq:graph-combinatorics-3}). By (\ref{Eq:to-prove-point-3}) it will be enough to show that 
\[
d \coloneqq \textup{gcd}(\{ v_{j} - u_{j} \colon s > j \in \omega \} \cup \{ k + w_{j^{\ast}} - n \colon s^{\ast} > j^{\ast} \in \omega \})
\]
is a divisor of $g$. To this end, with each $q \in At$ we associate a set
\[
M_{q} \coloneqq \{ m \in \omega \colon \Box^{m}q \in \Gamma^{-} \}.
\]
If $M_{q} \ne \emptyset$, we denote by $m_{q}$ its minimum element. Otherwise, we set $m_{q} \coloneqq \infty$. Observe that
\begin{equation}\label{Eq:final-infity}
t = m_{p}.
\end{equation}

\begin{Claim}\label{Claim:gcd-is-positive}
$d \geq 1$.
\end{Claim}

\begin{proof}
Suppose the contrary, with a view to contradiction. Then $d = 0$. By definition of $\textup{gcd}$, this implies $\Delta = \nabla = \emptyset$. As $\Delta = \emptyset$, for every $q \in At(\LL)$,
\begin{equation}\label{Eq:infity-1}
\text{either }M_{q} = \{ m_{q} \} \text{ or }m_{q} = \infty.
\end{equation}
Moreover, as $\nabla = \emptyset$,
\begin{equation}\label{Eq:infity-2}
m_{\textbf{c}_{i}} = \infty.
\end{equation}

Consider the sets
\begin{align*}
A_{1} &\coloneqq \{ \Box^{r}q \colon q \in At(\LL), r \in \omega, \text{ and }m_{q} = \infty \}\\
A_{2} &\coloneqq \{ \Box^{r}q \colon q \in At(\LL) \smallsetminus \{ \textbf{c}_{i} \}, m_{q} \in \omega, \text{ and }k + m_{q} -1 \geq r \in \omega \}.
\end{align*}
Moreover, define $A \coloneqq A_{1} \cup A_{2}$. Observe that
\begin{equation}\label{Eq:infity-3}
At(\LL) \subseteq A.
\end{equation}
To prove this, consider $q \in At(\LL)$. If $m_{q} = \infty$, then $q \in A_{1} \subseteq A$. Then suppose that $m_{q} \ne \infty$. By (\ref{Eq:infity-2}) $q \ne \textbf{c}_{i}$, whence $q \in At(\LL) \smallsetminus \{ \textbf{c}_{i} \}$. Moreover, $k \geq 1$ as by assumption $k > n \geq 0$. Thus $0 \leq k + m_{q} -1$. Consequently, $q = \Box^{0}q \in A_{2} \subseteq A$. This establishes (\ref{Eq:infity-3}).

We shall endow the set $A$ with the structure of an algebra $\A$ of type $\LL$. For all ordinal $j < \alpha$ we define $\textbf{c}_{j}^{\A} \coloneqq \textbf{c}_{j} \in A$. This can be done by (\ref{Eq:infity-3}). Moreover, we define $\Box^{\A} \colon A \to A$ as
\begin{displaymath}
\Box^{\A}\gamma \coloneqq \left\{\begin{array}{@{\,}ll}
\Box \gamma & \text{if $\Box \gamma \in A$}\\
\Box^{n} \textbf{c}_{i} & \text{otherwise}\\
\end{array} \right.
\end{displaymath}
for every $\gamma \in A$.
Observe that $\Box^{\A}$ is well-defined by (\ref{Eq:infity-2}). Moreover, from (\ref{Eq:infity-2}) it follows
\begin{equation}\label{Eq:infity-3/4}
\underbrace{\Box^{\A} \dots \Box^{\A}}_{r\text{-times}}\textbf{c}_{i} = \Box^{r}\textbf{c}_{i}, \text{ for all }r \in \omega.
\end{equation}
Similarly, from from (\ref{Eq:infity-1}) it follows that for all $q \in At(\LL)$,
\begin{equation}\label{Eq:infity-3/4b}
\text{if }m_{q} \in \omega\text{, then }\underbrace{\Box^{\A} \dots \Box^{\A}}_{(k + m_{q})\text{-times}}q = \Box^{n}\textbf{c}_{i}.
\end{equation}
Now, recall by (\ref{Eq:infity-3}) that $\Var \subseteq A$. Bearing this in mind, let $h \colon \Fm \to \A$ be the unique homomorphism such that $v(x) = x$ for every $x \in \Var$. We shall prove that
\begin{equation}\label{Eq:infity-4}
h(\Box^{k}\gamma) = h(\Box^{n} \textbf{c}_{i})
\end{equation}
for every $\gamma \in \Gamma^{-}$. To this end, consider $\gamma \in \Gamma^{-}$. There are $q \in At(\LL)$ and $j \in \omega$ such that $\gamma = \Box^{j}q$. Together with $\Box^{j}q \in \Gamma^{-}$, (\ref{Eq:infity-1}, \ref{Eq:infity-2})  imply $m_{q} = j \in \omega$. Hence, applying (\ref{Eq:infity-3/4}, \ref{Eq:infity-3/4b}), we get
\[
h(\Box^{k} \gamma)= h(\Box^{k}\Box^{m_{q}}q) =\underbrace{\Box^{\A} \dots \Box^{\A}}_{(k+m_{q})\text{-times}}h(q) = \underbrace{\Box^{\A} \dots \Box^{\A}}_{(k+m_{q})\text{-times}}q = \Box^{n} \textbf{c}_{i} = h(\Box^{n}\textbf{c}_{i}),
\]
establishing (\ref{Eq:infity-4}).

By (\ref{Eq:infity-4}), the kernel $\textup{Ker}(h)$ contains the generators of $\theta(\Gamma^{-}, \btau)$, whence $\theta(\Gamma^{-}, \btau) \subseteq \textup{Ker}(h)$. Together with (\ref{Eq:infity-if-you-say-so}), this yields $\langle \Box^{k}\Box^{h}p, \Box^{n}\textbf{c}_{i}\rangle \in \theta(\Gamma^{-}, \btau) \subseteq \textup{Ker}(h)$, i.e., $h(\Box^{k}\Box^{h}p) = h(\Box^{n}\textbf{c}_{i})$. Now, recall that $\Box^{t}p \in \Gamma^{-}$. By (\ref{Eq:infity-1}) this implies $m_{p} = t \in \omega$ and, therefore, $h = m^{p} + g$ (as $h = t + g$). By the definition of $\Box^{\A}$ and (\ref{Eq:infity-3/4}, \ref{Eq:infity-3/4b}) we obtain
\begin{align*}
h(\Box^{h}\Box^{k}p) &= h(\Box^{g + k + m_{p}}p) = \underbrace{\Box^{\A} \dots \Box^{\A}}_{g\text{-times}}\underbrace{\Box^{\A} \dots \Box^{\A}}_{(k + m_{p})\text{-times}}p\\
&=\underbrace{\Box^{\A} \dots \Box^{\A}}_{g\text{-times}}\Box^{n}\textbf{c}_{i} = \Box^{g+n}\textbf{c}_{i}.
\end{align*}
Together with $h(\Box^{k}\Box^{h}p) = h(\Box^{n}\textbf{c}_{i})$ and (\ref{Eq:infity-3/4}), this implies $\Box^{n}\textbf{c}_{i} = \Box^{g+n}\textbf{c}_{i}$, whence $g = 0$. Thus $\Box^{h}p = \Box^{t+g}p = \Box^{t}p \in \Gamma^{-} \subseteq \Gamma$. But this contradicts the assumption that $\Box^{h}p \notin \Gamma$. Hence we reached a contradiction, as desired.
\end{proof}

Now, consider the sets
\begin{align*}
B_{1} &\coloneqq \{ \Box^{r}q \colon q \in At(\LL) \smallsetminus \{ \textbf{c}_{i} \}, r \in \omega, \text{ and }m_{q} = \infty \}\\
B_{2} &\coloneqq \{ \Box^{r}q \colon q \in At(\LL) \smallsetminus \{ \textbf{c}_{i} \}, m_{q} \in \omega, \text{ and } k + m_{q} -1 \geq r \in \omega \}\\
B_{3} &\coloneqq \{ \Box^{r}\textbf{c}_{i} \colon n+d-1 \geq r \in \omega \}.
\end{align*}
Moreover, define $B \coloneqq B_{1} \cup B_{2} \cup B_{3}$. Observe that
\begin{equation}\label{Eq:infity-5}
At(\LL) \subseteq B.
\end{equation}
To prove this, consider $q \in At(\LL)$. We have two cases: either $q \ne \textbf{c}_{i}$ or $q = \textbf{c}_{i}$. If $q \ne \textbf{c}_{i}$, we repeat the argument used to establish (\ref{Eq:infity-3}). Then we consider the case where $q = \textbf{c}_{i}$. By Claim \ref{Claim:gcd-is-positive}  we get $0 \leq n + d -1$, whence $\textbf{c}_{i} \in B_{3} \subseteq B$. This establishes (\ref{Eq:infity-5}).

Let $\B$ be the algebra of type $\LL$ obtained endowing $B$ with the interpretation of the symbols in $\LL$ defined in the proof of Claim \ref{Claim:gcd-is-positive} for the algebra $\A$. Observe that the interpretation of constant symbols is well-defined by  (\ref{Eq:infity-5}). Similarly, from Claim \ref{Claim:gcd-is-positive} it follows that $n \leq n + d - 1$ and, therefore, that $\Box^{n} \textbf{c}_{i} \in B_{3} \subseteq B$. This guarantees that the interpretation of $\Box$ in $\B$ is also well-defined.

Observe that for all $r \in \omega$,
\begin{equation}
\underbrace{\Box^{\B} \dots \Box^{\B}}_{r\text{-times}}\textbf{c}_{i} = \begin{cases}
\Box^{r}\textbf{c}_{i} & \text{if $r \leq n$}\\
\Box^{n + a} \textbf{c}_{i} & \text{if $r > n$ and $a$ is the remainder of $\frac{r - n}{d}$.}\\
  \end{cases} \label{eq:torus1}
\end{equation}
Similarly, for every $q \in At(\LL) \smallsetminus \{ \textbf{c}_{i} \}$ such that $m_{q} \in \omega$, and every $r \in \omega$,
\begin{equation}
\underbrace{\Box^{\B} \dots \Box^{\B}}_{r\text{-times}}q = \begin{cases}
\Box^{r}q & \text{if $r \leq k + m_{q} - 1$}\\
\Box^{n + a} \textbf{c}_{i} & \text{if $r > k + m_{q}$ and $a$ is the remainder of $\frac{r - k - m_{q}}{d}$.}\\
\end{cases} \label{eq:torus2}
\end{equation}

Now, recall by (\ref{Eq:infity-5}) that $\Var \subseteq B$. Bearing this in mind, let $h \colon \Fm \to \B$ be the unique homomorphism such that $v(x) = x$ for every $x \in \Var$. We shall prove that
\begin{equation}\label{Eq:infity-6}
h(\Box^{k}\gamma) = h(\Box^{n} \textbf{c}_{i})
\end{equation}
for every $\gamma \in \Gamma^{-}$. To this end, consider $\gamma \in \Gamma^{-}$. There are $q \in At(\LL)$ and $j \in \omega$ such that $\gamma = \Box^{j}q$. There are two cases: either $q \ne \textbf{c}_{i}$ or $q = \textbf{c}_{i}$. 

First suppose that $q \ne \textbf{c}_{i}$. As $\Box^{j}q = \gamma \in \Gamma$, necessarily $m_{q} \leq j$. Observe that
\begin{equation}\label{Eq:infity7}
j - m_{q} = r \cdot d\text{, for some }r \in \omega.
\end{equation}
To prove this, notice that if $j = m_{q}$, we are done taking $r \coloneqq 0$. Then suppose $j \ne m_{q}$, i.e., $m_{q} < j$. In this case, $\langle \Box^{m_{q}}q, \Box^{j}q\rangle \in \Delta$. Consequently, $j - m_{q}$ is a multiple of $d$, as desired. This establishes (\ref{Eq:infity7}). Together with (\ref{eq:torus1}, \ref{eq:torus2}), this yields
\begin{align*}
h(\Box^{k}\gamma) &= h(\Box^{k}\Box^{j}q) = h(\Box^{j-m_{q}}\Box^{k + m_{q}}q) = h(\Box^{r \cdot d}\Box^{k + m_{q}}q)\\
& = \underbrace{\Box^{\B} \dots \Box^{\B}}_{(r \cdot d)\text{-times}}\underbrace{\Box^{\B} \dots \Box^{\B}}_{(k + m_{q})\text{-times}}h(q) = \underbrace{\Box^{\B} \dots \Box^{\B}}_{(r \cdot d)\text{-times}}\underbrace{\Box^{\B} \dots \Box^{\B}}_{(k + m_{q})\text{-times}}q\\
&= \Box^{n} \textbf{c}_{i} = h(\Box^{n} \textbf{c}_{i}).
\end{align*}

Then we consider the case where $q = \textbf{c}_{i}$. In this case, an argument similar to the one detailed above (where $\nabla$ takes the role of $\Delta$) shows that there is an $r \in \omega$ such that $r \cdot d = k + j - n$. Consequently, by (\ref{eq:torus1}), we obtain
\begin{align*}
h(\Box^{k}\gamma) &= h(\Box^{k}\Box^{j}\textbf{c}_{i}) = h(\Box^{k+j-n}\Box^{n}\textbf{c}_{i}) = h(\Box^{r \cdot d}\Box^{n}\textbf{c}_{i})\\
&=\underbrace{\Box^{\B} \dots \Box^{\B}}_{(r \cdot d)\text{-times}}\underbrace{\Box^{\B} \dots \Box^{\B}}_{n\text{-times}}h(\textbf{c}_{i}) = \underbrace{\Box^{\B} \dots \Box^{\B}}_{(r \cdot d)\text{-times}}\underbrace{\Box^{\B} \dots \Box^{\B}}_{n\text{-times}}\textbf{c}_{i}\\
&= \Box^{n} \textbf{c}_{i} = h(\Box^{n} \textbf{c}_{i}),
\end{align*}
thus establishing (\ref{Eq:infity-6}).

Lastly, by (\ref{Eq:infity-6}) the kernel $\textup{Ker}(h)$ contains the generators of $\theta(\Gamma^{-}, \btau)$, whence $\theta(\Gamma^{-}, \btau) \subseteq \textup{Ker}(h)$. By (\ref{Eq:infity-if-you-say-so}) this yields $\langle \Box^{k}\Box^{h}p, \Box^{n}\textbf{c}_{i}\rangle \in \theta(\Gamma^{-}, \btau) \subseteq \textup{Ker}(h)$, i.e., $h(\Box^{k}\Box^{h}p) = h(\Box^{n}\textbf{c}_{i})$. We obtain
\[
\underbrace{\Box^{\B} \dots \Box^{\B}}_{(g + m_{p} + k)\text{-times}}p = \underbrace{\Box^{\B} \dots \Box^{\B}}_{(g + t + k)\text{-times}}p = \underbrace{\Box^{\B} \dots \Box^{\B}}_{(h + k)\text{-times}}p = h(\Box^{k}\Box^{h}p) = h(\Box^{n}\textbf{c}_{i})= \Box^{n} \textbf{c}_{i}.
\]
The equalities above are justified as follows. First first one follows from (\ref{Eq:final-infity}), the second from (\ref{Eq:to-prove-point-3}), the third one is obvious, the fourth one was justified right before the above display, and the last one follows from (\ref{eq:torus1}).

Now, as $m_{p} = t \in \omega$ and $p \ne \textbf{c}_{i}$, the value of $\underbrace{\Box^{\B} \dots \Box^{\B}}_{(g + m_{p} + k)\text{-times}}p$ can be computed according to condition (\ref{eq:torus2}). Together with the above display, this guarantees that $g$ is a multiple of $d$. Hence, condition (\ref{Eq:graph-combinatorics-3}) holds, as desired.

To prove the ``if'' part, suppose that conditions (\ref{Eq:graph-combinatorics-1}, \ref{Eq:graph-combinatorics-2}, \ref{Eq:graph-combinatorics-3}) hold. Consider $s > j \in \omega$. Observe that $\Box^{u_{j}}q_{j}, \Box^{v_{j}}q_{j} \in \Gamma$. Consequently,
\[
\langle \Box^{k}\Box^{u_{j}}q_{j}, \Box^{n} \textbf{c}_{i}\rangle, \langle \Box^{k}\Box^{v_{j}}q_{j}, \Box^{n} \textbf{c}_{i}\rangle \in \theta(\Gamma, \btau).
\]
Bearing in mind that $u_{j} < v_{j}$, this easily implies
\begin{equation}\label{Eq:gcd-infity1}
\Box^{v_{j}- u_{j}}\Box^{n} \textbf{c}_{i} \equiv \Box^{v_{j}-u_{j}}\Box^{k}\Box^{u_{j}}q_{j} = \Box^{v_{j}}\Box^{k}q_{j} \equiv \Box^{n} \textbf{c}_{i}\mod{\theta(\Gamma, \btau)}.
\end{equation}
Furthermore, consider any $s^{\ast} > j^{\ast} \in \omega$ and recall that $\Box^{w_{j^{\ast}}} \textbf{c}_{i} \in \Gamma$. Consequently,
\[
\langle \Box^{k} \Box^{w_{j^{\ast}}} \textbf{c}_{i}, \Box^{n} \textbf{c}_{i}\rangle \in \theta(\Gamma, \btau).
\]
Again, this implies
\begin{equation}\label{Eq:gcd-infity2}
\Box^{k + w_{j^{\ast}} - n} \Box^{n} \textbf{c}_{i} = \Box^{k + w_{j^{\ast}}} \textbf{c}_{i} \equiv \Box^{n} \textbf{c}_{i}\mod{\theta(\Gamma, \btau)}.
\end{equation}

Now, from (\ref{Eq:gcd-infity1}, \ref{Eq:gcd-infity2}) it follows that for all $r \in \omega$,
\begin{equation}\label{Eq:gcd-infity-final}
\langle\Box^{n} \textbf{c}_{i}, \Box^{r \cdot d}\Box^{n} \textbf{c}_{i}\rangle \in \theta(\Gamma, \btau),
\end{equation}
where
\[
d \coloneqq \textup{gcd}(\{ v_{j} - u_{j} \colon s > j \in \omega \} \cup \{ k + w_{j^{\ast}} - n \colon s^{\ast} > j^{\ast} \in \omega \}).
\]
Hence there is an $r \in \omega$ such that
\[
\Box^{k}\Box^{h}p = \Box^{g} \Box^{t} \Box^{k}p \equiv \Box^{g} \Box^{n} \textbf{c}_{i} = \Box^{r \cdot d} \Box^{n} \textbf{c}_{i} \equiv \Box^{n} \textbf{c}_{i}\mod{\theta(\Gamma, \btau)}.
\]
The steps in the above display are justified as follows. The first follows from the assumption that $h = t+g$. The second from the assumption that $\Box^{t}p \in \Gamma$ and, therefore, $\langle \Box^{k}\Box^{t}p, \Box^{n} \textbf{c}_{i}\rangle \in \theta(\Gamma, \btau)$. The third from the assumption that there is some $r \in \omega$ such that $g = r \cdot d$, and the last one from (\ref{Eq:gcd-infity-final}). The above display concludes the proof.
\end{proof}

\begin{Lemma}\label{Lem:graph-based-combinatorics2}
Let $\LL$ be a graph-based language, $\Gamma \subseteq Fm_{\LL}$, and 
\[
\btau(x) = \{ \Box^{k}x \thickapprox \Box^{n} \textbf{c}_{i} \}
\]
for some nonnegative integers $n < k$ and ordinal $i < \alpha$. Moreover, let $h \in \omega$. Then $\langle \Box^{k}\Box^{h} \textbf{c}_{i}, \Box^{n}\textbf{c}_{i}\rangle \in \theta(\Gamma, \btau)$ if and only if there are $s, s^{\ast} \in \omega$, $\{ q_{j} \colon s > j \in \omega \} \subseteq At(\LL)$, and $\{ u_{j} \colon s > j \in \omega \} \cup \{ v_{j} \colon s >  j \in \omega \} \cup \{ w_{j^{\ast}} \colon s^{\ast} > j^{\ast} \in \omega \} \subseteq \omega$ such that 
\benroman
\item\label{Eq:graph-combinatorics-1-c} $u_{j} < v_{j}$, for all $s > j \in \omega$;
\item\label{Eq:graph-combinatorics-2-c} $\Box^{u_{j}}q_{j}, \Box^{v_{j}}q_{j}, \Box^{w_{j^{\ast}}}\textbf{c}_{i} \in \Gamma$, for all $s > j \in \omega$ and $s^{\ast} > j^{\ast} \in \omega$;
\item\label{Eq:graph-combinatorics-3-c} $k + h - n$ is a multiple of
\[
\textup{gcd}(\{ v_{j} - u_{j} \colon s > j \in \omega \} \cup \{ k + w_{j^{\ast}} - n \colon s^{\ast} > j^{\ast} \in \omega \}).
\]
\eroman
\end{Lemma}

\begin{proof}
We begin by proving the ``only if'' part. To this end, suppose that $\langle \Box^{k}\Box^{h} \textbf{c}_{i}, \Box^{n}\textbf{c}_{i}\rangle \in \theta(\Gamma, \btau)$. We define $\Delta$, $\nabla$, and $m_{q}$ as in the proof of Lemma \ref{Lem:graph-based-combinatorics1}. Similarly, we enumerate $\Delta$ and $\nabla$ as in the proof of Lemma \ref{Lem:graph-based-combinatorics1}. Observe that conditions (\ref{Eq:graph-combinatorics-1-c}, \ref{Eq:graph-combinatorics-2-c}) hold for these definitions of $u_{j}, v_{j}, w_{j^{\ast}}$, and $q_{j}$.

Then define
\[
d \coloneqq \textup{gcd}(\{ v_{j} - u_{j} \colon s > j \in \omega \} \cup \{ k + w_{j^{\ast}} - n \colon s^{\ast} > j^{\ast} \in \omega \}).
\]
To conclude this part of the proof, it suffices to establish condition (\ref{Eq:graph-combinatorics-3-c}), i.e., to prove that $k + h - n$ is a multiple of $d$. To this end, we rely on the following:
\begin{Claim}\label{Claim:gcd-positive-c}
$d \geq 1$.
\end{Claim}
\begin{proof}
Suppose the contrary, with a view to contradiction. We replicate the first part of the proof of Claim \ref{Claim:gcd-is-positive} by defining an algebra $\A$ and a homomorphism $h \colon \Fm \to \A$ such that $\textup{Ker}(h)$ extends $\theta(\Gamma, \btau)$. Together with the assumption $\langle \Box^{k}\Box^{h} \textbf{c}_{i}, \Box^{n}\textbf{c}_{i}\rangle \in \theta(\Gamma, \btau)$, the fact that $\theta(\Gamma, \btau) \subseteq \textup{Ker}(h)$ implies
\[
\underbrace{\Box^{\A} \dots \Box^{\A}}_{(k+h)\text{-times}} \textbf{c}_{i} = h(\Box^{k}\Box^{h} \textbf{c}_{i}) = h(\Box^{n}\textbf{c}_{i}) = \underbrace{\Box^{\A} \dots \Box^{\A}}_{n\text{-times}} \textbf{c}_{i}.
\]
By (\ref{Eq:infity-3/4}) this yields $k + h = n$. But this contradicts the fact that $k > n$. Hence we reached a contradiction, as desired.
\end{proof}
Using Claim \ref{Claim:gcd-positive-c}, we replicate the construction of the algebra $\B$ and pf the homomorphism $h \colon \Fm \to \B$ in the proof of Lemma \ref{Lem:graph-based-combinatorics1}, and obtain $\textup{Ker}(h) \supseteq\theta(\Gamma, \btau)$. Together with the assumption $\langle \Box^{k}\Box^{h} \textbf{c}_{i}, \Box^{n}\textbf{c}_{i}\rangle \in \theta(\Gamma, \btau)$, this implies
\[
\underbrace{\Box^{\B} \dots \Box^{\B}}_{(k+h)\text{-times}} \textbf{c}_{i} = h(\Box^{k}\Box^{h} \textbf{c}_{i}) = h(\Box^{n}\textbf{c}_{i}) = \underbrace{\Box^{\B} \dots \Box^{\B}}_{n\text{-times}} \textbf{c}_{i}.
\]
As $k > n$, we get $h+k > n$. Thus, together with (\ref{eq:torus1}), the above display implies that $k+h-n$ is a multiple of $d$.

To prove the ``if'' part, suppose that conditions (\ref{Eq:graph-combinatorics-1-c}, \ref{Eq:graph-combinatorics-2-c}, \ref{Eq:graph-combinatorics-3-c}) hold. As in the proof of Lemma \ref{Lem:graph-based-combinatorics1}, we obtain that for all $r \in \omega$,
\begin{equation}\label{Eq:gcd-infity-final-c}
\langle\Box^{n} \textbf{c}_{i}, \Box^{r \cdot d}\Box^{n} \textbf{c}_{i}\rangle \in \theta(\Gamma, \btau),
\end{equation}
where
\[
d \coloneqq \textup{gcd}(\{ v_{j} - u_{j} \colon s > j \in \omega \} \cup \{ k + w_{j^{\ast}} - n \colon s^{\ast} > j^{\ast} \in \omega \}).
\]
Hence we obtain that for some $r \in \omega$,
\[
\Box^{k}\Box^{h}\textbf{c}_{i} = \Box^{k + h - n}\Box^{n}\textbf{c}_{i} = \Box^{r \cdot d} \Box^{n} \textbf{c}_{i} \equiv \Box^{n} \textbf{c}_{i} \mod \theta(\Gamma, \btau).
\]
The steps in the above display are justified as follows. The first is sound because by assumption $k > n$. The second follows from the assumption that there is some $r \in \omega$ such that $k + h - n = r \cdot d$, and the last one from (\ref{Eq:gcd-infity-final-c}). The above display concludes the proof.
\end{proof}

\bibliographystyle{plain}
\end{document}